\documentclass{amsart}
\usepackage{amsrefs}
\usepackage{a4wide}
\input xypic
\newdir{ >}{{}*!/-9pt/\dir{>}}

\renewcommand{\O}       {\mathcal{O}}
\newcommand{\OC}        {\mathcal{O}_C}
\newcommand{\MC}        {\mathcal{M}_C}
\newcommand{\OD}        {\mathcal{O}_D}
\newcommand{\OS}        {\mathcal{O}_S}
\newcommand{\OT}        {\mathcal{O}_T}
\newcommand{\OV}        {\mathcal{O}_V}
\newcommand{\OX}        {\mathcal{O}_X}
\newcommand{\OY}        {\mathcal{O}_Y}
\newcommand{\OZ}        {\mathcal{O}_Z}

\newcommand{\Div}       {\operatorname{Div}}
\newcommand{\End}       {\operatorname{End}}

\newcommand{\Hom}       {\operatorname{Hom}}
\newcommand{\Map}       {\operatorname{Map}}
\newcommand{\Nil}       {\operatorname{Nil}}
\newcommand{\Pic}       {\operatorname{Pic}}
\newcommand{\Sub}       {\operatorname{Sub}}
\newcommand{\Sym}       {\operatorname{Sym}}
\newcommand{\adj}       {\operatorname{adj}}
\newcommand{\ann}       {\operatorname{ann}}
\newcommand{\conn}      {\operatorname{conn}}
\newcommand{\graph}     {\operatorname{graph}}
\newcommand{\img}       {\operatorname{image}}
\newcommand{\ord}       {\operatorname{ord}}
\newcommand{\res}       {\operatorname{res}}
\newcommand{\spec}      {\operatorname{spec}}
\newcommand{\spf}       {\operatorname{spf}}
\newcommand{\trace}     {\operatorname{trace}}
\newcommand{\trf}	{\operatorname{tr}}

\newcommand{\al}        {\alpha}
\newcommand{\bt}        {\beta} 
\newcommand{\gm}        {\gamma}
\newcommand{\dl}        {\delta}
\newcommand{\ep}        {\epsilon}
\newcommand{\zt}        {\zeta}
\newcommand{\tht}       {\theta}
\newcommand{\sg}        {\sigma}
\newcommand{\lm}        {\lambda}
\newcommand{\om}        {\omega}

\newcommand{\Tht}       {\Theta}

\newcommand{\gmb}       {\overline{\gamma}}
\newcommand{\dlb}       {\overline{\delta}}
\newcommand{\phb}       {\overline{\phi}}
\newcommand{\lmb}       {\overline{\lambda}}

\newcommand{\Gm}        {\Gamma}
\newcommand{\Dl}        {\Delta}
\newcommand{\Sg}        {\Sigma}
\newcommand{\tSg}       {\widetilde{\Sigma}}
\newcommand{\Sgi}       {\Sigma^\infty}
\newcommand{\Sgip}      {\Sigma^\infty_+}
\newcommand{\Om}        {\Omega}
\newcommand{\bOm}       {\overline{\Omega}}

\newcommand{\F}         {{\mathbb{F}}}
\newcommand{\C}         {{\mathbb{C}}}
\newcommand{\Fp}        {{\mathbb{F}_p}}          
\newcommand{\N}         {{\mathbb{N}}}
\newcommand{\Q}         {{\mathbb{Q}}}
\newcommand{\Z}         {{\mathbb{Z}}}
\newcommand{\Zp}        {{\mathbb{Z}_p}}          
\newcommand{\Zpl}       {{\mathbb{Z}_{(p)}}}      

\newcommand{\CA}        {{\mathcal{A}}}
\newcommand{\CC}        {{\mathcal{C}}}
\newcommand{\CE}        {{\mathcal{E}}}
\newcommand{\CF}        {{\mathcal{F}}}
\newcommand{\CG}        {{\mathcal{G}}}
\newcommand{\CL}        {{\mathcal{L}}}
\newcommand{\CO}        {{\mathcal{O}}}
\newcommand{\PP}        {{\mathcal{P}}}
\newcommand{\CS}        {{\mathcal{S}}}
\newcommand{\CU}        {{\mathcal{U}}}
\newcommand{\CV}        {{\mathcal{V}}}

\newcommand{\Ab}        {\overline{A}}
\newcommand{\Bb}        {\overline{B}}
\newcommand{\Cb}        {\overline{C}}
\newcommand{\Db}        {\overline{D}}
\newcommand{\Ib}        {\overline{I}}
\newcommand{\Jb}        {\overline{J}}
\newcommand{\Rb}        {\overline{R}}
\newcommand{\Sb}        {\overline{S}}
\newcommand{\xb}        {\overline{x}}
\newcommand{\ub}        {\overline{u}}

\newcommand{\aff}       {\mathbb{A}}
\newcommand{\haf}       {\widehat{\mathbb{A}}}
\newcommand{\hot}       {\widehat{\otimes}}
\newcommand{\invlim} {\operatornamewithlimits{\underset{\longleftarrow}{lim}}}
\newcommand{\colim}  {\operatornamewithlimits{\underset{\longrightarrow}{lim}}}
\newcommand{\mxi}       {\mathfrak{m}}
\newcommand{\nxi}       {\mathfrak{n}}
\newcommand{\op}        {\oplus}
\newcommand{\ot}        {\otimes}
\newcommand{\psb}[1]    {[\![#1]\!]}
\newcommand{\sse}       {\subseteq}
\newcommand{\st}        {\;|\;}
\newcommand{\tm}        {\times}
\newcommand{\Smash}     {\wedge}
\newcommand{\hC}        {\widehat{C}}
\newcommand{\hE}        {\widehat{E}}
\newcommand{\hGa}       {\widehat{G}_a}
\newcommand{\hK}        {\widehat{K}}
\newcommand{\CP}        {\mathbb{C}P}
\newcommand{\CPi}       {\mathbb{C}P^\infty}
\newcommand{\MG}        {\mathbb{G}_m}
\newcommand{\sm}        {\setminus}
\newcommand{\tC}        {\widetilde{C}}
\newcommand{\tD}        {\widetilde{D}}
\newcommand{\tE}        {\widetilde{E}}

\newcommand{\tI}        {\widetilde{I}}
\newcommand{\tJ}        {\widetilde{J}}
\newcommand{\tph}       {\tilde{\phi}}
\newcommand{\tbt}       {\tilde{\bt}}
\newcommand{\tmu}       {\tilde{\mu}}
\newcommand{\tc}        {\tilde{c}}
\newcommand{\tf}        {\tilde{f}}
\newcommand{\tg}        {\tilde{g}}
\renewcommand{\th}      {\tilde{h}}
\newcommand{\tj}        {\tilde{\jmath}}
\newcommand{\oB}        {B^\circ}
\newcommand{\bd}        {\mathbf{d}}
\newcommand{\bD}        {\mathbf{D}}
\newcommand{\tG}        {\widetilde{G}}
\newcommand{\tdim}      {\widetilde{\dim}}
\newcommand{\two}       {\mathbf{2}}
\newcommand{\BU}[1]     {BU\langle #1\rangle}
\newcommand{\tF}        {\widetilde{F}}
\newcommand{\tdet}      {\widetilde{\det}}

\newcommand{\un}[1]     {\underline{#1}}
\newcommand{\ov}[1]     {\overline{#1}}
\newcommand{\ip}[1]     {\langle #1\rangle}

\newcommand{\xra}       {\xrightarrow}

\newcommand{\bigWedge}  {\bigvee}
\newcommand{\bigSmash}  {\bigwedge}

\renewcommand{\:}{\colon}

\newtheorem{theorem}{Theorem}[section]

\newtheorem{lemma}[theorem]{Lemma}
\newtheorem{proposition}[theorem]{Proposition}
\newtheorem{corollary}[theorem]{Corollary}
\theoremstyle{definition}
\newtheorem{remark}[theorem]{Remark}
\newtheorem{definition}[theorem]{Definition}
\newtheorem{construction}[theorem]{Construction}
\newtheorem{example}[theorem]{Example}
\newtheorem{notation}[theorem]{Notation}
\newtheorem{convention}[theorem]{Convention}

\begin{document}
\title{Multicurves and equivariant cohomology}
\author{N.~P.~Strickland}
\address{
Department of Pure Mathematics\\
University of Sheffield\\
Sheffield S3 7RH\\
UK
}
\email{N.P.Strickland@sheffield.ac.uk}

\date{\today}
\bibliographystyle{abbrv}



\begin{abstract}
 Let $A$ be a finite abelian group.  We set up an algebraic framework
 for studying $A$-equivariant complex-orientable cohomology theories
 in terms of a suitable kind of equivariant formal groups.  We compute
 the equivariant cohomology of many spaces in these terms, including
 projective bundles (and associated Gysin maps), Thom spaces, and
 infinite Grassmannians.
\end{abstract}

\maketitle 

\section{Introduction}
\label{sec-intro}

Let $A$ be a finite abelian group.  In this paper, we set up an
algebraic framework for studying $A$-equivariant complex-orientable
cohomology theories in terms of a suitable kind of formal groups.  In
part, this is a geometric reformulation of earlier work of Cole,
Greenlees, Kriz and others on equivariant formal group
laws~\cites{co:cor,cogrkr:efg,cogrkr:uec,kr:zec}.  However, the theory
of divisors, residues and duality for multicurves is new, and forms a
substantial part of the present paper.  Although we focus on the
finite case, many results can be generalised to compact abelian Lie
groups.  On the other hand, we have evidence that nonabelian groups
will need a completely different theory.

We now briefly recall the nonequivariant theory, using the language of
formal schemes.  We will follow the conventions and terminology
developed in~\cite{st:fsfg}.  Let $E$ be an even periodic cohomology
theory, and put $S=\spec(E^0)$ and
$C=\spf(E^0\CPi_+)=\colim_n\spec(E^0\CP^n_+)$.  Some basic facts are
as follows.
\begin{itemize}
 \item[(a)] $C$ is a formal group scheme over $S$.
 \item[(b)] If we forget the group structure, then $C$ is isomorphic
  to the formal affine line $\haf^1_S$ as a formal scheme over $S$; in
  other words, $C$ is a formal curve over $S$.
 \item[(c)] For many interesting spaces $X$, the formal scheme
  $\spf(E^0X)$ has a natural description as a functor of $C$; for
  example, we have $\spf(E^0BU(d))=C^d/\Sg_d=\Div_d^+(C)$, the formal
  scheme of effective divisors of degree $d$ on $C$.  Similar
  descriptions are known for $\Om U(d)$, $\C P^d$, $BSU$,
  Grassmannians, flag varieties, toric varieties and so on.
 \item[(d)] If $M$ is a compact complex manifold then the ring $E^0M$
  has Poincar\'e duality: there is a map $\tht\:E^0M\to E^0$ such that
  the pairing $(a,b)\mapsto\tht(ab)$ is perfect.  When formulated this
  way, the map $\tht$ is not quite canonical; we need to build in a
  twist by a certain line bundle to cure this.  There is a formula for
  $\tht$ (due to Quillen) in terms of residues of differential forms
  on $C$.
\end{itemize}

Now let $\CU=\CU_A$ be a complete $A$-universe, and let $\CS_A$ be the
category of $A$-spectra indexed on $\CU$ (as in~\cite{lemast:esh}).
Consider an $A$-equivariant commutative ring spectrum $E\in\CS_A$ that
is periodic and orientable in a sense to be made precise later.  In
this context, the right analogue of $\CPi$ is the projective space
$P\CU$.  This has an evident $A$-action.  We put $S=\spec(E^0)$ and
$C=\spf(E^0P\CU)$.  This is again a formal group scheme over $S$, but
it is no longer a formal curve.  This appears to create difficulties
with~(c) above, because we no longer have a good hold on $C^d/\Sg_d$
or a good theory of divisors on $C$.

Our first task is to define the notion of a \emph{formal multicurve}
over $S$, and to show that $C$ is an example of this notion.  Later we
will develop an extensive theory of formal multicurves and their
divisors, and show that many statements about generalized cohomology
can be made equivariant by replacing curves with multicurves.  

\begin{example}
 The simplest nontrivial example of a formal multicurve is the formal
 scheme $C=\spf(R)$, where $R$ is the completion of $\Z[x]$ at the
 ideal generated by the polynomial $f(x)=x^2-2x$, or equivalently
 $R=\Z\psb{y}[x]/(x^2-2x-y)$.  If we consider formal schemes as
 functors from rings to sets, then 
 \[ C(A) = \{a\in A \st f(a) \text{ is nilpotent }\}. \]
 One can check that $R/2=(\Z/2)\psb{x}$ but
 $R/3\simeq(\Z/3)\psb{y}\tm(\Z/3)\psb{y}$. 
\end{example}

\begin{example}
 Let $A$ be a finite abelian group as before, let $R(A)=\Z[A^*]$ be
 the complex representation ring, and put 
 \[ f(u)=\prod_{\al\in A^*}(u-\al) \in R(A)[u]. \]
 Let $R_1$ be the completion of $R(A)[u]$ at the ideal generated by
 $f(u)$, and put $C_1=\spf(R_1)$.  Then $C_1$ is again a formal
 multicurve.  Moreover, it is known that $R_1$ is the equivariant
 complex $K$-theory of $P\CU_A$, so this is an instance of the
 topological situation mentioned above.  This is a theorem of Cole,
 Greenlees and Kriz, which will be discussed further in
 Section~\ref{subsec-MG}. 
\end{example}

\begin{remark}
 It is not true that every formal multicurve has the form
 $\spf(k[x]^\wedge_f)$, but for many purposes one can reduce to that
 special case.  This will be discussed further in
 Section~\ref{sec-embeddings}. 
\end{remark}

\subsection{Outline of the paper}

\begin{itemize}
 \item In Section~\ref{sec-multicurves} we give the most basic
  definitions about multicurves.
 \item In Section~\ref{sec-differential} we set up a theory of
  differential forms on multicurves.
 \item In Section~\ref{sec-projective} we discuss the topology and
  geometry of projective spaces, being careful to use explicit and
  natural methods that transfer automatically to the equivariant
  setting.
 \item In Section~\ref{sec-orientability} we set up the basic theory
  of equivariantly complex orientable ring spectra, and we show that a
  suitable class of such spectra give rise to equivariant formal
  groups.
 \item A nonequivariant even periodic ring spectrum $\hE$ gives two
  different cohomology theories on $A$-spaces, namely the forgetful
  theory $X\mapsto\hE^0(X)$ and the Borel theory
  $X\mapsto\hE^0(X_{hA})$.  In Section~\ref{sec-simple} we discuss the
  (very simple) corresponding ways to convert nonequivariant formal
  groups to equivariant formal groups.
 \item In Section~\ref{sec-algebraic} we show how suitable algebraic
  groups give rise to equivariant formal groups, and we discuss the
  relevance to equivariant $K$-theory and elliptic cohomology.
 \item In Section~\ref{sec-product} we introduce a new class of EFGs
  (those of ``product type'') and show that all EFGs over fields are
  in this class.  This makes it easy to classify EFGs over
  algebraically closed fields, which leads to a discussion of
  equivariant Morava $K$-theories.
 \item It is well-known that the category of rational even periodic
  spectra is contravariantly equivalent to that of formal groups over
  rational rings.  In Section~\ref{sec-rational} we explain and prove
  the corresponding fact for equivariant formal groups.
 \item In Section~\ref{sec-pushout} we introduce another algebraic
  construction on EFGs.  This is used in Section~\ref{sec-E-theory} to
  construct universal deformations of EFGs over fields, and to relate
  them to equivariant Morava $E$-theory.
 \item Various theorems are known to the effect
  that $E^0EG$ is a completion of $E^0(\text{point})$, for suitable
  groups $G$ and equivariant ring spectra $E$.  The earliest such
  result is the Atiyah-Segal theorem for equivariant
  $K$-theory~\cite{atse:ekc}; for some more recent results,
  see~\cite{grma:lct} and the papers referenced there.  In
  Section~\ref{sec-completion} we prove a theorem of this type,
  following the argument of Greenlees~\cite{gr:aie} but interpreting
  everything in terms of EFGs.
 \item In Section~\ref{sec-counterexample} we discuss a salutory
  counterexample.
 \item In Section~\ref{sec-divisors} we return to the algebraic theory
  of divisors on multicurves.  We give two different definitions of
  divisors, and prove the crucial fact that the weaker definition is
  in fact equivalent to the stronger one.
 \item In Section~\ref{sec-embeddings} we define what it means for a
  multicurve to be \emph{embeddable}, and show how to reduce various
  questions to this special case.
 \item In Sections~\ref{sec-symmetric} and~\ref{sec-classify} we
  discuss classification of divisors and the isomorphism
  $\Div_d^+(C)=C^d/\Sg_d$.  Many results are the same as in the
  nonequivariant case, but the proofs are much more technical (and
  involve the reductions established in
  Section~\ref{sec-embeddings}).
 \item In Section~\ref{sec-local-structure} we study the local
  structure of the formal scheme $\Div_d^+(C)=C^d/\Sg_d$, showing that
  in a suitable sense it is a formal manifold of dimension $d$.
 \item In Section~\ref{sec-grassmann} we give a reorganised proof of a
  theorem of Cole, Greenlees and Kriz~\cite{cogrkr:uec} showing that
  the formal scheme associated to the Grassmannian $G_d\CU$ is
  precisely $\Div_d^+(C)$.  This allows us to build a dictionary
  between vector bundles and divisors, just as in the nonequivariant
  case. 
 \item In Section~\ref{sec-thom} we use and extend these ideas to
  give an analogue of the Projective Bundle Theorem, and descriptions
  of the cohomology of flag bundles and Grassmann bundles.
 \item In Section~\ref{sec-duality} we develop an algebraic theory
  for certain kinds of rings with Poincar\'e duality, and apply it to
  understand residues, Gysin maps and duality for projective bundles.
 \item In Section~\ref{sec-grassmann-more} we return to the
  study of Grassmannians.  Greenlees's definition of equivariant
  connective $K$-theory~\cite{gr:evr} gives canonical equivariant
  analogues of the spaces $\Z\tm BU$, $BU$ and $BSU$, and we describe
  the equivariant cohomology of these.
 \item Finally, in Section~\ref{sec-transfer} we show how to construct
  a Mackey functor from an arbitrary $A$-equivariant formal group $C$
  over $S$, and thus a map from the Burnside ring of $A$ to $\OS$.  In
  the case where $C$ arises from a ring spectrum, this map will agree
  with the usual one defined using equivariant topology.
\end{itemize}

\section{Multicurves}
\label{sec-multicurves}

\begin{definition}
 Let $X=\spf(R)$ be a formal scheme, and let $Y$ be a subscheme of
 $X$.  We say that $Y$ is a \emph{regular hypersurface} if
 $Y=\spf(R/J)$ for some ideal $J=I_Y\leq R$ that is a free module of
 rank one over $R$.  Equivalently, there should be a regular element
 $f\in R$ such that the vanishing locus $V(f)=\spf(R/f)$ is precisely
 $Y$.
\end{definition}

Let $S=\spec(k)$ be an affine scheme.

\begin{definition}\label{defn-multicurve}
 A \emph{formal multicurve} over $S$ is a formal scheme $C$ over $S$
 such that
 \begin{itemize}
  \item[(a)] $C=\spf(R)$ for some formal ring $R$
  \item[(b)] There exists a regular element $y\in R$ such that for all
   $k\geq 0$, the ideal $Ry^k$ is open in $R$, and $R/y^k$ is a
   finitely generated free module over $\OS$, and $R=\invlim_k R/y^k$.
  \item[(c)] The diagonal subscheme $\Dl\subset C\tm_SC$ is a regular
   hypersurface.
 \end{itemize}
 A generator $d$ for the ideal $I_\Dl$ will be called a
 \emph{difference function} for $C$ (because $d(a,b)=0$ iff
 $(a,b)\in\Dl$ iff $a=b$).  We will choose a difference function
 $d$, but as far as possible we will express our results in a form
 independent of this choice.  An element $y$ as in~(b) will be called
 a \emph{good parameter} on $C$.
\end{definition}

\begin{remark}\label{rem-formal-base}
 If $S$ is a formal scheme, then we can write $S=\colim_\al S_\al$ for
 some filtered system of affine schemes, and formal schemes over $S$
 are the same as compatible systems of formal schemes over the $S_\al$
 by~\cite{st:fsfg}*{Proposition 4.27}.  In the rest of this paper, we
 will generally work over an affine base but will silently use this
 result to transfer definitions and theorems to the case of a formal
 base where necessary.
\end{remark}

The formal affine line $\haf^1_S=\spf(k\psb{x})$ is a formal
multicurve, and the category of formal multicurves is closed under
disjoint union.  Conversely, condition~(c) implies that the module
$\Om^1_{C/S}=I_\Dl/I_\Dl^2$ is free of rank one over $R=\OC$, so
formal multicurves may be thought of as being smooth and
one-dimensional.  Similarly, if $y$ is a good parameter then $R$ is a
finitely generated projective module over $k\psb{y}$, which means that
$C$ admits a finite flat map to $\haf^1_S$, again indicating a
one-dimensional situation.  If $k$ is an algebraically closed field,
we shall see later that every small formal multicurve over $S$ is a
finite disjoint union of copies of $\haf^1_S$.

\begin{remark}\label{rem-Om-quot}
 Note that $I_\Dl$ is the kernel of the multiplication map
 $\mu\:R\hot R\xra{}R$, which is split by the map $a\mapsto a\ot 1$.
 It follows that $I_\Dl$ is topologically generated by elements of the
 form $a\ot b-ab\ot 1$.  We also see by similar arguments that for any
 ideal $J\leq R$, the kernel of the multiplication map
 $(R/J)\ot(R/J)\xra{}R/J$ is just the image of $I_\Dl$ and thus is
 generated by $d$.
\end{remark}

\begin{construction}\label{cons-top-basis}
 Let $C=\spf(R)$ be a formal scheme over $S=\spec(k)$, and let $y\in
 R$ be a good parameter.  (We are not assuming that $C$ is a formal
 multicurve, but we assume that $y$ satisfies part~(b) of
 Definition~\ref{defn-multicurve}.) Choose elements
 $e_0,\dotsc,e_{n-1}\in R$ whose residue classes give a basis for the
 free module $R/y$ over $k$.  Define elements $e_i\in R$ for all
 $i\in\N$ by $e_{nj+k}=y^je_k$.  Define $\mu_r\:\prod_{i=0}^{nr-1}k\to
 R/y^r$ by $\mu_r(\un{a})=\sum_{i=0}^{nr-1}a_ie_i$.  Define
 $\mu_\infty\:\prod_{i=0}^\infty k\to R$ by
 $\mu_\infty(\un{a})=\sum_ia_ie_i$.  (This is meaningful because
 $R=\invlim_iR/y^i$.)
\end{construction}

\begin{lemma}\label{lem-top-basis}
 The maps $\mu_r$ and $\mu_\infty$ are isomorphisms.  (In other words,
 $\{e_i\}_{i\in\N}$ is a topological basis for $R$.)
\end{lemma}
\begin{proof}
 Define $\mu'_r\:\left(\prod_{i=0}^{nr-1}k\right)\tm R\to R$ by 
 $\mu'_r(\un{a},b)=\left(\sum_ia_ie_i\right)+y^rb$.  As $y$ is a regular
 element and the elements $e_0,\dotsc,e_{n-1}$ form a basis for $R/y$
 we see that $\mu'_1$ is an isomorphism, so $R=k^n\op R$.  We can
 substitute this int itself to see that $R=k^{2n}\op R$; by a
 slightly more precise version of the same argument, we see that
 $\mu'_2$ is an isomorphism.  We can extend this inductively to see
 that $\mu'_r$ is an isomorphism for all $r$.  We then reduce mod
 $y^r$ to see that $\mu_r$ is an isomorphism, and pass to the limit to
 see that $\mu_\infty$ is an isomorphism.
\end{proof}

\begin{definition}
 A \emph{formal multicurve group} over $S$ is a formal multicurve over
 $S$ with a commutative group structure.
\end{definition}

In the presence of a group structure, axiom~(c) can be modified.
\begin{definition}
 Let $C$ be a commutative formal group scheme over $S$.  A
 \emph{coordinate} on $C$ is a regular element $x\in\OC$ whose
 vanishing locus is the zero-section.  Clearly, such an $x$ exists iff
 the zero-section is a regular hypersurface.
\end{definition}
\begin{remark}\label{rem-xbar}
 If $x$ is a coordinate, then so is the function $\xb$ defined by
 $\xb(a)=x(-a)$.  
\end{remark}

\begin{proposition}\label{prop-group-coord}
 Let $C$ be a formal group scheme over $S$ satisfying axioms~(a)
 and~(b) in Definition~\ref{defn-multicurve}.  Then $C$ is a formal
 multicurve iff the zero-section $S\xra{}C$ is a regular hypersurface.
 More precisely, if $x$ is a coordinate on $C$, then the function
 $d(a,b)=x(b-a)$ defines a difference function, and if $d$ is a
 difference function, then the function $x(b)=d(0,b)$ is a
 coordinate. 
\end{proposition}

The proof relies on the following basic lemma.

\begin{lemma}\label{lem-res-reg}
 Let $C$ be a formal multicurve, and let $f\:X\xra{}C$ be any map of
 schemes.  Then the function $d'(x,b)=d(f(x),b)$ on $X\tm_SC$ is
 regular in $\O_{X\tm_SC}$.
\end{lemma}
\begin{proof}
 We have a short exact sequence as follows:
 \[ R\hot R \xra{\tm d} R\hot R \xra{\mu} R. \]
 We regard $R\hot R$ as a module over $R$ via the map
 $t\mapsto t\ot 1$.  The map $\mu$ is then $R$-linearly split by the
 map $t\mapsto t\ot 1$, so the sequence remains exact after applying
 the functor $\O_X\hot_R(-)$.  The resulting sequence is just
 \[ \O_{X\tm_SC} \xra{\tm d'} \O_{X\tm_SC} \xra{} \O_X, \]
 which proves the lemma.
\end{proof}

\begin{corollary}\label{cor-sec-regular}
 Let $C\xra{q}S$ be a formal multicurve, and let $S\xra{u}C$ be a
 section.  Then the subscheme $uS\subset C$ is a regular
 hypersurface, and the ideal $I_{uS}$ is generated by the function 
 $d'(c)=d(u(q(c)),c)$, or equivalently
 \[ d' = (C = S\tm_SC \xra{u\tm 1} C\tm_SC \xra{d} \aff^1). \]
\end{corollary}
\begin{proof}
 Take $X=S$ in the lemma.
\end{proof}

\begin{proof}[Proof of Proposition~\ref{prop-group-coord}]
 First suppose that the zero section is a regular hypersurface, so we
 can choose a coordinate $x$.  It follows easily from axiom~(b) that
 $R\simeq\prod_{k=0}^\infty\OS$ as topological $\OS$-modules, so
 $R\hot R=\prod_{k=0}^\infty R$ as $R$-modules, so $1\ot x$ is a
 regular element in $R\ot R$.  If we regard $x$ as a function on $C$,
 this says that the function $x_1\:(a,b)\mapsto x(b)$ is a regular
 element of $\O_{C\tm_SC}$, whose vanishing locus is precisely the
 closed subscheme where $b=0$.  The map $s\:(a,b)\mapsto(a,b-a)$ is an
 automorphism of $C\tm_SC$, and $s^*x_1$ is the function
 $d(a,b)=x(b-a)$.  As $s$ is an automorphism, we see that $d$ is
 regular and its vanishing locus is the subscheme where $a=b$, or in
 other words the diagonal.

 The converse is the case $u=0$ of Corollary~\ref{cor-sec-regular}. 
\end{proof}

To formulate the definition of an equivariant formal group, we need
some basic notions about divisors.
\begin{definition}\label{defn-divisor}
 A \emph{divisor} on $C$ is a scheme of the form $D=\spec(\OC/J)$,
 where $J$ is an open ideal generated by a single regular element, and
 $\OC/J$ is a finitely generated projective module over $\OS$.  Thus
 $D$ is a regular hypersurface in $C$ and is finite and very flat over
 $S$.  Strictly speaking, we should refer to such subschemes as
 \emph{effective divisors}, but we will have little need for more
 general divisors in this paper.

 If $D_i=\spf(R/J_i)$ is a divisor for $i=0,1$ then we put
 $D_0+D_1:=\spf(R/(J_0J_1))$, which is easily seen to be another
 divisor.

 The \emph{degree} of $D$ is the rank of $\OD$ over $k$.  Note that
 this need not be constant, but that $S$ can be split as a finite
 disjoint union of pieces over which $D$ has constant degree.

 If $T$ is a scheme over $S$, then a divisor on $C$ over $T$ means a
 divisor on the formal multicurve $T\tm_SC$ over $T$.
\end{definition}

Note that if $D$ is an effective divisor of degree one, then the
projection $D\xra{\pi}S$ is an isomorphism, so the map
$S\xra{\pi^{-1}}D\subset C$ is a section of $C$.  Conversely, if
$u\:S\xra{}C$ is a section, then (by Corollary~\ref{cor-sec-regular})
the image $uS$ is a divisor of degree one, which is conventionally
denoted by $[u]$.

In the case of ordinary formal curves, it is well-known that there is
a moduli scheme $\Div_d^+(C)$ for effective divisors of degree $d$ on
$C$, and that it can be identified with the symmetric power
$C^d/\Sg_d$.  Analogous facts are true for multicurves, but much more
difficult to prove.  We will return to this in
Section~\ref{sec-divisors}. 

Let $A$ be a finite abelian group (with the group operation written
additively).  We write $A^*$ for the dual group $\Hom(A,\Q/\Z)$.
This gives us a group scheme 
\[ A^* \tm S = \coprod_{\al\in A^*}S =\spec(\prod_{\al\in A^*}\OS)
\]
over $S$.

\begin{definition}\label{defn-neighbourhood}
 Let $X=\spf(R)$ be a formal scheme, and let $Y=\spf(R/J)$ be a closed
 formal subscheme.  We say that $X$ is a \emph{formal neighbourhood}
 of $Y$ if $R$ is isomorphic to $\invlim_mR/J^m$ as a topological
 ring, or equivalently $X=\colim_m\spf(R/J^m)$, which essentially
 means that every point in $X$ is infinitesimally close to $Y$.
\end{definition}

\begin{definition}\label{defn-efg}
 An \emph{$A$-equivariant formal group} or \emph{$A$-efg} over a scheme
 $S$ is a formal multicurve group $C$ over $S$, together with a
 homomorphism $\phi\:A^*\tm S\xra{}C$, such that $C$ is the formal
 neighbourhood of the divisor
 \[ [\phi(A^*)] := \sum_{\al\in A^*}[\phi(\al)] \subset C. \]
\end{definition}
\begin{remark}
 Choose a coordinate $x$ on $C$, and put $d(a,b)=x(b-a)$.  For any
 $\al\in A^*$ we have a function $x_\al$ on $C$ defined by
 $x_\al(a)=x(a-\phi(\al))=d(\phi(\al),a)$.  More precisely, $x_\al$ is
 the composite
 \[ C = S\tm_SC \xra{\phi(\al)\tm_S1} C\tm_SC
     \xra{\text{subtract}} C \xra{x} \aff^1.
 \]
 The vanishing locus of $x_\al$ is the divisor $[\phi(\al)]$, so the
 vanishing locus of the product $y:=\prod_\al x_\al$ is the divisor
 $[\phi(A^*)]$.  We see using Corollary~\ref{cor-sec-regular} that $y$
 is a regular element in $\OC$.  The final condition in
 Definition~\ref{defn-efg} says that $y$ is topologically nilpotent.
 It is not hard to deduce that $y$ is a good parameter on $C$.
\end{remark}

\begin{proposition}\label{prop-monic}
 Let $f$ be a monic polynomial of degree $d>0$ over $\OS$, and let $R$
 be the completion of $\OS[x]$ at $f$.  Then the scheme
 $C=\spf(R)=\colim_kV(f^k)\sse\aff^1_S$ is a formal multicurve.
\end{proposition}
\begin{proof}
 Condition~(b) is clear, because $\{x^i\st i<dj\}$ is a basis for
 $R/(f^j)$ over $\OS$.  Next, observe that
 \[ R\simeq
     \OS\psb{y}[x]/(f(x)-y) = \OS\psb{y}\{x^i\st i<d\},
 \]
 so
 \[ R\hot R\simeq
     \OS\psb{y_0,y_1}[x_0,x_1]/(f(x_0)-y_0,f(x_1)-y_1) =
     \OS\psb{y_0,y_1}\{x_0^ix_1^j\st i,j<d\}.
 \]
 It is clear that $y_1-y_0$ is not a zero-divisor in this ring, and
 $y_1-y_0=f(x_1)-f(x_0)$ which is divisible by $x_1-x_0$, so $x_1-x_0$
 is also not a zero-divisor.  The multiplication map
 $\mu\:R\hot R\to R$ is surjective, split by the map
 $\sg\:R\to R\hot R$ given by $\sg(f)=f\ot 1$.  The kernel $I_\Dl$ is
 the same as the image of $1-\sg\mu$, or in other words the set of
 functions of the form $f(y_0,y_1,x_0,x_1)-f(y_0,y_0,x_0,x_0)$.  All
 such functions lie in the ideal generated by $(y_1-y_0,x_1-x_0)$ and
 $x_1-x_0$ divides $y_1-y_0$ so $I_\Dl=(x_1-x_0)$, so~(c) holds.
\end{proof}
\begin{definition}
 We say that a formal multicurve $C$ over $S$ is \emph{embeddable} if
 it has the form $\colim_kV(f^k)$ as above for some monic polynomial
 $f$.  
\end{definition}

We will see in Section~\ref{sec-embeddings} that any formal multicurve
can be made embeddable by a faithfully flat base change; this allows
us to reduce many questions to the embeddable case.

\begin{lemma}\label{lem-alg-closed}
 Suppose that $k=\OS$ is an algebraically closed field, and that $C$
 is a formal multicurve over $S$.  Then $C$ is a finite disjoint union
 of copies of $\haf^1_S$, and is embeddable.
\end{lemma}
\begin{proof}
 Let $y\in R$ be a good parameter.  Then the ring $\Rb:=R/y$ is a
 finite-dimensional algebra over the field $k$, so it splits as a
 finite product of local algebras.  As $R$ is complete at $(y)$ we can
 lift this splitting to $R$, which splits $C$ as a disjoint union, say
 $C=C_1\amalg\ldots\amalg C_r$.  One can check that each $C_i$ is
 a formal multicurve.  Put $R_i=\O_{C_i}$, so $R=R_1\tm\ldots\tm R_r$.
 Let $y_i$ be the component of $y$ in $R_i$ and put $\Rb_i=R_i/y_i$,
 so $\Rb=\Rb_1\tm\ldots\tm\Rb_r$.  Moreover, $\Rb_i$ is local, with
 maximal ideal $\mxi_i$ say.  As $k$ is algebraically closed we see
 that $\Rb_i/\mxi_i=k$.  This gives an augmentation
 $u_i^*\:R_i\xra{}k$, or equivalently a section $u_i\:S\xra{}C_i$.  It
 follows from Corollary~\ref{cor-sec-regular} that the kernel of
 $u_i^*$ is generated by a single regular element, say $x_i$.  This
 means that the image of $x_i$ in $\Rb_i$ generates $\mxi_i$, so
 $\Rb_i=k\op\Rb_ix_i$.  Next, the descending chain of
 finite-dimensional $k$-spaces $(x_i^n)$ must eventually stabilise,
 say with $(x_i^n)=(x_i^{n+1})$.  This means that $x_i^n=x_i^{n+1}y$
 for some $y$, so $(1-x_iy)x_i^n=0$ but
 $1-x_iy\in\Rb_i\sm(x_i)=\Rb_i^\tm$ so $x_i^n=0$.  We can combine this
 with the fact that $\Rb_i=k\op\Rb_ix_i$ to see that 
 $\Rb_i\simeq k[x_i]/x_i^m$ for some $m\leq n$, and thus $y_i$ divides
 $x_i^m$.  On the other hand, we clearly have $u^*(y_i)=0$ so $x_i$
 divides $y_i$.  One can now check that the obvious map
 $k[x_i]\to R_i$ extends to give a map $k\psb{x_i}\to R_i$ which is an
 isomorphism, so $C_i\simeq\haf^1_S$.

 Finally, as $k$ is algebraically closed, it is certainly infinite, so
 we can choose distinct elements $\lm_1,\ldots,\lm_r\in k$ say.  If we
 put $f(x)=\prod_i(x-\lm_i)$ we find that the completion of $k[x]$ at
 $(f)$ is isomorphic to $\prod_{i=1}^r k\psb{x}$ and thus to $\OC$.
 This proves that $C$ is embeddable.
\end{proof}

\begin{lemma}
 Let $k$ be an algebraically closed field, and let $A$ be a
 finite-dimensional local $k$-algebra whose maximal ideal is generated
 by a single element $x$.  Then $A=k[x]/x^m$ for some $m$.
\end{lemma}
\begin{proof}
 As $A$ is finite-dimensional the powers of $x$ cannot be linearly
 independent, so there exists a monic polynomial $f(t)\in k[t]$ with
 $f(x)=0$.  We can factor this as $f(t)=t^ng(t)$ with $g(0)\neq 0$.
 Now $g(x)=g(0)\neq 0\pmod{\mxi}$ so $g(x)\not\in\mxi$ so
 $g(x)\in A^\tm$, so $x^n=0$.  
\end{proof}

We next explain how to find a topological basis for the ring $\OC$,
when $(C,\phi)$ is an $A$-equivariant formal group.  Put $n=|A|=|A^*|$
and choose an enumeration of the elements of $A^*$, say
$A^*=\{\al_0,\dotsc,\al_{n-1}\}$.  More generally, for $i\in\N$ we
can write $i=nj+k$ with $0\leq k<n$ and we define $\al_i=\al_k$.  We
then define $e_i\in\OC$ by
\[ e_i(a) = \prod_{j<i} x(a-\phi(\al_j))
          = \prod_{j<i} x_{\al_{j}}(a). 
\]
Note that $e_n$ is just the good parameter $y=\prod_{\al\in A^*}x_\al$
and more generally $e_{nj+k}=y^je_k$, so $e_i\to 0$ in the $y$-adic
topology on $\OC$.  We can thus define a map
$\mu\:\prod_{i\in\N}\OS\to\OC$ by $\mu(t)=\sum_it_ie_i$.  This also
induces a map $\mu_i\:\prod_{j<i}\OS\to\OC/e_i$.

\begin{proposition}\label{prop-top-basis}
 The maps $\mu_i$ and $\mu$ are isomorphisms.
\end{proposition}
\begin{proof}
 We can define a map $\nu_0\:\OS\tm\OC\to\OC$ by $\nu_0(t,u)=t+ux$, and it
 follows from the definition of a coordinate that this is a bijection.
 We can twist this by the translation action of $\phi(\al)$ to see
 that the map $\nu_\al(t,u)=t+ux_\al$ is also a bijection.

 Now define $\mu'_i\:(\prod_{j<i}\OS)\tm\OC\to\OC$ by
 $\mu'_i(t,u)=(\sum_{j<i}t_je_j)+ue_i$.  We see that $\mu'_0=1$ and
 $\mu'_{i+1}=\mu'_i\circ(1\tm\nu_{\al_i})$; it follows that $\mu'_i$ is
 an isomorphism for all $i$.  We can reduce mod $e_i$ to see that
 $\mu_i$ is an isomorphism, and then pass to the limit to see that
 $\mu$ is an isomorphism.
\end{proof}
\begin{remark}
 If we knew in advance that $\mu_n$ was an isomorphism, we could use
 Lemma~\ref{lem-top-basis} to see that $\mu_i$ is an isomorphism for
 all $i$ including $i=\infty$.  However, that approach does not save
 us any effort, as the proof for $\mu_n$ is no easier than the proof
 for $\mu_i$ for all $i$.
\end{remark}

\section{Differential forms}
\label{sec-differential}

We next recall some basic ideas about differential forms, and record
some formulae that will be useful later in our study of residues.

Given a formal multicurve $C$ over $S$, we put
\[ \Om=\Om^1_{C/S}=I_\Dl/I_\Dl^2, \]
and call this the module of differential forms on $C$.  

We also put $\Dl_2=\spf(\O_{C\tm_SC}/I_\Dl^2)$, and regard this as the
second-order infinitesimal neighbourhood of $\Dl$ in $C\tm_SC$.
In these terms, $\Om$ is the module of functions on $\Dl_2$ that
vanish on $\Dl$.

Given a difference function $d\in I_\Dl$, we let $\al$ be the image
of $d$ in $\Om$; this generates $\Om$ freely as a module over $\OC$,
so we can regard $\Om$ as a trivialisable line bundle on $C$.

For any function $f\in\OC$, we write $\bd f$ for the image of
$1\ot f-f\ot 1$ in $\Om$, or equivalently the function
$(a,b)\mapsto f(b)-f(a)$ on $\Dl_2$.  As usual, we have the Leibniz
rule 
\[ \bd(fg) = f\bd(g) + g \bd(f). \]

Now suppose that $C$ has a commutative group structure.  In
particular, this gives a zero-section $Z\subset C$, and we write
$Z_2=\spec(\OC/I_Z^2)$ and 
\[ \om=I_Z/I_Z^2 =
    \{ \text{ functions on $Z_2$ that vanish on $Z$ }\}.
\]
The map $b\mapsto(0,b)$ gives an inclusion $Z_2\xra{}\Dl_2$ and thus a
map $\Om\xra{}\om$, which in turn gives an isomorphism $\Om|_Z=\om$ of
line bundles on $S$.  The image of $\bd f$ under this map is the
element $\bd_0f$ corresponding to the function $b\mapsto f(b)-f(0)$ on
$Z_2$.  If $x$ is a coordinate on $C$, then $\bd_0x$ generates $\om$
freely as a module over $\OS$.

Next, for any function $f\in\OC$ we define a function $\bD f$ on
$\Dl_2$ by
\[ (\bD f)(a,b) = f(b-a) - f(0). \]
This construction gives a map $\bD\:\OC\xra{}\Om$.  If $x$ is a
coordinate then $\bD x$ is the restriction of the usual difference
function $d(a,b)=x(b-a)$ to $\Dl_2$, so it is a generator of $\Om$.

It is easy to see that $\bD f$ depends only on $\bd_0f$, and thus that
$\bD$ induces an $\OS$-linear inclusion $\om\xra{}\Om$, right inverse
to the restriction map $\Om\xra{}\Om|_Z=\om$.  A differential form is
said to be \emph{invariant} if it lies in the image of this map.

By extension of scalars, we obtain an $\OC$-linear map
$\OC\ot_{\OS}\om\xra{}\Om$, sending $f\ot\bd_0g$ to $f\bD g$.  In
particular, it sends $f\ot\bd_0x$ to $f\bD x$, and so is an
isomorphism.

\section{Equivariant projective spaces}
\label{sec-projective}

We now start to build a connection between multicurves and
$A$-equivariant topology (where $A$ is a finite abelian group).
Naturally, this involves the generalised cohomology of the projective
spaces of representations of $A$.  In this section, we assemble some
facts about the homotopy theory of such projective spaces.  

For $\al\in A^*=\Hom(A,\Q/\Z)$ we write $L_\al$ for $\C$ with $A$
acting by $a.z=e^{2\pi i\al(a)}z$.  In particular, $L_0$ has trivial
action, and $L_\al\ot L_\bt=L_{\al+\bt}$.  For any finite-dimensional
representation $V$, we put
\[ V[\al] =
    \{v \in V \st av=e^{2\pi i\al(a)}v \text{ for all } a\in A\}
      \simeq \Hom_A(L_\al,V)\ot L_\al.
\] 
It is well-known that $V=\bigoplus_\al V[\al]$ and
$\Hom_{\C[A]}(V,W)=\bigoplus_\al\Hom_\C(V[\al],W[\al])$.  It follows
that if there exists an equivariant linear embedding $V\xra{}W$, then
the space of such embeddings is connected, giving a canonical homotopy
class of maps $PV\xra{}PW$ of projective spaces.

We write $\CU[\al]=L_\al\ot\C^\infty$, and
$\CU=\CU_A=\bigoplus_\al\CU[\al]$, so $\CU$ is a complete $A$-universe.
We write $P\CU$ for the projective space associated to $\CU$, which
has a natural $A$-action.  By the previous paragraph, for any
finite-dimensional representation $V$, there is a canonical map
$PV\xra{}P\CU$ up to homotopy.  Similarly, the space of equivariant
linear isometries $\CU\ot\CU\xra{}\CU$ is contractible, which gives a
canonical homotopy class of maps $P\CU\tm P\CU\xra{}P\CU$, making
$P\CU$ an abelian group up to equivariant homotopy.  We can choose a
conjugate-linear equivariant automorphism $\chi\:\CU\xra{}\CU$, and
the resulting map $P\CU\xra{}P\CU$ is the negation map for our group
structure.  

It is well-known that $P\CU$ is the classifying space for equivariant
complex line bundles.  More precisely, for any $A$-space $X$, we write
$\Pic_A(X)$ for the group of isomorphism classes of equivariant
complex line bundles over $X$.  Let $T$ denote the tautological line
bundle over $P\CU$, so $T\in\Pic(P\CU)$.  Then for any $A$-space $X$,
the construction $[f]\mapsto[f^*T]$ gives a group isomorphism
$[X,P\CU]^A\simeq\Pic_A(X)$.  Note that we regard $T$ as the universal
example; some other treatments in the literature use the dual bundle
$T^*=\O(1)$ instead.

Note that $A$ acts by scalars on $\CU[\al]$, and thus acts as the
identity on $P\CU[\al]\subset P\CU$.  Moreover, the map $L\mapsto
L_\al\ot L$ gives a homeomorphism $\CPi=P(\C^\infty)\xra{}P\CU[\al]$.
Using this, we have a homeomorphism
$(P\CU)^A=\coprod_{\al}P\CU[\al]=A^*\tm\CPi$, and thus a bijection
$\pi_0((P\CU)^A)=A^*$, which is easily seen to respect the natural
group structures.  Thus, the group structure on $P\CU$ gives a
translation action (up to homotopy) of $A^*$ on $P\CU$.  
\begin{definition}\label{defn-tau} 
 We write $\tau_\al\:P\CU\xra{}P\CU$ for translation by an element
 $\al\in A^*$.
\end{definition}

For various purposes we will need to use an $A$-fixed basepoint in
$P\CU$.  We have embeddings $L_\al\xra{}\CU[\al]\xra{}\CU$, and
$PL_\al$ is an $A$-fixed point.  Any other fixed point lies in the
same component of $(P\CU)^A$ as $PL_\al$ for some $\al$, so it can be
replaced by $PL_\al$ for most purposes.  Moreover, the map $\tau_\al$
gives a homotopy equivalence of pairs
$(P\CU,PL_\bt)\xra{}(P\CU,PL_{\al+\bt})$.  Where not otherwise stated,
we use $PL_0$ as the basepoint.

\begin{proposition}\label{prop-proj-diag}
 Let $V$, $W$ and $X$ be unitary representations of $A$, where $V$ and
 $W$ have finite dimension and $X$ is a colimit of finite-dimensional
 subrepresentations.  Put $U=V\op W\op X$.  Then there is a
 homotopy-commutative diagram as follows, in which the maps marked $q$
 are the obvious collapses, the maps marked $j$ are the obvious
 inclusions, and $\dl$ is the diagonal map.
 \[ \xymatrix{
  {PU} 
   \rrto^{\dl}
   \dto_{q_{V\op W}} & & 
  {PU\tm PU}
   \dto^{q_V\Smash q_W} \\
  {PU/P(V\op W)} \rto_(0.35){\dlb} & 
  {P(V\op X)/PV\Smash P(W\op X)/PW} \rto_(0.58){j\Smash j} & 
  {PU/PV\Smash PU/PW}
  } \]
 Moreover, if $\dim(X)=1$ then $\dlb$ is just the standard
 homeomorphism
 \[ S^{\Hom(X,V\op W)} = S^{\Hom(X,V)} \Smash S^{\Hom(X,W)}. \]
 All maps and homotopies are natural for isometric embeddings of
 $V$, $W$ and $X$.
\end{proposition}
\begin{remark}
 For any ring spectrum $E$, the above diagram gives a map 
 \[ \dlb^*\: E^*(P(V\op X),PV) \ot E^*(P(W\op X),PW)
     \xra{} E^*(PU,P(V\op W)).
 \]
 In his unpublished thesis~\cite{co:cor}, Cole writes $a*b$ for
 $\dlb^*(a\ot b)$.  The idea of using this construction seems to be
 original to that thesis; our approach differs only in being somewhat
 more geometric.
\end{remark}
\begin{proof}
 Assume for the moment that $X$ is finite-dimensional.  We start by
 defining a map 
 \[ \gmb\:PU/P(V\op W)\xra{}PU/PV\Smash PU/PW, \]
 which will be homotopic to $(j\Smash j)\circ\dlb$.  For
 $u=(v,w,x)\in U^\tm := U\sm\{0\}$ we put
 \[ s = s(u) = (\|w\|-\|v\|)/(\|v\|+\|w\|+\|x\|). \]
 Note that $s(u)\in [-1,1]$, and $s(\lm u)=s(u)$ for all
 $\lm\in\C^\tm$, and $s(u)\geq 0$ iff $\|w\|\geq\|v\|$.  We next
 define $\al,\bt\:U^\tm\xra{}U$ by
 \begin{align*}
  \al(v,w,x) &=
   \begin{cases}
    ((1-s)v, sw,x) & \text{ if } s\geq 0 \\
    (v,0,x)        & \text{ if } s\leq 0
   \end{cases} \\
  \bt(v,w,x) &=
   \begin{cases}
    (0,w,x)        & \text{ if } s\geq 0 \\
    (-sv,(1+s)w,x) & \text{ if } s\leq 0. 
   \end{cases}
 \end{align*}
 Note that $\al(\lm u)=\lm\al(u)$ and similarly for $\bt$.

 We claim that when $u\neq 0$, the line joining $u$ to $\al(u)$ never
 passes through $0$ (so in particular $\al(u)\neq 0$).  Indeed, if
 $s\leq 0$, then the points on the line have the form $(v,tw,x)$ for
 $0\leq t\leq 1$.  Thus, the line can only pass through zero if
 $v=x=0$.  The relation $s\leq 0$ means that $\|w\|\leq\|v\|=0$, so
 $w=0$ as well, contradicting the assumption that $u\neq 0$.  In the
 case $s>0$, the points on the line have the form
 $((1-ts)v,(1-t+ts)w,x)$.  As $s>0$ and $0\leq t\leq 1$ we have
 $1-t+ts>0$.  For the line to pass through zero we must thus have
 $x=w=0$, and the relation $s\geq 0$ means that $\|v\|\leq\|w\|=0$,
 again giving a contradiction.  Similarly, the line from $u$ to
 $\bt(u)$ never passes through $0$.

 It follows that $\al$ and $\bt$ induce self-maps of $PU$ that are
 homotopic to the identity, so the map
 $\gm=(\al,\bt)\:PU\xra{}PU\tm PU$ is homotopic to the diagonal map
 $\dl$.  

 Next, note that if $u\in V\op W$, then for $s\geq 0$ we have
 $\gm(u)\in U\tm W$, and for $s\leq 0$ we have $\gm(u)\in V\tm U$.  It
 follows that the induced map on projective spaces has
 \[ \gm(P(V\op W)) \sse (PU\tm PW) \cup (PV\tm PU), \]
 so there is an induced map 
 \[ \gmb\: PU/P(V\op W) \xra{} PU/PV \Smash PU/PW. \] 
 As $\gm$ is homotopic to $\dl$, we see that
 $\gmb\circ q_{V\op W}\simeq(q_V\Smash q_W)\circ\dl$.

 To construct the map $\dlb$, we need a slightly different model.
 Clearly 
 \[ PU\sm P(V\op W) =
    (V\tm W\tm X^\tm)/\C^\tm =
    (V\tm W\tm S(X))/S^1, 
 \] 
 and $PU/P(V\op W)$ is the one-point compactification of this.
 Similarly, $P(V\op X)/PV\Smash P(W\op X)/PW$ is the one-point
 compactification of the space $(V\tm S(X))/S^1\tm(W\tm S(X))/S^1$.
 We can thus define $\dlb$ by giving a proper map
 \[ V\tm W \tm S(X) \xra{} V \tm S(X) \tm W\tm S(X) \]
 with appropriate equivariance.  The map in question just sends
 $(v,w,x)$ to $(v,x,w,x)$. 

 If $X$ is one-dimensional and $(v,x)\in V\tm S(X)$ then we have a
 linear map $\al\:X\xra{}V$ given by $\al(x)=v$, which does not change
 if we multiply $(v,x)$ by an element of $S^1$.  This gives a
 homeomorphism $(V\tm S(X))/S^1=\Hom(X,V)$, and thus
 $P(V\op X)/PV=S^{\Hom(X,V)}$.  It is easy to see that with this
 identification, $\dlb$ is just the standard homeomorphism 
 \[ S^{\Hom(X,V\op W)} = S^{\Hom(X,V)} \Smash S^{\Hom(X,W)}. \]

 We now show that $(j\Smash j)\circ\dlb\simeq\gmb$.  Put 
 \[ T = \{((v_0,w_0,x_0),(v_1,w_1,x_1))\in U^2 \st
           \|(w_0,x_0)\| = \|(v_1,x_1)\| = 1 \},
 \]
 so that $PU/PV \Smash PU/PW$ is the one-point compactification of
 $T/(S^1\tm S^1)$.  Define maps 
 \[ \tht_t \: V\tm W\tm S(X) \xra{} T \]
 for $0\leq t\leq 1$ by
 \[ \tht_t(v,w,x) = \begin{cases}
     \left(
      \frac{((1-st)v,stw,x)}{\|(stw,x)\|},
      (0,w,x)
     \right) & \text{ if } s \geq 0 \\
     \left(
      (v,0,x),
      \frac{(-stv,(1+st)w,x)}{\|(-stv,x)\|}
     \right) & \text{ if } s \leq 0,
    \end{cases}
 \]
 where $s=(\|w\|-\|v\|)/(\|v\|+\|w\|+\|x\|)$ as before.  (Note that
 both clauses give $\tht_t(v,w,x)=((v,0,x),(0,w,x))$ if $s=0$, so the
 two clauses are consistent.)

 We claim that the maps $\tht_t$ are proper.  To see this, put
 \[ \nu((v_0,w_0,x_0),(v_1,w_1,x_1)) = \max(\|v_0\|,\|w_1\|), \]
 and $T_k=\{t\in T\st\nu(t)\leq k\}$.  Every
 compact subset of $T$ is contained in some $T_k$, so it will be
 enough to show that $\tht_t^{-1}T_k$ is compact.  In the case
 $s\geq 0$ we have $0\leq 1-st\leq 1$ and $\|(stw,x)\|\geq\|x\|=1$ so
 $\|((1-st)v/\|(stw,x)\|)\|\leq\|v\|\leq\|w\|$, so
 $\nu(\tht_t(v,w,x))=\|w\|$.  Similarly, when $s\leq 0$ we have
 $\nu(\tht_t(v,w,x))=\|v\|$, so in general
 $\nu(\tht_t(v,w,x))=\max(\|v\|,\|w\|)$.  It follows immediately that
 $\tht_t$ is proper, and we get an induced family of maps
 \[ \tht_t \: PU/P(V\op W) \xra{} PU/PV\Smash PU/PW. \]
 We see from the definitions $\tht_0=(j\Smash j)\circ\dlb$ and
 $\tht_1=\gmb$.  The proposition follows easily (for the case where
 $X$ has finite dimension).

 If $X$ has infinite dimension, we apply the above to all finite
 dimensional subrepresentations of $X$.  We see by inspection that all
 constructions pass to the colimit, so the conclusion is valid for $X$
 itself. 
\end{proof}

By an evident inductive extension, we obtain the following:
\begin{corollary}\label{cor-proj-diag}
 Let $L_1,\ldots,L_d$ be one-dimensional representations of $A$, and
 let $X$ be as above.  Put $Y=\bigoplus_iL_i$ and $U=Y\op X$.  Then
 there is a homotopy-commutative diagram as follows:
 \[ \xymatrix{
  {PU} \dto_{q} \rrto^{\dl} & & 
  {PU^r} \dto^{q} \\
  {PU/PY}  \rto_(0.35){\dlb} & 
  {\bigSmash_i P(L_i\op X)/PL_i} \rto_(0.6){j} &
  {\bigSmash_i PU/PL_i}
  } \]
 Moreover, if $\dim(X)=1$ then $\dlb$ is just the standard
 homeomorphism
 \[ S^{\Hom(X,Y)} = \bigSmash_i S^{\Hom(X,L_i)}. \qed \]
\end{corollary}

We mention one more useful special case.
\begin{corollary}\label{cor-module}
 For any ring spectrum $E$, the group $E^*(P(V\op W),PV)$ is naturally
 a module over $E^*PW$.
\end{corollary}
\begin{proof}
 Take $W=0$ in Proposition~\ref{prop-proj-diag} to get a map
 \[ \dlb\:P(V\op X)/PV\to (P(V\op X)/PV)\Smash PX_+, \]
 and thus a map 
 \[ \dlb^*\:E^0PX\ot E^0(P(V\op X),PV) \to E^0(P(V\op X),PV). \]
 If we identify $P(V\op X)/PV$ as the one-point compactification of
 $(V\tm S(X))/S^1$ then the formula is just $\dlb([v,x])=([v,x],[x])$;
 this is clearly coassociative and counital, so the corresponding map
 in cohomology gives a module structure.  A slight change of notation
 recovers the stated corollary.
\end{proof}

We conclude with some further miscellaneous observations about the
space $P\CU$.
\begin{proposition}\label{prop-PU-FEACP}
 The space $P\CU$ is equivariantly equivalent to $F(EA_+,\CPi)$ (where
 $\CPi$ is the usual space with trivial $A$-action).  Equivalently,
 $P\CU$ is the second space in the Borel cohomology spectrum
 $F(EA_+,H)$, so $[X,P\CU]^A=H^2(X_{hA})$ for any $A$-space $X$.
 Moreover, the space $\Om P\CU$ is equivariantly equivalent to $S^1$
 with the trivial action.
\end{proposition}
\begin{proof}
 There is an evident inclusion
 $\CPi=P(\CU^A)\xra{}(P\CU)^A\xra{}P\CU$, which is a nonequivariant
 equivalence.  It follows that the resulting map
 $F(EA_+,\CPi)\xra{}F(EA_+,P\CU)$ is an equivariant equivalence (see
 Lemma~\ref{lem-FEG}).  On the other hand, the collapse map
 $EA_+\xra{}S^0$ gives a map
 $j\:P\CU\xra{}F(EA_+,P\CU)\simeq F(EA_+,\CPi)$.  We claim that this
 is an equivalence.  Indeed, if we take fixed points for a subgroup
 $A_0\leq A$ we get a map $A_0^*\tm\CPi\xra{}F((BA_0)_+,\CPi)$ of
 commutative $H$-spaces.  It is clear that
 \[ \pi_k(A_0^*\tm\CPi) = 
     \begin{cases}
      A_0^* & \text{ if } k = 0 \\
      \Z    & \text{ if } k = 2 \\
      0     & \text{ otherwise. }
     \end{cases}
 \]
 On the other hand, we have 
 \[ \pi_kF((BA_0)_+,\CPi) = [\Sg^k(BA_0)_+,K(\Z,2)] = H^{2-k}BA_0. \]
 This clearly vanishes for $k>2$ and gives $\Z$ for $k=2$.  The short
 exact sequence $\Z\xra{}\Q\xra{}\Q/\Z$ gives long exact sequences of
 cohomology groups, using which we find that $H^1BA_0=0$ and
 $H^2BA_0=A_0^*$.  This shows that $\pi_*F((BA_0)_+,\CPi)$
 is abstractly isomorphic to $\pi_*(A_0^*\tm\CPi)$.  With a little
 more work one sees that the isomorphism is induced by $j$, and the
 first part of the proposition follows.

 We now see that
 \[ \Om P\CU \simeq \Om F(EA_+,\CPi) = F(EA_+,\Om\CPi) = F(EA_+,S^1).
 \]
 As above we find that 
 \[ \pi_k(F(EA_+,S^1)^{A_0}) = H^{1-k}BA_0 = 
     \begin{cases} 
      \Z & \text{ if } k=1 \\
      0  & \text{ otherwise. }
     \end{cases}
 \]
 It follows that the obvious map $S^1\xra{}F(EA_+,S^1)$ is an
 equivariant equivalence.
\end{proof}

For the convenience of the reader we record a proof of a standard
lemma that was used above.
\begin{lemma}\label{lem-FEG}
 If $f\:X\to Y$ is an $A$-equivariant based map and a nonequivariant
 homotopy equivalence, then the induced map $F(EA_+,X)\to F(EA_+,Y)$
 is an equivariant homotopy equivalence.
\end{lemma}
\begin{proof}
 Let $\CC$ be the category of those based $A$-spaces $C$ for which
 $f_*\:F(C,X)\to F(C,Y)$ is an equivariant equivalence.  For any based
 $A$-spaces $P,Q$ we have 
 \[ [P,F(A_+\Smash Q,X)]^A = [A_+\Smash P\Smash Q,X]^A =
     [P\Smash Q,X].
 \]
 Using this and the Yoneda lemma we see that $A_+\Smash Q\in\CC$.
 Moreover, the category $\CC$ is closed under homotopy pushouts, so it
 follows by cellular induction that it contains all finite free based
 $A$-CW-complexes.  Moreover, $\CC$ is closed under telescopes, so one
 can take a colimit over skeleta to see that $\CC$ contains all free
 based $A$-CW-complexes.  In particular we have $EA_+\in\CC$, as
 required. 
\end{proof}

\begin{proposition}
 Let $T$ be the tautological line bundle over $P\CU$, and let $S(T^n)$
 be the unit circle bundle in the $n$'th tensor power of $T$.  Then
 $S(T^n)$ is equivariantly equivalent to $F(EA_+,B(\Z/n))$.
\end{proposition}
\begin{proof}
 First, we let $\Z/n$ act freely on the contractible space $S(\CU)$ by
 multiplication by $n$'th roots of unity, so $S(\CU)/(\Z/n)$ is a
 model for $B(\Z/n)$.  Given a point $v\in S(\CU)$ we have a line
 $L=\C v\in P(\CU)$ and an element $v^{\ot n}\in L^n$, giving a
 point $(L,v^{\ot n})\in S(T^n)$ that depends only on the $\Z/n$-orbit
 of $v$.  This construction gives a homeomorphism
 $S(\CU)/(\Z/n)\to S(T^n)$, so $S(T^n)$ is also a model
 (nonequivariantly) for $B(\Z/n)$. 

 We now analyse the equivariant picture.  Suppose that $(L,u)\in
 S(T^n)$ is fixed by a subgroup $A_0\leq A$.  We see that $A_0$ acts
 on $L$ by some character $\al\in A_0^*$, so $A_0$ acts on $u$ by
 $n\al$, but $u$ is fixed so $n\al=0$.  Given that $n\al=0$, we see
 that every point in $L^n$ is fixed by $A_0$.  Using this, we see that
 $S(T^n)^{A_0}=A_0^*[n]\tm B(\Z/n)$, where $A_0^*[n]$ denotes the
 subgroup of points of order $n$ in $A_0^*$.  Using this, we find that
 $\pi_*S(T^n)^{A_0}=H^{1-*}(BA_0;\Z/n)$, and the claim follows by the
 same method as in the previous proposition.
\end{proof}

It is also useful to describe $S(\CU)$ as the universal space for a
family of subgroups, as in the following definition.
\begin{definition}\label{defn-EF}
 Let $G$ be a finite group, and let $\CF$ be a collection of subgroups
 such that 
 \begin{itemize}
  \item[(a)] $1\in\CF$
  \item[(b)] If $K\leq H\in\CF$ then $K\in\CF$
  \item[(c)] If $H\in\CF$ and $g\in G$ then $gHg^{-1}\in\CF$.
 \end{itemize}
 Then $E\CF$ is the $G$-space characterised up to unique equivariant
 homotopy equivalence by the following properties:
 \begin{itemize}
  \item[(A)] $(E\CF)^H$ is contractible if $H\in\CF$
  \item[(B)] $(E\CF)^H=\emptyset$ if $H\not\in\CF$.
 \end{itemize}
 We also write $\tE\CF$ for the unreduced suspension $\tSg E\CF$, or
 equivalently the cofibre of the collapse map $E\CF_+\to S^0$. 
\end{definition}
(For the existence and uniqueness of $E\CF$, see~\cite{el:sfp}.
Alternatively, uniqueness is straightforward by obstruction theory,
and we will have concrete models to prove existence in all cases that
we need.)

\begin{proposition}\label{prop-SU-EF}
 Put $\CF=\{B\leq S^1\tm A\st B\cap S^1=\{1\}\}$, which is a family of
 subgroups of $S^1\tm A$.  Then the unit sphere $S(\CU)$ is a model
 for $E\CF$, and so $P\CU=(E\CF)/S^1$.
\end{proposition}
\begin{proof}
 First, we let $S^1\subset\C^\tm$ act on $S(\CU)$ by multiplication,
 and let $A$ act in the usual way.  These actions commute and so give
 an action of $S^1\tm A$.  We need only check that $S(\CU)$ has the
 characterizing property of $E\CF$, or in other words that $S(\CU)^B$
 is contractible for $B\in\CF$ and empty for $B\not\in\CF$.  If
 $B\in\CF$ then $B\cap S^1$ is trivial so $B$ is the graph of a
 homomorphism $\phi\:A_0\xra{}S^1$ for some subgroup $A_0\leq A$.  Put
 \[ \CV = \{v\in\CU\st a.v=\phi(a)^{-1}v\text{ for all } a\in A_0\},
 \]
 so $S(\CU)^B=S(\CV)$.  As $\CU$ is a complete $A_0$-universe, we see
 that $\CV$ is infinite dimensional, and so $S(\CV)$ is contractible
 as required.  On the other hand, as $S^1$ acts freely on $S(\CU)$, it
 is clear that $S(\CU)^B=\emptyset$ whenever $B\not\in\CF$.
\end{proof}

\section{Equivariant orientability}
\label{sec-orientability}

Now let $E$ be a commutative $A$-equivariant ring spectrum.  We next
need to formulate suitable notions of orientability and periodicity
for $E$, and deduce consequences for the rings $E^*PV$.  Our results
differ from those of~\cite{co:cor} only in minor points of technical
detail.  We start by introducing some notation and auxiliary ideas.

\begin{convention}\label{conv-complete}
 All $A$-spectra are implicitly assumed to be indexed by a complete
 $A$-universe.
\end{convention}

\begin{notation}\label{notn-subscript}
 Given an $A$-equivariant spectrum $X$, we write $E^n(X)$ for the
 group $[X,\Sg^n E]^A$ of equivariant homotopy classes of equivariant
 maps.  In many parts of the literature this is written as $E^n_A(X)$,
 but we omit the subscript to avoid cluttering the notation.  If we
 have a $B$-spectrum $Y$ (for some subgroup $B\leq A$) then we will
 write $\res^A_B(E)^n(Y)$ or $E^n(A_+\Smash_BY)$ or $[Y,\Sg^nE]^B$ for
 the group that might otherwise be denoted $E^n_B(Y)$.  For specific
 choices of $E$ the subscript may reappear as part of the name of
 $E$.  For example, $A$-equivariant complex $K$-theory is represented
 by an $A$-spectrum called $KU_A$, and we write $KU^n_A(X)$ for
 $[X,\Sg^nKU_A]^A$.

 Similarly, we write $\pi_n(E)$ rather than $\pi^A_n(E)$ for the group
 $[S^n,E]^A=E^0(S^n)=E^{-n}(S^0)$.  When $X$ is a based $A$-space we
 usually write $\Sgi X$ (rather than $\Sgi_A X$) for the corresponding
 $A$-equivariant suspension spectrum.
\end{notation}

\begin{notation}\label{notn-reduced}
 Another potential source of clutter is the distinction between
 reduced and unreduced cohomology groups.  When $Y$ is an
 $A$-spectrum, we write $E^*(Y)=[Y,E]^A_{-*}$ as above.  When $X$ is
 an $A$-space, we write $E^*(X)$ for the unreduced groups, so
 $E^*(X)=E^*(\Sgip X)$.  If $X$ has an $A$-fixed basepoint then we can
 also define reduced groups $\tE^*(X)=E^*(\Sgi X)$.  Sometimes we will
 use notation that is slightly ambiguous, in that it could refer to a
 based space or to its suspension spectrum.  In such cases we will use
 either a tilde or an explicit $\Sgi$ to resolve the ambiguity.
\end{notation}

\begin{definition}\label{defn-univ-gen}
 Let $R$ be an $E$-algebra spectrum, and $M$ a module spectrum over
 $R$.  We say that $M$ is a free $R$-module if it is equivalent as an
 $R$-module to a wedge of (unsuspended) copies of $R$, or
 equivalently, there is a family of elements $e_i\in\pi_0M$ such that
 the resulting maps $\bigoplus_i[\Sg^nA/B_+,R]\xra{}[\Sg^nA/B_+,M]$
 are isomorphisms for all $n\in\Z$ and all $B\leq A$.  We say that
 such elements $e_i$ are \emph{universal generators} for for $\pi_0M$
 over $\pi_0R$.  We will often leave the identification of $R$ and $M$
 implicit.  For example, if we say that an element $e$ is a universal
 generator for $E^0(X,Y)$ over $E^0X$, we are referring to the case
 $R=F(X_+,E)$ and $M=F(X/Y,E)$.
\end{definition}

\begin{definition}\label{defn-complex-coord}
 Let $E$ an $A$-equivariant ring spectrum, and consider a class
 $x\in E^0(P\CU,PL_0)$.  For any $\al\in A^*$ we can embed
 $L_\al\op L_0$ in $\CU$, and thus restrict $x$ to get a class
 $u_{L_\al}\in E^0(P(L_\al\op L_0),PL_0)=\tE^0S^{L_\al}$.  This in
 turn gives an $E$-module map $m_\al\:\Sg^{L_\al}E\xra{}E$.

 We say that $x$ is a \emph{complex coordinate} for $E$ if for all
 $\al$ the map $m_\al$ is an equivalence, or equivalently $u_{L_\al}$
 generates $\Sg^{-L_\al}E$ as an $E$-module.  We say that $E$ is
 \emph{periodically orientable} if it admits such a coordinate.  We
 say that $E$ is \emph{evenly orientable} if in addition, the group
 $\pi_1^BE=[\Sg A/B_+,E]^A=E^{-1}(A/B)$ vanishes for all $B\leq A$. 
\end{definition}

From now on, we assume that $E$ is periodically orientable.  We choose
a complex coordinate $x$, but as far as possible we state our results
in a form independent of this choice.  We write $\xb=\chi^*x$, where
$\chi\:P\CU\xra{}P\CU$ is the negation map for the group structure.
It is easy to see that this is again a coordinate.

Recall that for any line bundle $L$ over $X$, there is an essentially
unique map $f_L\:X\xra{}P\CU$ with $f^*T\simeq L$ (where $T$ is the
tautological bundle over $P\CU$).  We define the \emph{Euler class} of
$L$ by
\[ e(L) = f^*_{L^*}(x) = f_L^*(\xb). \]
Thus, the element $x\in E^0P\CU$ is the Euler class of $T^*$, and
$\xb$ is the Euler class of $T$.
\begin{remark}\label{rem-euler-convention}
 There is some inconsistency in the literature about whether $e(L)$
 should be $f_L^*(x)$ or $f_L^*(\xb)$.  The convention adopted here is
 the opposite of that used in~\cite{st:fsfg}, but I believe that it is
 more common in other work and has some technical advantages.  The
 conventions used elsewhere in this paper are fixed by the following
 requirements. 
 \begin{itemize}
  \item[(a)] We have $e(V\op W)=e(V)e(W)$.
  \item[(b)] The Euler class of $V$ is the restriction of the Thom
   class in $\tE^0X^V$ to the zero section $X\subset X^V$. 
 \end{itemize}
 Our substitute for the nonequivariant theory of Chern classes will be
 more abstract, so we will not need sign conventions.  The r\^ole
 normally played by the Chern polynomial
 $\sum_{i+j=\dim(V)}\pm c_ix^j$ will be played by a certain element
 $f_V$; if $A=0$ and $V=\bigoplus_iL_i$ then
 $f_V=\prod_i(x+_Fe(L_i))$.  
\end{remark}

Next note that we can define 
\[ x_\al := \tau_{-\al}^*x\in E^0(P\CU,PL_\al).  \] 
(Here $\tau_\al$ is the translation map as in
Definition~\ref{defn-tau}.)  Because
\[ (\tau_{-\al}^*T)^* = (L_{-\al}\ot T)^* = L_\al\ot T^* =
    \Hom(T,L_\al),
\]
we have $x_\al=e(\Hom(T,L_\al))$.  If $L$ is a one-dimensional
representation isomorphic to $L_\al$, we also use the notation $x_L$
for $x_\al$.  We can identify $E^0(P(L_\bt\op L_\al),PL_\al)$ with
$\tE^0S^{L_{\bt-\al}}$, and we find that $x_\al$ restricts to
$u_{\bt-\al}$, which is a universal generator.

Now consider a finite-dimensional representation $V$ of $A$.  We have
a canonical homotopy class of embeddings $PV\xra{}P\CU$, and thus a
well-defined group $E^0(P\CU,PV)$.  We can write $V$ as
$\bigoplus_{i=1}^dL_i$, and Corollary~\ref{cor-proj-diag} gives
a map
\[ P\CU/PV \xra{} \bigSmash_i P\CU/PL_i \]
compatible with the diagonal.  Using this, we can pull back
$x_{L_1}\Smash\ldots\Smash x_{L_d}$ to get a class 
$x_V\in E^0(P\CU,PV)$ that maps to $\prod_ix_{\al_i}$ in $E^0P\CU$.  
Note that for any representation $W$ containing $V$ we can choose an
embedding $W\xra{}\CU$ and pull back $x_V$ along the resulting map
$PW\xra{}P\CU$ to get a class in $E^0(PW,PV)$, which we again denote
by $x_V$.

\begin{lemma}\label{lem-flag-gen}
 Let $V\leq W$ be complex representations of $A$, with $\dim(W/V)=1$.
 Then $x_V$ is a universal generator for $E^0(PW,PV)$.
\end{lemma}
\begin{proof}
 Write $V=L_1\op\ldots\op L_d$ as before, and $X=W\ominus V$, so
 $W=V\op X$ and $PW/PV=S^{\Hom(X,V)}=\bigSmash_iS^{\Hom(X,L_i)}$.
 Because $x$ is a complex coordinate, we know that
 $x_{L_i}\in E^0(PW,PL_i)$ restricts to a universal generator $v_i$ of
 $S^{\Hom(X,L_i)}$.  It follows from Corollary~\ref{cor-proj-diag}
 that $x_V=\prod_iv_i\in\tE^0S^{\Hom(X,V)}=E^0(PW,PV)$, and this is
 easily seen to be a universal generator.
\end{proof} 

\begin{corollary}\label{cor-flag-basis}
 Let $0=U_0<U_1<\ldots<U_d=U$ be representations of $A$ with
 $\dim(U_i)=i$.  Then $\{x_{U_i}\st i<d\}$ is a universal basis for
 $E^0PU$ over $E^0$.
\end{corollary}
\begin{proof}
 This follows by an evident induction from the lemma.
\end{proof}
\begin{remark}\label{rem-xbar-basis}
 As $\xb$ is another coordinate, it gives rise to another universal
 basis $\{\xb_{U_i}\st i<d\}$ for $E^0PU$, which is sometimes more
 convenient.
\end{remark}

We record separately some easy consequences that are independent of
the choice of flag $\{U_i\}$:
\begin{proposition}\label{prop-EPV}
 Let $U$ be a $d$-dimensional representation of $A$.  Then
 \begin{itemize}
  \item[(a)] $F(PU_+,E)$ is a free module of rank $d$ over $E$.
  \item[(b)] If $U=V\op W$ then the restriction
   map $F(PU_+,E)\xra{}F(PV_+,E)$ is split surjective.  The kernel is
   a free module of rank one over $F(PW_+,E)$, generated by
   $x_V$. \qed  
 \end{itemize}
\end{proposition}

We now put $S=\spec(E^0)$ and $R=E^0P\CU$ and $C=\spf(R)$.  We must
show that $C$ is an equivariant formal group over $S$.

We first exhibit a topological basis for $R$.  This will be
essentially the same as in Proposition~\ref{prop-top-basis}, but we
cannot appeal to that result because we do not yet know that we have a
multicurve.  We can list the elements of $A^*$ as 
\[ A^* = \{ \al_0=0,\al_1,\ldots,\al_{n-1} \} \]
(where $n=|A|$), and then define $\al_k$ for all $k\geq 0$ by
$\al_{ni+j}=\al_j$.  We then have an evident filtration
\[ 0 = V_0 < V_1 < V_2 < \ldots < \CU = \colim_k V_k \]
where $V_k=\bigoplus_{j<k}L_{\al_j}$.  If we put $e_k=x_{V_k}$ we find
that $\{e_i\st 0\leq i<k\}$ is a universal basis for $E^0PV_k$, and it
follows by an evident limiting argument that $\{e_i\st i\geq 0\}$ is a
universal topological basis for $E^0P\CU$, giving an isomorphism
$F(P\CU_+,E)=\prod_kE$.  If we put $y=x_{\C[A]}=x_{V_n}=e_n$, it is
easy to see that $e_{ni+j}=y^ie_j$, and it follows that $E^0P\CU$ is a
free module over $E^0\psb{y}$ with basis $\{e_i\st i<n\}$.  Thus,
conditions~(a) and~(b) in Definition~\ref{defn-multicurve} are
satisfied.

Next, we have 
\[ F(P\CU^2_+,E) = F(P\CU_+,F(P\CU_+,E)) = 
    F(P\CU_+,\prod_j E) = \prod_{i,j} E. 
\]
By working through the definitions, we deduce that the elements
$e_i\ot e_j$ form a universal topological basis for 
$E^0(P\CU\tm P\CU)$, so $E^0(P\CU\tm P\CU)=R\hot R$, so 
$\spf(E^0(P\CU\tm P\CU))=C\tm_SC$.  As $P\CU$ is a
commutative group up to equivariant homotopy, we now see that $C$ is a
commutative formal group scheme over $S$.

Now note that $e_1$ is just the coordinate $x$, and this divides $e_k$
for all $k>0$.  In particular it divides $y$, which is a regular
element in $R$, so $x$ is also a regular element.  It is also now easy
to see $x$ generates the ideal $E^0(P\CU,PL_0)$, which is just the
augmentation ideal in the Hopf algebra $R$, so the vanishing locus of
$x$ is the zero-section in $C$.  Thus $x$ is a coordinate on $C$,
showing (via Proposition~\ref{prop-group-coord}) that $C$ is in fact a
formal multicurve group.

Next, recall that $\pi_0((P\CU)^A)=A^*$, which gives a map
$A^*\xra{}P\CU$ of groups up to homotopy, and thus a map 
$\phi\:A^*\tm S\xra{}C$ of formal group schemes over $S$.  By working
through the definitions, we see that the image of the section
$\phi(\al)$ is the closed subscheme $\spec(E^0PL_\al)=\spec(R/x_\al)$,
so the divisor $D:=\sum_\al[\phi(\al)]$ is
\[ \spec(R/\prod_\al x_\al) = \spec(R/y) = E^0P\C[A]. \]
As $y$ is topologically nilpotent, we see that any function on $C$
that vanishes on $D$ is topologically nilpotent, so $C$ is a formal
neighbourhood of $D$.  We have thus proved the following result:

\begin{theorem}\label{thm-coho-efg}
 Let $E$ be a periodically orientable $A$-equivariant ring spectrum.
 Then the scheme $C:=\spf(E^0P\CU)$ is an $A$-equivariant formal group
 over $S:=\spec(E^0)$.  \qed
\end{theorem}

\begin{remark}\label{rem-om-top}
 We have $I_0=\{f\in\OC\st f(0)=0\}=E^0(P\CU,P\C)$, and thus
 $I_0^2=E^0(P\CU,P(\C\op\C))$, and thus 
 \[ \om = I_0/I_0^2 = E^0(P(\C\op\C),P\C) = \tE^0S^2 = \pi_2E. \]
\end{remark}

It will be helpful to record the naturality properties of the above
construction.

\begin{definition}\label{defn-gamma-functor}
 We write $\PP\CO_A$ for the category whose objects are periodically
 orientable $A$-equivariant ring spectra, and whose morphisms are
 homotopy classes of $A$-equivariant ring maps.  We also write
 $\CE\CO_A$ for the evenly orientable subcategory.

 Next, we write $\CG_A$ for the category of triples $(S,C,\phi)$,
 where $S$ is an affine scheme and $(C,\phi)$ is an $A$-EFG over $S$.
 The morphisms from $(S,C,\phi)$ to $(S',C',\phi')$ are the pairs
 $(f,\tf)$ where 
 \begin{itemize}
  \item $f$ is a map $S\to S'$ of schemes;
  \item $\tf$ is an isomorphism $C\to f^*C'$ of formal group schemes
   over $S$;
  \item $\phi'=\tf\circ\phi$.
 \end{itemize}

 The construction $E\mapsto(\spec(E^0),\spf(E^0(P\CU_A)),\phi)$ then
 defines a contravariant functor $\Gm\:\PP\CO_A\to\CG_A$.
\end{definition}

\section{Simple examples}
\label{sec-simple}

Let $\hC$ be a nonequivariant formal group over a scheme $S$, so $\hC$
is the formal neighbourhood of its zero section.  For any finite
abelian group $A$, we can of course let $\phi\:A^*\tm S\xra{}\hC$ be the
zero map, and this gives us an $A$-equivariant formal group.  More
generally, \emph{any} homomorphism $A^*\tm S\xra{}\hC$ will give an
$A$-efg, although often there will not be any homomorphisms other than
zero. 

Now suppose that $\hC$ is the formal group associated to a
nonequivariant even periodic ring spectrum $\hE$.  We then have an
$A$-equivariant ring spectrum $E=F(A_+,\hE)$ (which the Wirthm\"uller
isomorphism also identifies with $A_+\Smash\hE$).  This satisfies
$E^*X=\hE^*\res(X)$, where $\res\:\CS_A\xra{}\CS_0$ is the restriction
functor.  It follows easily that $E$ is periodically orientable, and
that the associated equivariant formal group is just $\hC$, equipped
with the zero map $\phi\:A^*\tm S\xra{}\hC$ as above.

For a slightly more subtle construction, suppose we allow $S$ to be a
formal scheme, and assume that some prime $p$ is topologically
nilpotent in $\OS$.  Suppose also that the formal group $\hC$ has
finite height $n$.  Put $S'=\Hom(A^*,\hC)$; it is well-known that
$\O_{S'}$ is a free module of rank $|A_{(p)}|^n$ over $\OS$, so $S'$
is finite and flat over $S$.  By definition, $S'$ is the universal
example of a formal scheme $T$ over $S$ equipped with a homomorphism
from $A^*$ to the group of maps $T\xra{}\hC$ of formal schemes over
$S$, or equivalently the group of sections of $T\tm_SC$ over $T$.  If
we put $C'=S'\tm_S\hC$, there is thus a tautological map $\phi\:A^*\tm
S'\xra{}C'$.  Here $C'$ is an ordinary formal group over $S'$ and thus
is the formal neighbourhood of its zero section.  It follows that
$(C',\phi)$ is automatically an $A$-equivariant formal group over
$S'$.

Now suppose we have a $K(n)$-local even periodic ring spectrum $\hE$.
We give the ring $\pi_0\hE$ the natural topology as
in~\cite{host:mkl}*{Section 11} --- in most cases of interest, this is
the same as the $I_n$-adic topology.  We then put $S=\spf(\pi_0\hE)$
and $\hC=\spf(\hE^0\CPi)$, which gives an ordinary formal group of
height $n$ over $S$.  Let $EA$ denote a contractible space with free
$A$-action, and put $E=F(EA_+,\hE)$.  This is a commutative
$A$-equivariant ring spectrum, with $E^*X=\hE^*X_{hA}$, where $X_{hA}$
denotes the homotopy orbit space or Borel construction.  In
particular, we have $E^0(\text{point})=\hE^0BA$.  A character 
$\al\in A^*$ gives a map $\spf((B\al)^*)\:\spf(\hE^0BA)\to\hC$, and by
letting $\al$ vary we get a map $\spf(\hE^0BA)\to\Hom(A^*,C)$.  By
reduction to the cyclic case one can check that this is an
isomorphism; see~\cite{hokura:ggc}*{Proposition 5.12} for details.

Next, observe that we have an $A$-equivariant inclusion
$P\CU[0]\xra{}P\CU$, which is nonequivariantly a homotopy equivalence,
so the map $EA\tm P\CU[0]\xra{}EA\tm P\CU$ is an equivariant homotopy
equivalence.  It follows that
$E^*P\CU=E^*P\CU[0]=\hE^*(BA\tm\CPi)=\hE^*BA\ot_{\hE^*}\hE^*\CPi$, and
thus that $\spf(E^0P\CU)=\Hom(A^*,\hC)\tm_S\hC$.  This shows that the
equivariant formal group associated to $E$ is just the pullback
$C'=S'\tm_S\hC$ as discussed above.

\section{Formal groups from algebraic groups}
\label{sec-algebraic}

We now show how to pass from algebraic groups (in particular, elliptic
curves or the multiplicative group) to equivariant formal groups.

\subsection{The multiplicative group}
\label{subsec-MG}

Let $S=\spec(k)$ be a scheme, and consider the group scheme 
$\MG\tm S=\spec(k[u,u^{-1}])$ over $S$.  Suppose we are given a
homomorphism $\phi$ from $A^*\tm S$ to $\MG\tm S$ of group schemes
over $S$, or equivalently a homomorphism $\phi\:A^*\xra{}k^\tm$ of
abstract groups.  We can then form the divisor 
\[ D = \sum_\al[\phi(\al)] = \spec(k[u^{\pm 1}]/y), \]
where $y=\prod_\al(1-u/\phi(\al))$.  It is convenient to observe that
$u$ is invertible in $k[u]/y$ and thus in $k[u]/y^m$ for all $m$, so
$D$ can also be described as $\spec(k[u]/y)$.  We then define $C$ to
be the formal neighbourhood of $D$ in $\MG\tm S$, so
\[ C = \colim_m\, \spec(k[u]/y^m) = \spf(k[u]^\wedge_y), \] 
which is an embeddable formal multicurve.  It is easy to see that this
is a subgroup of $\MG\tm S$ and is an equivariant formal group, with
coordinate $x=1-u$.

The universal example of a ring with a map $A^*\xra{}k^\tm$ is
$k=\Z[A^*]$, which can be identified with the representation ring
$R(A)$.  Thus, the universal example of a scheme $S$ with a map
$A^*\tm S\xra{}\MG\tm S$ as above is $S=\Hom(A^*,\MG)=\spec(R(A))$.
We can apply the above construction in this tautological case to get
an equivariant formal group $C$ over $\Hom(A^*,\MG)$.  Explicitly, if
we let $v_\al\in\Z[A^*]$ be the basis element corresponding to 
$\al\in A^*$ and put $y=\prod_\al(1-uv_{-\al})\in \Z[A^*][u]$, then
$C=\spf(\Z[A^*][u]^\wedge_y)$.

\begin{theorem}[Cole-Greenlees-Kriz]\label{thm-KA}
 The $A$-efg associated to the equivariant complex $K$-theory spectrum
 $K_A$ is isomorphic to the $A$-efg $C$ over $\Hom(A^*,\MG)$
 constructed above.
\end{theorem}
\begin{proof}
 This is just a geometric restatement of~\cite{cogrkr:efg}*{Section 6}.
 It is proved by identifying $K_A^*P\CU$ with $K_{A\tm S^1}^*E\CF_+$
 (where $\CF=\{B\leq A\tm S^1\st B\cap S^1=\{1\}\}$ as in
 Proposition~\ref{prop-SU-EF}) and applying a
 suitable completion theorem.
\end{proof}

\subsection{Elliptic curves}
\label{subsec-elliptic}

We now carry out the same program with the multiplicative group
replaced by an elliptic curve (with some technical conditions assumed
for simplicity).  The resulting equivariant formal groups should be
associated to equivariant versions of elliptic cohomology.  See
Definition~\ref{defn-elliptic-spectrum} and subsequent comments for
more detail.

Suppose that we are given a ring $k$ and an element
$\lm\in k$, and that $2$, $\lm$ and $1-\lm$ are invertible in $k$.
Let $\tC$ be the elliptic curve given by the homogeneous cubic
$y^3=x(x-z)(x-\lm z)$, so the zero element is $O=[0:1:0]$, and the
points $P:=[0:0:1]$, $Q:=[1:0:1]$ and $R:=[\lm:0:1]$ are the three
points of exact order two in $\tC$.  Define rational functions $t$ and
$r$ on $\tC$ by $t([x:y:z])=x/y$ and $r([x:y:z])=z/y$.  One checks
that the subscheme $U=\tC\sm\{P,Q,R\}$ is the affine curve with
equation $r=t(t-r)(t-\lm r)$, and that on $U$, the function $t$ has a
simple zero at $O$ and no other poles or zeros.

Now let $A$ be an abelian group of odd order $n$, and let
$\phi\:A^*\tm S\xra{}\tC$ be a homomorphism.  Define
$V=\bigcap_\al(U+\phi(\al))$, which is an affine open subscheme of
$U$.
\begin{lemma}
 For each $\bt\in A^*$, the section $\phi(\bt)\:S\xra{}\tC$ actually
 lands in $V$.
\end{lemma}
\begin{proof}
 We first show that for all $\gm\in A^*$, the section $\phi(\gm)$
 lands in $U$.  Put $D=[P]+[Q]+[R]$, so $U=\tC\sm D$.  Let $T$ be the
 closed subscheme of points $s\in S$ where $\phi(\gm)(s)\in D$; we
 must show that $T=\emptyset$.  As $n$ is odd and $D$ is the divisor
 of points of exact order $2$, we see that multiplication by $n$ is
 the identity on $D$, but of course $n.\phi(\al)=O$.  We conclude that
 over $T$ we have $O\in D$.  As $2$ is invertible in $k$ we know that
 $O$ and $D$ are disjoint, so $T=\emptyset$ as required.  

 We now apply this to $\gm=\bt-\al$ to deduce that
 $\phi(\bt)\in U+\phi(\al)$.  This holds for all $\al$, so
 $\phi(\bt)\in V$ as claimed.
\end{proof}

We now define $C$ to be the formal neighbourhood of the divisor
$D=\sum_\al[\phi(\al)]$ in $V$.  If we put
$s(a)=\prod_\al t(a-\phi(\al))$ then $s\in\OV$ and the vanishing locus
of $s$ is just $D$, so we have $\OC=(\OV)^\wedge_s$.  Using this, we
see that $C$ is an equivariant formal group, with coordinate $t$ and
good parameter $s$.

Now suppose instead that we are given a curve $\tC$ over $S$ as above,
but not the map $\phi\:A^*\tm S\xra{}\tC$.  We can then consider the
scheme $S_1=\Hom(A^*,\tC)$, which is easily seen to be a closed
subscheme of $\Map(A^*,U)$ and thus affine.  We can thus pull back
$\tC$ to get a curve $\tC_1$ over $S_1$ equipped with a tautological
map $\phi\:A^*\tm S_1\xra{}\tC_1$, and we can carry out the previous
construction to get an equivariant formal group $C_1$ over $S_1$.

\begin{definition}\label{defn-elliptic-spectrum}
 An $A$-equivariant elliptic spectrum consists of an evenly orientable
 equivariant ring spectrum $E$ together with an elliptic curve $\tC$
 over an affine scheme $S$ and a compatible system of
 isomorphisms $\spec(E^0(A/B))=\Hom(B^*,\tC)$ and $\spf(E^0P\CU_A)=C$
 (where $C$ is constructed from $\tC$ as above).  
\end{definition}

We will not take the trouble to make this definition more precise, as
we will not use it very seriously.

In Section~\ref{sec-rational} we will construct equivariant elliptic
spectra associated to elliptic curves over $\Q$-algebras, as a simple
application of the general theory of rational $A$-spectra.

Elsewhere we have obtained partial results about the existence of
integral equivariant elliptic spectra.  These have unsatisfactory
indeterminacy and awkward technical hypotheses, but nonetheless they
are sufficient to make it clear that
Definition~\ref{defn-elliptic-spectrum} is the right one.  Much better
results (relying on a far-reaching theory of equivariant $E_\infty$
ring spectra and derived algebraic geometry) have been announced by
Jacob Lurie, but details have yet to appear.

\section{Equivariant formal groups of product type}
\label{sec-product}

A simple class of $A$-EFGs can be constructed as follows.  

\begin{definition}\label{defn-pt}
 Let $\hC$ be an ordinary, nonequivariant formal group, and let $B$ be
 a subgroup of $A$.  We then have a formal multicurve $C:=B^*\tm\hC$
 and a homomorphism
 \[ \phi := (A^* \xra{\text{res}} B^* \xra{\text{inc}} B^*\tm\hC=C), \]
 giving an $A$-efg.  Equivariant formal groups of this kind are said to
 be of \emph{product type}.  We call $B$ the \emph{core} of
 $(C,S,\phi)$.   
\end{definition}
We will show that EFGs over fields are of product type, and EFGs over
$\Q$-algebras are locally of product type.  Moreover, we will
introduce equivariant analogues of the Morava $K$-theory spectra, and
show that the associated EFGs are of product type.

\begin{proposition}\label{prop-product}
 An $A$-efg $(C,\phi)$ is of product type iff for every character
 $\al\in A^*$ with $\phi(\al)\neq 0$ in $C$ (or equivalently,
 $x(\phi(\al))\neq 0$ in $\OS$), the element $x(\phi(\al))$ is
 invertible in $\OS$.  (This is easily seen to be independent of the
 choice of coordinate.)
\end{proposition}
\begin{proof}
 First suppose that for all $\al$ with $\phi(\al)\neq 0$, the element
 $x(\phi(\al))$ is invertible.  The kernel of $\phi$ is a subgroup of
 $A^*$, so it necessarily has the form $\ann(B)$ for some $B\leq A$,
 so $\phi$ factors as $A^*\xra{\text{res}}B^*\xra{\psi}C$ for some
 $\psi$.  By assumption, $x(\psi(\bt))$ is invertible for all
 $\bt\in B^*\sm\{0\}$.

 Let $\hC=\{a\in C\st x(a)\text{ is nilpotent }\}$ be the formal
 neighbourhood of $0$ in $C$, and define $\sg\:B^*\tm\hC\xra{}C$ by
 $\sg(\bt,a)=\psi(\bt)+a$.  We need to show that $\sg$ is an
 isomorphism.  For this, we define $x_\bt(a)=x(a-\psi(\bt))$ and
 $y=\prod_{\bt\in B^*}x_\bt$ and $R=\OC$.  From the definition of an
 equivariant formal group, we know that $R=R^\wedge_y$, and it is
 clear that
 \[ \O_{B^*\tm\hC}=\prod_\bt R^\wedge_{x_\bt}. \]
 It will thus suffice to show that the natural map 
 \[ R^\wedge_y \xra{} \prod_\bt R^\wedge_{x_\bt} \] 
 is an isomorphism.  This will follow from the Chinese Remainder
 Theorem if we can check that the ideal $(x_\bt(a),x_\gm(a))$ contains
 $1$ whenever $\bt\neq\gm$.  This is clear because modulo that ideal,
 we have $\psi(\bt)=a=\psi(\gm)$, so $\psi(\bt-\gm)=0$, so
 $x(\psi(\bt-\gm))=0$, but $x(\psi(\bt-\gm))$ is invertible by
 assumption.  Thus, $C$ is of product type, as claimed.

 Conversely, suppose that $C$ is of product type.  The vanishing locus
 of $x$ is contained in $\{0\}\tm\hC$, so $x$ must be invertible on
 $(B^*\sm\{0\})\tm\hC$.  It follows immediately that when
 $\phi(\al)\neq 0$ we have $\phi(\al)\in(B^*\sm\{0\})\tm\hC$ and so
 $x(\phi(\al))$ is invertible, as required.
\end{proof}

\begin{corollary}\label{cor-product-field}
 Every $A$-equivariant formal group over a field is of product type.
\end{corollary}
\begin{proof}
 This is immediate from the proposition.
\end{proof}

We next show how groups of product type occur in topology.  For this
we need to use the geometric fixed point functors
$\phb^B\:\CS_A\xra{}\CS_0$ for $B\leq A$.  These preserve smash
products and satisfy $\phb^B\Sgi X=\Sgi X^B$ for based $A$-spaces
$X$.  (The definition and further properties of these functors will be
recalled in Section~\ref{sec-pushout}.)

\begin{theorem}
 Let $\hK$ be a nonequivariant even periodic cohomology theory, with
 associated formal group $\hC$ over $S$, and let $B$ be a subgroup of
 $A$.  Define a cohomology theory $K^*$ on $\CS_A$ by
 $K^*X=\hK^*\phb^BX$.  Then $K$ is evenly orientable, and the
 associated equivariant formal group is just $B^*\tm\hC$ over $S$.
\end{theorem}

\begin{proof}
 Note that $\phb^B S^V=S^{V^B}$ for any virtual complex representation
 $V$, and that $\phb^B\Sgi X=\Sgi X^B$ for any based $A$-space $X$.
 It follows that $\pi_1K=\pi_1\hK=0$ and that the periodicity
 isomorphism $F(S^{2n},\hK)=\hK$ gives an isomorphism 
 \[ K^*(X_+\Smash S^V) = \hK^*(X^B_+\Smash S^{V^B}) =
     \hK^*X^B_+ = K^*X_+
 \]
 of modules over $K^*X_+$.  This implies that $K$ is evenly periodic,
 with $K^0(\text{point})=\hK^0(\text{point})$ and thus
 $\spec(K^0(\text{point}))$ is the base scheme $S$ for $\hC$.  We also
 have 
 \[ K^0P\CU = \hK^0(P\CU)^B = \hK^*(B^*\tm\CPi) = \O_{B^*\tm\hC},  \]
 so the equivariant formal group associated to $K$ is just $B^*\tm\hC$
 as claimed.
\end{proof}

\begin{example}\label{eg-morava-k}
 Let $\hK=\hK(p,n)$ be the two-periodic version of Morava $K$-theory
 at a prime $p$, with height $n$.  We define an equivariant theory
 $K=K(p,n,B)$ as above; this is called \emph{equivariant Morava
 $K$-theory}.  In~\cite{st:ebc} we present evidence that these
 theories deserve this name, because they play the expected r\^ole in
 equivariant analogues of the Hopkins-Devinatz-Smith nilpotence
 theorems, among other things.  The same paper also explains the
 representing object for the theory $K$, and shows that we have
 natural isomorphisms as follows:
 \begin{align*}
  K_*(X\Smash Y) &= K_*(X)\ot_{K_*}K_*(Y) \\
  K^*X           &= \Hom_{K_*}(K_*X,K_*).
 \end{align*}
\end{example}

We now give a slight generalisation.
\begin{definition}\label{defn-lpt}
 We say that $(C,S,\phi)$ is \emph{locally of product type} if there is
 a splitting $S=\coprod_{B\leq A}S_B$ such that the restricted group
 $C_B=C\tm_SS_B$ is of product type with core $B$.
\end{definition}

\begin{definition}\label{defn-split-element}
 Let $R$ be a ring, and $r$ and element of $R$.  We say that $r$ is
 \emph{split} if there is an idempotent $e\in R$ with $Rr=Re$, or
 equivalently there is a splitting $R=R_0\tm R_1$ with respect to
 which $r\in\{0\}\tm R_1^\tm$.  
\end{definition}

\begin{definition}\label{defn-SB}
 Let $(C,S,\phi)$ be an $A$-equivariant formal group, and let $x$ be a
 coordinate on $C$.  For each $B\leq A$ we put 
 \[ I_B = (x(\phi(\al))\st \al\in\ann(B)\leq A^*) \leq \OS. \]
 We also put $S[B]=\spec(\OS/I_B)$, which can be described more
 invariantly as the largest closed subscheme $T\sse S$ over which
 $\phi\:A^*\tm T\to C\tm_ST$ factors through $B^*\tm T$.  
\end{definition}

\begin{proposition}\label{prop-lpt}
 Let $(C,S,\phi,x)$ be as above.  Then the following are equivalent:
 \begin{itemize}
  \item[(a)] $C$ is locally of product type.
  \item[(b)] For each $B\leq A$, the subscheme $S[B]\sse S$ is open as well
   as closed.
  \item[(c)] For each $\al\in A^*$ the element
   $e_\al=x(-\phi(\al))\in\OS$ is split. 
 \end{itemize}
\end{proposition}
\begin{proof}
 First suppose that~(a) holds, giving a decomposition
 $S=\coprod_BS_B$.  It follows that each $S_B$ is open and closed,
 and $S[B]=\coprod_{B'\leq B}S_{B'}$, which proves~(b).

 Next, suppose that~(b) holds.  For $\al\in A^*$ one checks that
 $S[\ker(\al)]$ is just the largest closed subscheme of $S$ where
 $\phi(\al)$ vanishes, or equivalently where $e_\al$ vanishes.
 As this is open as well as closed, it is the vanishing locus of an
 idempotent, say $f_\al$.  It follows that $e_\al$ and $f_\al$
 generate the same ideal, so $e_\al$ is split.  This proves~(c).

 Finally, suppose that~(c) holds.  This means that for each $\al$ we
 can split $S$ as $E_\al\amalg F_\al$ with $e_\al$ invertible on
 $E_\al$ and zero on $F_\al$.  Next, for any subset $U\sse A^*$ we put
 \[ F_U =
     \bigcap_{\al\in U} F_\al \cap \bigcap_{\al\not\in U} E_\al.
 \]
 We find that $F_U$ is both open and closed in $S$, and that $S$ is
 the disjoint union of all the sets $F_U$. 

 Now suppose that $U$ is not a subgroup of $A^*$, so there must exist
 $\al,\bt\in U$ with $\gm=\al+\bt\not\in U$.  Then we have $F_U\sse
 F_\al\cap F_\bt\cap E_\gm$.  This means that over $F_U\sse S$ we have
 $x(\phi(\al))=x(\phi(\bt))=0$ but $x(\phi(\gm))$ is invertible.  As
 $x$ is a coordinate we have $\phi(\al)=\phi(\bt)=0$, which implies
 $\phi(\gm)=0$, so $x(\phi(\gm))$ is zero as well as invertible.  It
 follows that $F_U$ is the empty scheme.

 Suppose instead that $U$ is a subgroup of $A^*$, so $U=\ann(B)$ for
 some subgroup $B\leq A$.  We then define $S_B$ to be $F_U$, and
 observe that $S=\coprod_BS_B$, and that the restriction
 $C_B=C\tm_SS_B$ is of product type with core $B$.  This proves~(a).
\end{proof}

\section{Equivariant formal groups over rational rings}
\label{sec-rational}

We next prove an equivariant analogue of the well-known fact that all
formal groups over a $\Q$-algebra are additive.  We write $\hGa$ for
the ordinary additive formal group over $S$.  If we consider formal
schemes over $S$ as functors in the usual way, this sends an
$\OS$-algebra $R$ to the set $\Nil(R)$ of nilpotents in $R$.  Given a
free module $L$ of rank one over $\OS$ (or equivalently, a
trivialisable line bundle over $S$), we can instead consider the
functor $R\mapsto L\ot_{\OS}\Nil(R)$, which we denote by $L\ot\hGa$.
This gives a formal group over $S$, noncanonically isomorphic to
$\hGa$.  If $C$ is a formal multicurve group over $S$, then the
cotangent spaces to the fibres give a trivialisable line bundle
$\om_C$ on $S$.  This is easily seen to be the same as $\om_{\hC}$,
where $\hC$ is the formal neighbourhood of zero, as usual.  From now
on we just write $\om$ for this module.  If $S$ lies over $\spec(\Q)$
then the theory of logarithms for ordinary formal groups gives a
canonical isomorphism $\hC\xra{}\om^{-1}\ot\hGa$.

\begin{theorem}\label{thm-product-Q}
 Let $(C,\phi)$ be an $A$-equivariant formal group over a scheme $S$,
 such that the integer $n=|A|$ is invertible in $\OS$.  Then $C$ is of
 local product type.  Moreover, if $\OS$ is a $\Q$-algebra then
 $\hC\simeq\om_C^{-1}\ot\hGa$ and so
 \[ C \simeq \coprod S_B\tm_S B^*\tm (\om_C^{-1}\ot\hGa). \]
\end{theorem}
\begin{proof}
 Put $n=|A|$, and choose a coordinate $x$ on $C$.  For formal reasons
 we have $x(a+b)=x(a)+x(b)\pmod{x(a)x(b)}$ as functions on $C^2$, and
 it follows that $x(na)=v_n(a)x(a)$ for some function $v_n$ on $C$
 with $v_n(0)=n$.  (These functions will be studied in greater detail
 in Section~\ref{sec-transfer}.)  Put $Q=\OC/(x.v_n)$, so $\spf(Q)$ is
 just the closed subscheme $C[n]$ of points of order $n$ in $C$.  In
 $Q$ we have $xv_n=0$, and for any function $f\in\OC$ we have
 $f-f(0)\in(x)$ so $f.v_n=f(0).v_n$.  In particular, we can take
 $f=v_n$ to see that $v_n^2=nv_n$, so the function $u=1-v_n/n$ is
 idempotent in $Q$.  We have $u(0)=0$ so $u\in Qx$, but also $v_nx=0$
 so $ux=x$ so $x\in Qu$.  For any $\al\in A^*$ we know that
 $n\phi(\al)=0$ so $\phi(\al)\:S\to C[n]$, and we find that
 $x(\phi(\al))$ generates the same ideal in $\OS$ as $u(\phi(\al))$,
 so $x(\phi(\al))$ is split (in the sense of
 Definition~\ref{defn-split-element}).  It follows by
 Proposition~\ref{prop-lpt} that $C$ is of local product type.

 The rational statement now follows from the nonequivariant theory.
\end{proof}

The following slight extension can easily be proved in the same way.
\begin{proposition}
 Let $(C,\phi)$ be an $A$-equivariant formal group over a scheme $S$,
 such that $\OS$ is an algebra over $\Zpl$.  There is of course a
 unique splitting $A=A_0\tm A_1$, where $A_0$ is a $p$-group and $p$
 does not divide $|A_1|$.  Let $C_0\sse C$ be the formal neighbourhood
 of $[\phi(A_0^*)]$, and let $\phi_0\:A_0^*\xra{}C_0$ be the
 restriction of $\phi$.  Then there is a canonical decomposition
 $S=\coprod_{B\leq A_1}S_B$, and a corresponding decomposition
 \[ C \simeq \coprod S_B \tm_S C_0 \tm B^*, \]
 such that over $S_B$, the map $\phi$ is the product of $\phi_0$ and
 the restriction map $A_1^*\xra{}B^*$. \qed
\end{proposition}

We would like to understand how the splitting in
Theorem~\ref{thm-product-Q} works out when the equivariant formal
group comes from a ring spectrum.  For simplicity we will treat only
the rational case, although many parts of our analysis can also be
made to work assuming only that $|A|$ is invertible.  We start by
recalling an algebraic description of the category $\Q\CS_A$ of
rational $A$-spectra.  Let $\CV_*[A]$ denote the category of graded
modules over the group ring $\Q[A]$.  If $X\in\Q\CS_A$ and $B\leq A$
then the nonequivariant spectrum $\phb^B(X)$ has a homotopical action
of the Weyl group $A/B$, so $\pi_*(\phb^B(X))\in\CV_*[A/B]$.  Put
$\CA_A=\prod_{B\leq A}\CV_*[A/B]$, and define $\Phi\:\Q\CS_A\to\CA_A$
by
\[ \Phi(X)_B = \pi_*(\phb^B(X)). \]
We can make $\CA_A$ into a symmetric monoidal category by the rule
$(M\ot N)_B=M_B\ot N_B$.  For nonequivariant rational spectra $X$ and
$Y$ we have 
\[ \pi_*(X\Smash Y) = H_*(X\Smash Y;\Q) = 
    H_*(X;\Q)\ot H_*(Y;\Q) = \pi_*(X)\ot\pi_*(Y).
\]
Using this, we we that $\Phi$ takes smash products to tensor products.

\begin{theorem}\label{thm-rational-split}
 The functor $\Phi\:\Q\CS_A\to\CA$ is an equivalence.
\end{theorem}
As far as we know, the literature only contains a rather indirect
proof of this fact, going via the theory of Mackey functors.  Here we
give a slightly more direct argument.
\begin{proof}
 All spectra in this proof are implicitly rationalised; we will not
 indicate this explicitly in the notation.

 We have a natural map
 \[ \Phi_{XY}\:\Q\CS_A(X,Y)_*\to\CA_A(\Phi(X),\Phi(Y))_*, \]
 which we claim is an isomorphism.  

 Using Maaschke's theorem we see that all objects in $\CA$ are both
 injective and projective.  It is also clear that $\Phi$ preserves all
 coproducts.  It follows that both $\Q\CS_A(X,Y)_*$ and
 $\CA(\Phi(X),\Phi(Y))_*$ are cohomology theories of $X$.  If $X$ is
 finite then they are also homology theories in $Y$.  We can thus
 reduce easily to the case where $X$ and $Y$ are both finite.

 Now let $DY$ be the Spanier-Whitehead dual of $Y$.  The duality
 between $Y$ and $DY$ is encoded by the unit map
 $\eta\:S\to DY\Smash Y$ and the counit $\ep\:Y\Smash DY\to S$, which
 make the following diagrams commute:
 \[ \xymatrix{
     DY \rto^(0.3){\eta\Smash 1} \drto_1 &
     DY\Smash Y \Smash DY  \dto^{1\Smash\ep} & &
     Y \rto^(0.3){1\Smash\eta} \drto_1 &
     Y\Smash DY\Smash Y \dto^{\ep\Smash 1} \\
     & DY & & & Y
    }
 \]
 We can apply $\Phi$ to this to get a perfect duality between
 $\Phi(Y)$ and $\Phi(DY)$, which in turn gives a commutative square 
 \[ \xymatrix@C=4cm{
     \Q\CS_A(X,Y) \rto^{\Phi_{XY}} \dto_{\simeq} &
     \CA_A(\Phi(X),\Phi(Y)) \dto^{\simeq} \\
     \Q\CS_A(X\Smash DY,S^0) \rto_{\Phi_{X\Smash DY,S^0}} &
     \CA_A(\Phi(X\Smash DY),\Phi(S^0)).
    }
 \]
 Using this we reduce to the case $Y=S^0$.  We can then reduce further
 by cellular induction to the case $X=A/B_+$.  Here
 $\Q\CS_A(A/B_+,S^0)_0$ is the Burnside ring of $B$, which we will
 call $\Om(B)$.  Using the tom Dieck splitting
 $(\Sgi_A S^0)^B=\bigWedge_{C\leq B}\Sgi B(B/C)_+$ we also see that
 $\Q\CS_A(A/B_+,S^0)_n=H_n((\Sgi_AS^0)^B;\Q)=0$ for $n\neq 0$.  We
 also find that $\CA_A(\Phi(A/B_+),\Phi(S^0))=\prod_{C\leq B}\Q$, and
 the map $\Phi_{A/B_+,S^0}$ is just the usual fixed-point map, which
 is well-known to be an isomorphism.

 We now know that the functor $\Phi$ is full and faithful.  As
 $\Q\CS_A$ has all coproducts, and splittings for idempotents, we see
 that the essential image of $\Phi$ is closed under coproducts and
 retracts.  Moreover, all monomorphisms and epimorphisms in $\CA_A$
 are split, so $\img(\Phi)$ is closed under taking subobjects and
 quotient objects.  If $M\in\CA_A$ and $m\in M_{B,d}$ then there is a
 map $f_m\:\Phi(\Sg^dA/B_+)\to M$ sending an evident generator to
 $m$.  By taking a large direct sum of maps like this, we can express
 $M$ as a quotient of an object in the image of $\Phi$.  It follows
 that $\Phi$ is essentially surjective, as claimed.
\end{proof}

\begin{remark}\label{rem-CA-maps}
 The group $Y^mX=[X,Y]_{-m}=\CA_A(\Phi(X),\Phi(Y))_{-m}$ can also be
 described as
 $\prod_B\prod_iH^{m+i}(\phb^B(X);\pi_i(\phb^B(Y)))^{A/B}$.
\end{remark}

\begin{remark}\label{rem-torus-action}
 The action of $A/B$ on $\pi_*(\phb^B(X))$ plays a central r\^ole
 here, so it is useful to know when this action is trivial.  Choose an
 embedding of $A$ in a torus $T$, and suppose that $X$ arises by
 applying the restriction functor $\res^T_A$ to a $T$-spectrum which
 we also call $X$.  Then the action map $A/B\to[\phb^B(X),\phb^B(X)]$
 factors through $\pi_0(T/B)=1$, so $A/B$ acts trivially on all
 homotopical invariants of $\phb^B(X)$.  

 In particular, elementary character theory tells us that any
 complex representation $V$ of $A$ admits a compatible $T$-action, so
 the above remarks apply to all spaces constructed functorially from
 $V$ such as the projective space $PV$, the sphere $S(V)$ and so on.
\end{remark}

\begin{proposition}\label{prop-QEO-split}
 Let $E$ be an evenly orientable rational $A$-spectrum, let
 $(C,S,\phi)$ be the corresponding formal group, and let $S_B\sse S$
 be as in Theorem~\ref{thm-product-Q}.  Put $j_B=\pi_0(\phb^BE)$ and
 $k_B=j_B^{A/B}$ and $\om=\pi_2E$.  Then $S_B=\spec(k_B)$, and for all
 $X$ we have 
 \[ E^m(X) =
     \prod_{B\leq A}\prod_{d\in\Z}
       H^{m+2d}(\phb^B(X);j_B\ot\om^d)^{A/B}.
 \]
\end{proposition}
\begin{proof}
 Using Theorem~\ref{thm-rational-split} and Remark~\ref{rem-CA-maps}
 we have 
 \[ E^m(X) =
     \prod_B\prod_i H^{m+i}(\phb^B(X);\pi_i(\phb^B(E)))^{A/B}.
 \]
 We now take $m=1$ and $X=(A/B')_+$.  As $E$ is evenly periodic, the
 left hand side vanishes.  On the right, the factor where $B=B'$ is
 \[ H^0(A/B';\pi_{-1}(\phb^{B'}(E)))^{A/B'} = \pi_{-1}(\phb^{B'}(E)); \]
 this must therefore vanish as well.  As $E$ is periodically
 orientable we find that
 $\pi_{i+2}(\phb^B(E))\simeq\pi_i(\phb^B(E))\ot_{E_0}\om$.   It
 follows that $\pi_{2i}(\phb^B(E))=j_B\ot_{E^0}\om^i$ and
 $\pi_{2i+1}(\phb^B(E))=0$.  Using this we get 
 \[ E^m(X) =
     \prod_{B\leq A}\prod_{d\in\Z}
       H^{m+2d}(\phb^B(X);j_B\ot\om^d)^{A/B}
 \]
 as claimed.  We apply this to the spectrum $X=\Sgip P\CU$, noting
 that the action of $A/B$ on $\phb^B(X)=\Sgip(B^*\tm\CPi)$ is
 homotopically trivial by Remark~\ref{rem-torus-action}.  This gives 
 \[ E^0P\CU =
     \prod_B\prod_d\Map(B^*,H^{2d}(\CPi;k_B\ot\om^d))
 \]
 and thus
 \[ C \simeq \coprod_B B^*\tm\spec(k_B)\tm\hGa. \]
 We also find that the map $\phi\:A^*\tm S\to C$ is given over
 $\spec(k_B)$ by the restriction $A^*\to B^*$ and the obvious
 inclusion $B^*\to B^*\tm\hGa$, so $S_B=\spec(k_B)$.
\end{proof}

\begin{construction}\label{cons-Delta}
 Let $(C,S,\phi)$ be an $A$-EFG as in Theorem~\ref{thm-product-Q}.
 Put $k_B=\O_{S_B}$, which is a $\Q$-algebra.  Let $\om$ be the
 cotangent space to $C$ at zero, and put $\om_B=\om|_{S_B}$, which is
 a free module of rank one over $k_B$.  Define a multiplicative
 cohomology theory $E^*$ on $A$-spectra by 
 \[ E^m(X) = \prod_B\prod_d H^{m+2d}(\phb^B(X);\om_B^d)^{A/B}. \]
 (Here we use the trivial action of $A/B$ on $\om_B^d$, so the action
 on the cohomology group comes solely from the action on $\phb^B(X)$.)
 Write $E=\Dl(C,S,\phi)$ for the spectrum that represents this
 cohomology theory.  
\end{construction}

\begin{proposition}
 $\Dl(C,S,\phi)$ is evenly orientable, with
 $\Gm\Dl(C,S,\phi)\simeq(C,S,\phi)$. 
\end{proposition}
\begin{proof}
 Put $E=\Dl(C,S,\phi)$.  From the definitions we have
 $E^0=\prod_BH^0(1;k_B)^{A/B}=\prod_Bk_B=\OS$ and
 \[ E^0P\CU = \prod_B\prod_d H^{2d}(B^*\tm\CPi;\om_B^d)
     = \prod_B\prod_d\prod_{\bt\in B^*} H^{2d}(\CPi;\om^d).
 \] 
 Choose a generator $u$ of $\om$, and let $y$ be the standard
 generator of $H^2(\CPi)$.  Put 
 \[ x_{B,d,\bt} = 
     \begin{cases} 
      1 & \text{ if } d=0,\; \bt\neq 0 \\
      y\ot u & \text { if } d=1,\; \bt=0 \\
      0 & \text{ otherwise. }
    \end{cases}
 \]
 This defines a class $x\in E^0P\CU$.  If we quotient out the ideal
 generated by $x$ then all factors where $\bt\neq 0$ or $d>0$ are
 killed, and this leaves $(E^0P\CU)/x=\prod_Bk_B=E^0$.  

 Now consider a character $\al$ of $A$ and the subspace
 $P(L_0\op L_\al)\subset P\CU$.  If $B\leq\ker(\al)$ we find that
 \[ \phb^BP(L_0\op L_\al) = 
     \{0\}\tm\CP^1\subset B^*\tm\CPi = 
      \phb^BP\CU.
 \]
 On the other hand, if $B\not\leq\ker(\al)$ we find that 
 \[ \phb^BP(L_0\op L_\al) = \{0,\al|_B\} \sse B^* \subset
      B^*\tm\CPi = \phb^BP\CU.
 \]
 Either way, we see that the reduced cohomology of this space with
 coefficients in $k_B^*=\bigoplus_d\om_B^d$ is freely generated over
 $k_B^*$ by the image of $x$.  By taking the product over all $B$, we
 see that the image of $x$ generates $E^*(P(L_0\op L_\al),PL_0)$ as a
 module over $E^*$.  We leave it to the reader to check in the same
 way that $E^*((A/A')\tm P(L_0\op L_\al),A/A'\tm PL_0)$ is freely
 generated by $x$ over $E^*(A/A')$, so $x$ has the required universal
 generating property and is a coordinate for $E$.  This shows that $E$
 is periodically orientable, and it is immediate from the definitions
 that $E^{-1}(A/A')=0$, so in fact it is evenly orientable.

 It is clear by construction that 
 \[ \spf(E^0 P\CU) = \coprod_BB^*\tm S_B \tm (\om_B^{-1}\ot\hGa)
      \simeq C,
 \] 
 so $\Gm\Dl$ is the identity functor.
\end{proof}

\begin{proposition}\label{prop-gsr}
 There is a natural map $\xi_E\:\Dl\Gm(E)\to E$ for $E\in\Q\CE\CO_A$.
 Moreover, the following are equivalent:
 \begin{itemize}
  \item[(a)] $E^0(A/B)=E^0/(x(\phi(\al))\st \al\in(A/B)^*\leq A^*)$
   (for any coordinate $x$ and any subgroup $B\leq A$).
  \item[(b)] The unit map $E^0\to E^0(A/B)$ is surjective (for any
   subgroup $B\leq A$).
  \item[(c)] $A/B$ acts trivially on $E^0(A/B)$ (for any subgroup
   $B\leq A$).
  \item[(d)] $\xi_E$ is an isomorphism.
 \end{itemize}
\end{proposition}
\begin{remark}\label{rem-torus-action-ii}
 If $E$ is the restriction of a $T$-spectrum for some torus
 $T\supseteq A$ then we see as in Remark~\ref{rem-torus-action}
 that~(c) holds, and so the other conditions hold as well.  In
 particular, this will hold automatically for the rationalisations of
 popular examples such as $K$-theory and cobordism, which are
 constructed in a uniform way for all compact Lie groups.
\end{remark}

\begin{proof}
 Let $E$ be a rational, evenly periodic $A$-equivariant ring spectrum.
 Put $j_B=\pi_0(\phb^BE)$ and $k_B=j_B^{A/B}$ and $\om=\pi_2E$ as in
 Proposition~\ref{prop-QEO-split}.   The inclusions $k_B\to j_B$ give
 an inclusion 
 \[ \xi_{E,X} \: \prod_B H^*(\phb^B(X);k_B\ot\om^*)^{A/B} \to 
     \prod_B H^*(\phb^B(X);j_B\ot\om^*)^{A/B}, 
 \]
 or equivalently $\Dl\Gm(E)^*(X)\to E^*(X)$.  By the Yoneda lemma,
 this is induced by a map $\Dl\Gm(E)\to E$ of equivariant ring spectra.

 Next, it is trivial that (a) implies (b), and a naturality argument
 shows that~(b) implies~(c).  Now note that 
 \[ E^0(A/A') = \prod_{B\leq A'} H^0(A/A';j_B)^{A/B} 
              = \prod_{B\leq A'} j_B^{A'/B}, 
 \]
 and the factor indexed by $B=A'$ is just $j_{A'}$.  It follows that
 $A$ acts trivially on $E^0(A/A')$ for all $A'$ iff it acts trivially
 on $j_{A'}$ for all $A'$ iff $k_{A'}=j_{A'}$ for all $A'$.  If this
 holds, it is clear that $\xi_E$ is an isomorphism.  Thus~(c)
 implies~(d).  Finally, suppose that~(d) holds.  Note that
 $\Dl\Gm(E)^0(A/A')=\prod_{B\leq A'}k_B$, so if $\xi_E$ is 
 an isomorphism we must have $k_B=j_B^{A'/B}$ for all
 $B\leq A'\leq A$.  By considering the case $A'=B$ we see that
 $j_B=k_B$, so the action of $A$ on $j_B$ is trivial.  This gives
 $E^0=\prod_Bk_B$ and $E^0(A/A')=\prod_{B\leq A'}k_B$.  We also know
 that the map $\phi\:A^*\tm S\to C$ is given over $S_B=\spec(k_B)$ by
 the restriction $A^*\to B^*$ and the inclusion $B^*\to B^*\tm\hC$,
 so the largest closed subscheme where $\phi((A')^*)=0$ is just
 $\coprod_{B\leq A'}S_B$, which is the spectrum of the ring
 $\prod_{B\leq A'}k_B=E^0(A/A')$.  On the other hand, this closed
 subscheme can also be described as the spectrum of the ring 
 \[ k'_{A'} = E^0/(x(\phi(\al))\st \al\in(A/A')^*\leq A^*). \]
 We must therefore have $E^0(A/A')=k'_{A'}$, so~(a) holds.
\end{proof}

The following corollary follows directly.

\begin{corollary}
 Let $\Q\CE\CO'_A$ be the full subcategory of $\Q\CE\CO_A$ generated
 by objects for which the equivalent conditions of
 Proposition~\ref{prop-gsr} are satisfied.  Then the functors $\Gm$
 and $\Dl$ give an equivalence between $\Q\CE\CO'_A$ and $\Q\CG_A$.
 \qed
\end{corollary}

\begin{example}
 Suppose we have an elliptic curve $\tC$ over a scheme $S=\spec(k)$ as
 in Section~\ref{subsec-elliptic}, and suppose that $k$ is a
 $\Q$-algebra.  We then define an equivariant formal group $C$ over
 $S_1=\Hom(A^*,\tC)$ as described in Section~\ref{subsec-elliptic},
 and use this to construct an evenly periodic $A$-equivariant ring
 spectrum $E=\Dl(C)$.  Using Proposition~\ref{prop-gsr} we see that
 $E^0(A/B)=\O_{\Hom(B^*,\tC)}$ for all $B\leq A$, so we have an
 elliptic spectrum as defined in
 Definition~\ref{defn-elliptic-spectrum}. 
\end{example}

\section{Equivariant formal groups of pushout type}
\label{sec-pushout}

We next consider a slightly different generalization of the notion of
a group of product type.  

\begin{definition}
 Suppose we have a subgroup $B\leq A$ and a formal multicurve group
 $C'$, with a map $\phi'\:(A/B)^*\xra{}C'$ making it an
 $A/B$-equivariant formal group.  There is an evident embedding
 $(A/B)^*\xra{}A^*$, which we can use to form a pushout
 \[ \xymatrix{
  {(A/B)^*} \rto^{\phi'} \ar@{ >->}[d] &  {C'} \dto \\
  {A^*} \rto_{\phi} & {C.}
  } \]
 If we choose a transversal $T$ to $(A/B)^*$ in $A^*$, then the
 underlying scheme of $C$ is just $\coprod_{\al\in T}C$.  This implies
 that the formation of the pushout is compatible with base change, and
 that $C$ is an $A$-equivariant formal group.  Formal groups
 constructed in this way are said to be of \emph{pushout type}.  (The
 case where $\phi'=0$ evidently gives groups of product type.)
\end{definition}

We next examine how formal groups of this kind can arise in
equivariant topology.  For this, we need to recall in more detail the
various different change of group functors and fixed-point functors
for $A$-spectra.

Given a homomorphism $\zt\:B\xra{}A$, there is a pullback functor
$\zt^*\:\CS_A\xra{}\CS_B$, which preserves smash products and function
spectra.  (Note that if $\zt$ is not injective, then $\zt^*\CU_A$ is
not a complete $B$-universe, so the definition of $\zt^*$ contains an
implicit change of universe.)  If $\zt$ is the inclusion of a subgroup
then $\zt^*$ is called \emph{restriction} and written $\res^A_B$.
This functor has a left adjoint written $X\mapsto A_+\Smash_BX$, and a
right adjoint written $X\mapsto F_B(A_+,X)$.  These two adjoints are
actually isomorphic, by the generalized Wirthm\"uller
isomorphism~\cite{lemast:esh}*{Theorem II.6.2}.  

If $\zt$ is the projection $A\xra{}A/B$ then $\zt^*$ is called
\emph{inflation}.  This has a right adjoint functor
$\lm^B\:\CS_A\xra{}\CS_{A/B}$, which we call the Lewis-May fixed point
functor.  The adjunction is discussed in~\cite{lemast:esh}*{Section
II.7}; there $\lm^BX$ is written $X^B$.  One can check
that the following square commutes up to natural isomorphism:
\[ \xymatrix{
 {\CS_A} \rto^{\lm^C} \dto_{\res^A_B} &
 {\CS_{A/C}} \dto^{\res^{A/C}_{B/C}} \\
 {\CS_B} \rto_{\lm^C} & {\CS_{B/C}.}
 } \]
It will be convenient to write
\[ \lmb^B = \res^{A/B}_0\lm^B = \lm^B\res^A_B \:
    \CS_A \xra{} \CS_0. 
\]
The usual equivariant homotopy groups of $X$ are defined by
$\pi^B_*X=\pi_*\lmb^BX$.  The functors $\lm^B$ and $\lmb^B$ do not
preserve smash products, and there is no sense in which $\lm^B$ acts
as the identity on $B$-fixed objects.

Lewis and May also introduce another functor
$\phi^B\:\CS_A\xra{}\CS_{A/B}$, called the \emph{geometric fixed point
  functor}.  To explain the definition, let $V$ be a complex
representation of $A$.  We write $\chi_V$ for the usual inclusion
$S^0\xra{}S^V$, which can be regarded as an element of the
$R(A)$-graded homotopy ring $\pi_*S^0$ in dimension $-V$.  It is
easily seen to be zero if $V^A\neq 0$, but it turns out to be nonzero
otherwise.  It is also clear that $\chi_{V\op W}=\chi_V\chi_W$.

By dualizing the standard cofibration
$S(V)_+\xra{}S^0\xra{\chi_V}S^V$, we see that $D(S(V)_+)$ deserves to
be called $S^0/\chi_V$.  On the other hand, we have
\[ S^0[\chi_V^{-1}] =
    \colim (S^0 \xra{\chi_V} S^V \xra{\chi_V} S^{2V} \xra{} \ldots) =
     S^{\infty V}.
\]
It follows that for any $X\in\CS_A$, the spectrum
$X[\chi_V^{-1}]=X\Smash S^{\infty V}$ is a Bousfield localization of
$X$, or more specifically, the finite localization away from the thick
ideal generated by $S^0/\chi_V$.  There is another characterization as
follows.  Let $\CF$ be the family of those subgroups $A'\leq A$ such
that $V^{A'}\neq 0$, and let $\CC$ be the thick ideal generated by
$\{A/A'_+\st A'\in\CF\}$.  It is not hard to see that
$S(\infty V)^{A'}$ is contractible for $A'\in\CF$ and empty for
$A'\not\in\CF$.  Thus $S(\infty V)$ is a model for the space $E\CF$
(as in Definition~\ref{defn-EF}) and 
$S^{\infty V}=\tSg S(\infty V)=\tSg E\CF=\tE\CF$.  It follows that
$X[\chi_V^{-1}]$ is the finite localization of $X$ away from $\CC$.
It also follows that $\CC$ is the same as the thick ideal generated by
$S^0/\chi_V$. 

Now fix a subgroup $B\leq A$, and take 
\[ V = V_{A,B} := \C[A]\ominus(\C[A]^B) =
    \bigoplus_{\al\in A^*\sm\ann(B)} L_\al.
\]
In this context, we write $\chi_{A,B}$ for $\chi_V$, and we also write
$\chi_A$ for $\chi_{A,A}$.  We also put 
$\CF[B]=\{C\leq A\st B\not\leq C\}$, and note that
$\tE\CF[B]=S^{\infty V}$.  The geometric fixed-point functor
$\phi^B\:\CS_A\xra{}\CS_{A/B}$ is defined by
\[ \phi^BX = \lm^B(X[\chi_{A,B}^{-1}]) = \lm^B(X\Smash\tE\CF[B]). \]
(In~\cite{lemast:esh} the functor $\phi^B$ is actually defined in a
different way, but the above description is proved as Theorem~II.9.8).
Let $\pi\:A\xra{}A/B$ be the projection.  One can check that $\phi^B$
preserves smash products~\cite{lemast:esh}*{Proposition 9.12}, the
composite
\[ \CS_{A/B} \xra{\pi^*} \CS_A \xra{\phi^B} \CS_{A/B} \]
is the identity~\cite{lemast:esh}*{Proposition 9.10}, and the following
diagram commutes:
\[ \xymatrix{
 {\CS_A} \rto^{\phi^C} \dto_{\res^A_B} & 
 {\CS_{A/C}} \dto^{\res^{A/C}_{B/C}} \\
 {\CS_B} \rto_{\phi^C} & {\CS_{B/C}.}
 } \]
Moreover, for any $A$-space $X$ we have $\phi^B\Sgi X=\Sgi X^B$
\cite{lemast:esh}*{Corollary 9.9}.  It will be convenient to write
\[ \phb^B = \res^{A/B}_0\phi^B = \phi^B\res^A_B \:
    \CS_A \xra{} \CS_0. 
\]
This again preserves smash products, and it is known that a spectrum
$X\in\CS_A$ satisfies $X=0$ iff $\phb^BX=0$ in $\CS_0$ for all 
$B\leq A$.  We will also need the following property:
\begin{lemma}\label{lem-F-phi}
 Suppose that $B\leq A$, and write $\chi=\chi_{A,B}$.  Then for
 $X,Y\in\CS_A$ there are natural equivalences
 \[ \lm^B F(X,Y[\chi^{-1}]) = 
    \lm^B F(X[\chi^{-1}],Y[\chi^{-1}]) = 
    F(\phi^BX,\phi^BY).
 \]
\end{lemma}
\begin{proof}
 First note that the map $W\xra{}W[\chi^{-1}]$ is an equivalence iff
 $W$ is concentrated over $B$ as defined in~\cite{lemast:esh}*{page
 109}.  Let $\CC$ be the category of such $W$, so we have
 functors $\phi^B=\lm^B\:\CC\xra{}\CS_{A/B}$ and
 $\psi\:\CS_{A/B}\xra{}\CC$ given by $\psi(Z)=(\pi^*Z)[\chi^{-1}]$.
 We see from~\cite{lemast:esh}*{Corollary II.9.6} that $\phi^B$ and
 $\psi$ are mutually inverse equivalences, and it follows that 
 \[ [X,Y[\chi^{-1}]]^B_* = [X[\chi^{-1}],Y[\chi^{-1}]]^B_* =
     [\phi^BX,\phi^BY]^{A/B}_*.
 \]
 Now consider $W\in\CS_{A/B}$ and replace $X$ by $(\pi^*W)\Smash X$ in
 the above.  We deduce that
 \[ [W,\lm^B F(X,Y[\chi^{-1}])]^{A/B}_* = 
    [W,\lm^B F(X[\chi^{-1}],Y[\chi^{-1}])]^{A/B}_* = 
    [W,F(\phi^BX,\phi^BY)]^{A/B}_*.
 \]
 The claim now follows by the Yoneda lemma.
\end{proof}

\begin{theorem}
 Let $E'$ be an $A/B$-equivariant periodically orientable ring
 spectrum, with associated equivariant formal group
 $(A/B)^*\xra{\phi'}C'$.  Let $\pi\:A\xra{}A/B$ be the projection, and
 put $E=(\pi^*E')[\chi_{A,B}^{-1}]$.  Then $E$ is an $A$-equivariant
 periodically orientable ring spectrum, and for all $X\in\CS_A$ we
 have
 \begin{align*}
  E_*X &= E'_*\phi^BX \\
  E^*X &= (E')^*\phi^BX.
 \end{align*}
 Moreover, the formal group associated to $E$ is the pushout of $C'$
 along the inclusion $(A/B)^*\xra{}A^*$.
\end{theorem}
\begin{proof}
 Because $\pi^*$ preserves smash products, it is clear that $\pi^*E'$
 is a commutative $A$-equivariant ring spectrum, and so the same is
 true of $E$.  We saw earlier that $\phi^B\pi^*=1$, so $\phi^BE=E'$.
 Also, we have $E\Smash X=E\Smash X[\chi_{A,B}^{-1}]$, so
 \[ \lm^B(E\Smash X)=\phi^B(E\Smash X)=\phi^B(E)\Smash\phi^B(X) = 
     E'\Smash\phi^B(X).
 \]
 We can apply $\lm^{A/B}$ to this to see that 
 $\lm^A(E\Smash X)=\lm^{A/B}(E'\Smash\phi^B(X))$, and by applying
 $\pi_*$ we deduce that $E_*X=E'_*\phi^BX$.  

 For the corresponding statement in cohomology, we see using
 Lemma~\ref{lem-F-phi} that
 $\lm^BF(X,E)=F(\phi^BX,\phi^BE)=F(\phi^BX,E')$.  We again apply the
 functor $\pi_*\lm^{A/B}(-)$ to see that $E^*X={E'}^*\phi^BX$, as
 claimed.  

 In particular, if $X$ is an $A$-space we have
 $\phi^B\Sgi X=\Sgi X^B$ and so $E^*X={E'}^*X^B$.  Thus, if we put
 $S=\spec({E'}^0(\text{point}))$, then $S$ is also the same as 
 $\spec(E^0(\text{point}))$.  We next consider the space $PV$, where
 $V$ is a representation of $A$.  We can split $V$ into isotypical
 parts for the action of $B$, say $V=\bigoplus_\bt V[\bt]$, where
 $V[\bt]$ is a sum of representations $L_\al$ with $\al|_B=\bt$.  We
 then have $(PV)^B=\coprod_\bt PV[\bt]$, and so
 $E^*PV=\prod_\bt{E'}^*PV[\bt]$.  Using this, it is easy to see that
 $E$ is periodically orientable.  Next, consider the space $P\CU_A$,
 so $(P\CU_A)^B=\coprod_\bt P(\CU_A[\bt])$.  The space $P(\CU_A[0])$
 is canonically identified with $P\CU_{A/B}$, so
 $\spf({E'}^0P\CU_A[0])=C'$.  For $\bt\neq 0$, we can choose $\tbt\in
 A^*$ extending $\bt$, and then tensoring with $L_{-\tbt}$ gives an
 equivalence $\tht\:P\CU_A[\bt]\simeq P\CU_{A/B}$.  If we change
 $\tbt$ by an element $\gm\in(A/B)^*$, then $\tht$ changes by the
 automorphism $\tau_{-\gm}$ of $P\CU_{A/B}$.  Using this, it is not
 hard to identify the curve
 $C=\spf(E^0P\CU_A)=\coprod\spf({E'}^0P\CU_A[\bt])$ with the pushout
 of $C'$ along the map $(A/B)^*\xra{}A^*$.
\end{proof}

\section{Equivariant Morava $E$-theory}
\label{sec-E-theory}

Let $\hC_0$ be the standard $p$-typical formal group of height $n$
over $S_0=\spec(\Fp)$.  We write $\hK$ for the two-periodic Morava
$K$-theory spectrum whose associated formal group is $\hC_0$, so
$\hK_*=\Fp[u^{\pm 1}]$ with $|u|=2$.  This formal group has a
universal deformation $\hC_1$ over
$S_1:=\spf(\Zp\psb{u_1,\ldots,u_{n-1}})$.  We write $\hE$ for the
corresponding Landweber-exact cohomology theory, and refer to it as
Morava $E$-theory.  Now suppose we have a finite abelian group $A$ and
a subgroup $B$.  We define $C_0=B^*\tm\hC_0$, which is an $A$-efg of
product type over $S_0$, associated to the equivariant Morava
$K$-theory $K^*X:=\hK^*\phb^BX$.  We can also define an
$A/B$-equivariant cohomology theory by $X\mapsto\hE^*X_{h(A/B)}$, as
in Section~\ref{sec-simple}.  The associated equivariant formal group
is $\hC_2=\hC_1\tm_{S_1}S$ over $S$, where $S=\Hom((A/B)^*,\hC_1)$.
We then perform the construction in Section~\ref{sec-pushout}.  This
gives an $A$-equivariant theory $E=E(p,n,B)$, defined by
\[ E^*X = \hE^*((\phi^B X)_{h(A/B)}), \]
whose associated equivariant formal group is the pushout of $\hC_2$
along the inclusion $(A/B)^*\xra{}A^*$.  We write $C$ for this
pushout, and we refer to $E$ as \emph{equivariant Morava $E$-theory}.
In~\cite{st:ebc} we give some evidence that this name is reasonable,
related to the theory of Bousfield classes and nilpotence.  Here we
give a further piece of evidence, based on formal group theory.  

We first note that $S_0$ is a closed subscheme of $S_1$, which is in
turn a closed subscheme of $S=\Hom((A/B)^*,\hC_1)$ (corresponding to
the zero homomorphism).  The restriction of $C$ to $S_1$ is just
$B^*\tm\hC_1$, and the restriction of this to $S_0$ is just $C_0$.
The inclusion $C_0\xra{}C$ corresponds to a ring map
$\OC\xra{}\O_{C_0}$, or equivalently $E^0P\CU\xra{}K^0P\CU$.  It can
be shown that this comes from a natural map $E^*X\xra{}K^*X$ of
cohomology theories.  Indeed, there is certainly a nonequivariant map
$q\:\hE\xra{}\hK$.  Moreover, up to homotopy there is a unique map
$A/B\xra{}E(A/B)$ of $A/B$-spaces, which gives a natural map
\[ \res(Y) = (A/B_+\Smash Y)/(A/B) \xra{}
             (E(A/B)_+\Smash Y)/(A/B) = Y_{h(A/B)} 
\]
for $A/B$-spectra $Y$.  If $Y=\phi^BX$ then $\res(Y)=\phb^BX$ and so
we get a map
\[ E^*X = \hE^*(\phi^BX)_{h(A/B)} \xra{} \hE^*\phb^B \xra{q_*}
    \hK^*\phb^BX = K^*X,
\]
as required.

\begin{definition}
 A \emph{deformation} of the $A$-efg $C_0$ over $S_0$ consists of an
 $A$-efg $C'$ over a base $S'$ together with a commutative square
 \[ \xymatrix{
  {C_0} \rto^{\tf} \dto & {C'} \dto \\
  {S_0} \rto_{f} & {S'}
  } \]
 such that
 \begin{itemize}
  \item[(a)] $f$ is a closed inclusion, and $S'$ is a formal
   neighbourhood of $f(S_0)$
  \item[(b)] $\tf$ induces an isomorphism $C_0\xra{}f^*C$ of $A$-efg's
   over $S_0$.
 \end{itemize}
 If $C'$ and $C''$ are deformations, a morphism between them means a
 commutative square
 \[ \xymatrix{
  {C'} \rto^{\tg} \dto & {C''} \dto \\
  {S'} \rto_{g} & {S''}
  } \]
 such that $\tg$ induces an isomorphism $C'\xra{}g^*C''$ of $A$-efg's
 over $S'$.  A \emph{universal deformation} means a terminal object in
 the category of deformations.
\end{definition}

As mentioned previously, the formal group $\hC_1$ associated to $\hE$
is the universal deformation of the formal group $\hC_0$ associated to
$\hK$.  Equivariantly, we have the following analogue.
\begin{theorem}\label{thm-defm}
 The $A$-equivariant formal group $C$ (associated to equivariant
 Morava $E$-theory) is the universal deformation of $C_0$ (associated
 to equivariant Morava $K$-theory).
\end{theorem}
\begin{proof}
 Suppose we have an $A$-efg $(C',\phi')$ over $S'$ equipped with maps
 $(f,\tf)$ making it a deformation of $C_0$.  We will identify $S_0$
 with $f(S_0)$ and thus regard it as a closed subscheme of $S'$.
 Similarly, we regard $C_0$ as the closed subscheme $C'|_{S_0}$ of
 $C'$.  Note that $S'$ is a formal neighbourhood of $S_0$, and it
 follows that $C'$ is a formal neighbourhood of $C_0$.  We choose a
 coordinate $x'$ on $C'$, and note that it restricts to give a
 coordinate on $C_0$.

 Now let $\hC'$ denote the formal neighbourhood of the zero section in
 $C'$.  We have $(\hC')|_{S_0}=\hC_0$, so we can regard $\hC'$ as a
 deformation of the ordinary formal group $\hC_0$.  As $\hC_1$ is the
 universal deformation of $\hC_0$, this gives us a pullback square
 \[ \xymatrix{
  {\hC'} \rto^{\tg} \dto & {\hC_1} \dto \\
  {S'}   \rto_{g} & {S_1.}
  } \]
 Next, suppose we have $\al\in(A/B)^*\subset A^*$, giving a section
 $\phi'(\al)$ of $C'$ and an element $x'(\phi'(\al))\in\O_{S'}$.  As
 $C'|_{S_0}=C_0=B^*\tm\hC_0$ and $\al|_B=0$ we have
 $\phi'(\al)|_{S_0}=0$, so $x'(\phi'(\al))$ maps to $0$ in
 $\O_{S_0}$.  As $S'$ is a formal neighbourhood of $S_0$, it follows
 that $x'(\phi'(\al))$ is topologically nilpotent in $\O_{S'}$, and
 thus that $\phi'(\al)$ is actually a section of $\hC'$.  Thus,
 $\tg\circ\phi'$ gives a map $(A/B)^*\xra{}\hC_1$, which is classified
 by a map $h\:S'\xra{}\Hom((A/B)^*,\hC_1)=S$.  The maps $\tg$ and $h$
 combine to give a map
 \[ \tht \: \hC' \xra{} h^*\hC = h^*(\hC_1\tm_{S_1}S) = g^*\hC_1. \]
 This can be regarded as an isomorphism of $A/B$-equivariant formal
 groups.  

 Next, the decomposition $C_0=B^*\tm\hC_0=\coprod_{\bt\in B^*}\hC_0$
 gives orthogonal idempotents $e_\bt\in\O_{C_0}$ with
 $\sum_\bt e_\bt=1$.  As $C'$ is a formal neighbourhood of $C_0$,
 these can be lifted to orthogonal idempotents in $\O_{C'}$, giving a
 decomposition $C'=\coprod_\bt C'_\bt$ say.  One can check that
 $C'_\bt=\phi'(\al)+\hC'$ for any $\al\in A^*$ with $\al|_B=\bt$, and
 it follows that $C'$ is just the pushout of the map
 $\phi'\:(A/B)^*\xra{}\hC'$ and the inclusion $(A/B)^*\xra{}A^*$.  It
 follows in turn that $\th$ extends to give an isomorphism
 $C'\xra{}g^*C$, and thus a morphism $C'\xra{}C$ of deformations.  All
 steps in this construction are forced, so one can check that this
 morphism is unique.  This means that $C$ is the universal deformation
 of $C_0$, as claimed.
\end{proof}

\section{A completion theorem}
\label{sec-completion}

Suppose we have an $A$-equivariant formal group $(C,\phi)$, and a
subgroup $B\leq A$, giving a subgroup $(A/B)^*\leq A^*$.  Let $S_0$ be
the closed subscheme of $S$ where $\phi((A/B)^*)=0$.  Equivalently, if
we put $e_\al=x(\phi(-\al))$ and $J=(e_\al\st\al\in(A/B)^*)$, then
$S_0=V(J)=\spec(\OS/J)$.  If we put $C_0=S_0\tm_SC$ then $\phi$
induces a map $\psi\:B^*\tm S_0\xra{}C_0$ making $C_0$ into a
$B$-equivariant formal group over $S_0$.  Next, we put
$S_1=\colim_m\spec(\OS/J^m)=\spf((\OS)^\wedge_J)$, the formal
neighbourhood of $S_1$ in $S$, and $C_1=S_1\tm_SC$.  This is an
$A$-equivariant formal group over $S_1$ for which $\phi((A/B)^*)$ is
infinitesimally close to $0$.

Now suppose that $C$ comes from an $A$-equivariant periodically
orientable theory $E$.  We would like to interpret $C_0$ and $C_1$
topologically.  

\begin{proposition}\label{prop-completion}
 Let $E_0$ be the $B$-spectrum $\res^A_B(E)$, representing the theory
 $E^*(A\tm_BY)$ for $B$-spaces $Y$.  Let $C'_0/S'_0$ be the associated
 $B$-equivariant formal group.  Then there is a map $S'_0\xra{}S_0$
 (which may or may not be an isomorphism) and an isomorphism
 $C'_0=C_0\tm_{S_0}S'_0$.
\end{proposition}
\begin{proof}
 We have $S'_0=\spec(\pi_0E_0)=\spec(E^0A/B)$, so there is a natural
 map $S'_0\xra{}S$.  Moreover, we have $P\CU_B\simeq\res^A_BP\CU_A$,
 which gives an isomorphism $A\tm_BP\CU_B\simeq A/B\tm P\CU_A$ and
 thus
 \[ E_0^0P\CU_B\simeq E^0(A/B\tm P\CU_A)=E^0(A/B)\ot_{E^0}E^0P\CU_A.
 \]
 This shows that the formal group for $E_0$ is just
 $C'_0:=C\tm_SS'_0$.  All that is left is to check that the map
 $S'_0\xra{}S$ factors through $S_0$, so $C'_0$ can also be described
 as $C_0\tm_{S_0}S'_0$.  To see this, note that $\phi$ comes from the
 inclusion $j\:A^*=\pi_0^AP\CU\xra{}P\CU$, so the corresponding map
 $\phi'_0$ over $S'_0$ comes from the map
 $1\tm j\:(A/B)\tm A^*\xra{}(A/B)\tm P\CU$.  Using the isomorphism
 \begin{align*}
  [(A/B)\tm A^*,(A/B)\tm P\CU]^A
    &= [A^*,(A/B)\tm P\CU]^B \\
    &= \Map(A^*,\pi^B_0((A/B)\tm P\CU)) \\
    &= \Map(A^*,(A/B)\tm B^*)
 \end{align*}
 we see that the restriction of $(1\tm j)$ to $(A/B)\tm(A/B)^*$ is
 null, so that $\phi'_0((A/B)^*)=0$ as claimed.
\end{proof}

If $E$ is the complex $K$-theory spectrum $KU_A$, then we saw earlier
that $S=\Hom(A^*,\MG)$ and so
\[ S_0 = \{\phi\in\Hom(A^*,\MG)\st \phi((A/B)^*)=0\} 
       = \Hom(B^*,\MG).
\]
On the other hand, it is well-known that $KU_A^*(A\tm_BY)=KU_B^*Y$ so
$E_0=KU_B$ so $S'_0=\Hom(B^*,\MG)=S_0$.  A similar argument works for
theories of the form $E^*X=\hE^*X_{hA}$ where $\hE$ is $K(n)$-local as
in Section~\ref{sec-simple}, in which case we have $S=\Hom(A^*,\hC)$
and $S_0=S'_0=\Hom(B^*,\hC)$.  At the other extreme, for theories of
the form $E^*X=\hE^*(\res^A_0(X))$, we have $S_0=S$ and
$S'_0=(A/B)\tm S$.

We next consider $C_1$.  Recall that there is an $A$-space $E[\leq B]$
characterised by the property that $E[\leq B]^C$ is contractible for
$C\leq B$ and empty for $C\not\leq B$.  The first approximation would
be to consider the ring spectrum $F(E[\leq B]_+,E)$.  However, as
$S_1$ is a formal scheme rather than an affine scheme, we need a
pro-spectrum rather than a spectrum.  The solution is to define
$F_\bullet(X_+,E)$ to be the pro-system of ring spectra
$F(X_{\al+},E)$, where $X_\al$ runs over finite subcomplexes of $X$,
and to put $E_1=F_\bullet(E[\leq B]_+,E)$.  The desired description of
$E_1^*P\CU$ is a kind of completion theorem in the style of
Atiyah-Segal, so we expect to need finiteness hypotheses.  However,
with these hypotheses, we have an exact result rather than an
approximate one as in the previous proposition.
\begin{theorem}\label{thm-completion}
 Suppose that $E^*(\text{point})$ is a Noetherian ring, and that
 $E^*(A/C)$ is finitely generated over it for all $C\leq A$.  Then the
 $A$-equivariant formal group associated to $E_1$ is $C_1$.
\end{theorem}
\begin{proof}
 This is essentially taken from~\cite{gr:aie}.  Choose generators
 $\al_1,\ldots,\al_r$ for $(A/B)^*$, let $L_i$ be the one-dimensional
 representation corresponding to $\al_i$ and let $\chi_i$ denote the
 inclusion $S^0\xra{}S^{L_i}$.  There is a canonical Thom class $u_i$
 in $E^0S^{L_i}$, and $\chi_i^*(u_i)$ is the Euler class
 $e_i=x(\phi(-\al_i))$.  One checks easily that the space
 $P:=\prod_iS(\infty L_i)$ is a model for $E[\leq B]$, and the spaces
 $T(m):=\prod_iS(mL_i)$ form a cofinal system of finite subcomplexes,
 so $E_1$ is equivalent to the tower of ring spectra
 $F(T(m)_+,E)=D(T(m)_+)\Smash E$.  Next, by taking the
 Spanier-Whitehead dual of the cofibration
 $S(mL_i)_+\xra{}S^0\xra{\chi_i^m}S^{mL_i}$, we see that
 $D(S(mL_i)_+)$ deserves to be called $S/\chi_i^m$, and so $D(T(m)_+)$
 deserves to be called $S/(\chi_1^m,\ldots,\chi_r^m)$.  This suggests
 that $\pi_*(E\Smash D(T(m)_+))$ should be $E_*/J_m$, where
 $J_m=(u_1^m,\ldots,u_r^m)\leq E_*$.  Unfortunately, 
 there are correction terms.  More precisely, the cofibration
 displayed above gives a two-stage filtration of $D(S(mL_i)_+)$ for
 each $i$, and by smashing these together we get a $(r+1)$-stage
 filtration of $D(T(m)_+)$, and thus a spectral sequence converging to
 $\pi_*(E\Smash D(T(m)_+))$.  The first page is easily seen to be the
 Koszul complex for the sequence $e^m_1,\ldots,e^m_r$, so the bottom
 line of the second page is $E_*/J_m$, and the remaining lines are
 higher Koszul homology groups.  The filtrations are compatible as $m$
 varies, so we get a spectral sequence in the abelian category of
 pro-groups converging to $\pi_*E_1$.  In the second page, the bottom
 line is the tower $\{E_*/J_m\}_{m\geq 0}$, and the remaining lines
 are pro-trivial by~\cite{gr:aie}*{Lemma 3.7}.  It follows that
 $\pi_*E_1\simeq\{E_*/J_m\}$ as pro-groups, and so the formal scheme
 corresponding to $\pi_0E_1$ is
 $\colim_m\spec(E^0/J_m)=\colim_m\spec(E^0/J^m)=S_1$.  We now replace
 $E$ by $F(P(n.\C[A])_+,E)$ and then take the limit as $n$ tends to
 infinity to conclude that $\spf(E_1^0P\CU)=C\tm_SS_1=C_1$ as
 claimed. 
\end{proof}

\begin{remark}\label{rem-dominant}
 Using the same circle of ideas one proves that the kernel of the map
 $E^0/J\xra{}E^0(A/B)$ is nilpotent, so the map $S'_0\xra{}S_0$ is
 dominant; compare~\cite{gr:aie}*{Theorem 1.4}.
\end{remark}

\section{A counterexample}
\label{sec-counterexample}

Here we exhibit a $\Z/2$-equivariant formal group $C$ with a number of
unusual properties, which are only possible because the base scheme
$S$ is not Noetherian.  The phenomena described here are the main
obstruction to our understanding of the equivariant Lazard ring.

For any $A$-equivariant formal group $(C,\phi)$, there is a natural
map $\psi\:A^*\tm\hC\xra{}C$ given by $\psi(\al,a)=\phi(\al)+a$.  As
$C$ is a formal neighbourhood of $[\phi(A^*)]$, it is natural to
expect that $\psi$ should be an epimorphism, or equivalently that the
map $\psi^*\:\OC\xra{}\prod_\al\O_{\hC}$ should be injective.  The key
feature of the example to be constructed here is that $\psi^*$ is not
in fact injective.

Start with $k_0=\F_2[e]$, let $M$ be the module 
$\F_2[e^{\pm 1}]/\F_2[e]$, and let $k$ be the square-zero extension
$k_0\oplus M$.  More explicitly, $k$ is generated over $k_0$ by
elements $u_1,u_2,\ldots$ subject to $eu_{i+1}=u_i$ (with $u_0$
interpreted as $0$) and $u_iu_j=0$.  Put $S=\spec(k)$.

Next, let $R$ be the completion of $k[x]$ at the element $y=x^2+ex$,
so $R=k\psb{y}\{1,x\}$, and put 
\[ C = \spf(R)
     = \{ x \in \aff^1_S \st x^2+ex \text{ is nilpotent } \}.
\]
This is a subgroup of $\aff^1_S$ under addition.  In the corresponding
Hopf algebra structure on $R$, the elements $x$ and $y$ are both
primitive.  There is a homomorphism $\phi\:\Z/2\xra{}C$ sending $0$ to
$0$ and $1$ to $e$.  The corresponding divisor is just $R/y$, and as
$R$ is complete at $y$, we deduce that $(C,\phi)$ is an equivariant
formal group.  

Next, we can define maps $\lm_0,\lm_a\:R\xra{}k\psb{t}$ by 
\begin{align*}
 \lm_0(x) &= t \\
 \lm_a(x) &= t+e \\
 \lm_0(y) &= \lm_a(y) = t^2 + te.
\end{align*}
The map $\psi^*\:\OC\xra{}\prod_\al\O_{\hC}$ is just the map
$(\lm_0,\lm_a)\:R\xra{}k\psb{t}\tm k\psb{t}$.  Now consider the
element 
\[ f =  \sum_{k\geq 0} u_{1-2^{k+1}} y^{2^k} \in R. \]
We then have
\begin{align*}
 \lm_0(f)  &= \sum_{k\geq 0} u_{1-2^{k+1}} (t^2 + et)^{2^k} \\
           &= \sum_{k\geq 0} u_{1-2^{k+1}} t^{2^{k+1}} +
              \sum_{k\geq 0} u_{1-2^{k+1}} e^{2^k} t^{2^k} \\
           &= \sum_{k\geq 0} u_{1-2^{k+1}} t^{2^{k+1}} +
              \sum_{k\geq 0} u_{1-2^k} t^{2^k} \\
           &= u_{1-2^0} t^{2^0} = u_0t = 0.
\end{align*}
We also have $\lm_a(f)=0$ by the same argument, so $\psi^*(f)=0$.

\section{Divisors}
\label{sec-divisors}

We now return to the purely algebraic theory of formal multicurves and
their divisors.

Recall that a divisor on $C$ is a regular hypersurface $D\sse C$ such
that $\OD$ is a finitely generated projective module over $\OS$, which
is discrete in the quotient topology.  We also make the following
temporary definition; one of our main tasks in this section is to show
(in Proposition~\ref{prop-weak-genuine}) that it is equivalent to the
preceeding one.  (For divisors of degree one, this follows from
Corollary~\ref{cor-sec-regular}.)
\begin{definition}\label{defn-weak-divisor}
 A \emph{weak divisor} on $C$ is a closed subscheme $D\subset C$ that
 is finite and very flat over $S$ (so $\OD$ is a discrete finitely
 generated projective module over $\OS$).  Thus, a weak divisor
 $D=\spf(R/J)$ is a divisor iff the ideal $J$ is open and generated by
 a regular element.  If $y$ is a good parameter on $C$, we note that
 $J$ is open iff $y^N\in J$ for $N\gg 0$.
\end{definition}
If $D_0=\spf(R/J_0)$ is a divisor and $D_1=\spf(R/J_1)$ is a weak
divisor then one checks that the scheme $D_0+D_1:=\spf(R/(J_0J_1))$ is
again a weak divisor.

\begin{definition}\label{defn-norm}
 Now suppose we have a map $q\:T\xra{}S$ of schemes which is finite
 and very flat, so that $\OT$ is a discrete finitely generated
 projective module over $\OS$.  If $g\in\O_T$ then multiplication by
 $g$ gives an $\OS$-linear endomorphism $\mu_g$ of $\OT$, whose
 determinant we denote by $N_q(g)$ or $N_{T/S}(g)$.
\end{definition}

\begin{definition}\label{defn-fD}
 Fix a difference function $d$ on $C$.  For any weak divisor $D$ on
 $C$ over, we can regard $d$ by restriction as a function on
 $D\tm_SC$.  We also have a projection $q\:D\tm_SC\xra{}C$, and we put 
 \[ f_D = N_q(d) = N_{(D\tm_SC)/C}(d) \in \OC. \]
\end{definition}

We will eventually show that $D=\spf(\OC/f_D)$.

\begin{remark}\label{rem-ordinary-difference}
 Consider the case where $C$ is an ordinary formal group, with
 coordinate $x$ and associated formal group law $F$.  We then have
 $\OC=\OS\psb{x}$ and $\O_{C\tm C}=\OS\psb{x_0,x_1}$, and we can take
 $d=x_1-_Fx_0$.  If $D$ has the form $\sum_i[u_i]$ for some family of
 sections $u_i$, then we have elements $a_i=x(u_i)\in\OS$ and we will
 see that $f_D=\prod_i(x-_Fa_i)$.  This is a unit multiple of the
 Chern polynomial $g_D=\prod_i(x-a_i)$, and it is familiar that
 $D=\spf(\OC/g_D)$, so $D=\spf(\OC/f_D)$ also.  In the multicurve
 case, one can still define $g_D$ (as the norm of the function
 $(a,b)\mapsto x(b)-x(a)$) and we find that it is divisible by
 $f_D$, but $g_D/f_D$ need not be invertible so $\OC/g_D\neq\OD$.
\end{remark}

\begin{lemma}\label{lem-det-inj}
 Let $R$ be a ring, $P$ a finitely generated projective $R$-module,
 and $\al$ an automorphism of $P$.  Then $\al$ is injective iff
 $\det(\al)$ is a regular element.
\end{lemma}
\begin{proof}
 After localising we may assume that $P=R^d$ for some $d$, and $\al$
 is represented by a $d\tm d$ matrix $A$.  If $\det(A)$ is regular,
 the equation $\adj(A)A=\det(A)I_d$ implies immediately that $\al$ is
 injective.  Conversely, suppose that $\al$ is injective.  As $P$ is
 flat, it follows that $\al^{\ot d}\:P^{\ot d}\xra{}P^{\ot d}$ is also
 injective.  It is easy to check with bases that $\lm^dP$ is naturally
 isomorphic to the image of the antisymmetrisation map
 $P^{\ot d}\xra{}P^{\ot d}$.  In particular, it embeds naturally in
 $P^{\ot d}$, and it therefore follows that $\lm^d\al$ is injective.
 On the other hand, $\lm^dP$ is an invertible $R$-module, so
 $\End(\lm^dP)=R$, and $\lm^d(\al)=\det(\al)$ under this isomorphism.
 It follows that $\det(\al)$ is regular as claimed.
\end{proof}
\begin{corollary}\label{cor-fD-regular}
 For any weak divisor $D$ on $C$, the element $f_D\in\OC$ is regular.
\end{corollary}
\begin{proof}
 Take $R=\OC$ and $P=\O_{D\tm_SC}$ and $\al=\mu_d$.  We know from
 Lemma~\ref{lem-res-reg} that $\al$ is injective, and the claim
 follows.
\end{proof}

\begin{lemma}\label{lem-norm-zero}
 Let $q\:T\xra{}S$ be finite and very flat, and let $g$ be a function
 on $T$.  If there is a section $u\:S\xra{}T$ such that $g\circ u=0$
 then $N_q(g)=0$.
\end{lemma}
\begin{proof}
 Put $J=\ker(u^*\:\OT\xra{}\OS)$, so $g\in J$.  We have a short exact
 sequence of $\OS$ modules $J\xra{}\OT\xra{u^*}\OS$, which is split by
 the map $q^*\:\OS\xra{}\OT$.  The sequence is preserved by $\mu_g$,
 and $\mu_g(\OT)=\OT.g\leq J$ so the induced map on the cokernel is
 zero.  Zariski-locally on $S$ we can choose bases adapted to the
 short exact sequence and it follows easily that $\det(\mu_g)=0$ as
 claimed.
\end{proof}
\begin{corollary}\label{cor-fD-zero}
 The function $f_D\in\OC$ vanishes on $D$.
\end{corollary}
\begin{proof}
 We have $f_D|_D=N_{q'}(d)$, where $q'\:D\tm_SD\xra{}D$ is the
 projection on the second factor.  The diagonal map
 $\dl\:D\xra{}D\tm_SD$ is a section of $q'$ with $d\circ\dl=0$, so
 $N_{q'}(d)=0$.
\end{proof}

\begin{lemma}\label{lem-norm-product}
 If $D=D_0+D_1$ (where $D_0,D_1$ are divisors) and $g\in\OD$
 then 
 \[ N_{D/S}(g) = N_{D_0/S}(g)N_{D_1/S}(g). \]
\end{lemma}
\begin{proof}
 Put $R=\OC$, and let the ideals corresponding to $D_i$ be $J_i=(f_i)$
 for $i=0,1$.  We then have a short exact sequence of $\OD$-modules as
 follows:
 \[ \O_{D_0} = R/f_0 \xra{\tm f_1} \OD=R/(f_0f_1)
     \xra{} \O_{D_1} = R/f_1.
 \]
 This is splittable, because $\O_{D_1}$ is projective over $\OS$.
 The map $\mu_g$ preserves the sequence, and it follows easily that
 $\det(\mu_g)=\det(\mu_g|\O_{D_0})\det(\mu_g|\O_{D_1})$, as required.
\end{proof}
\begin{corollary}\label{cor-fD-product}
 If $D=D_0+D_1$ as above then $f_D=f_{D_0}f_{D_1}$.
\end{corollary}
\begin{proof}
 Just change base to $C$ and take $g=d$.
\end{proof}

\begin{lemma}\label{lem-div-subtract}
 Suppose that $D$ is a weak divisor of degree $r$, that $D'$ is a
 divisor of degree $r'$, and that $D'\sse D$.  Then $D=D'+D''$ for
 some weak divisor $D''$ of degree $r-r'$.
\end{lemma}
\begin{proof}
 Put $J=I_D$ and $J'=I_{D'}$.  As $D'$ is a genuine divisor, we have
 $J'=Rf'$ for some regular element $f'\in R$.  As $D'\sse D$, we have
 $J\leq J'$.  Put $J''=\{g\in R\st f'g\in J\}\geq J$.  We then have a
 short exact sequence
 \[ R/J'' \xra{\tm f'} R/J \xra{} R/J'. \]
 As $R/J$ and $R/J'$ are projective modules of ranks $r$ and $r'$ over
 $k$, it follows that $R/J''$ is a projective module of rank $r-r'$.
 Thus, the scheme $D'':=\spf(R/J'')$ is a weak divisor.  From the
 definition of $J''$ we have $J'J''\leq J$.  Conversely, if $h\in J$
 then certainly $h\in J'=Rf'$ so $h=gf'$ for some $g\in R$.  From the
 definitions we have $g\in J''$, so $h\in J''J'$.  This shows that
 $J=J'J''$ and so $D=D'+D''$. 
\end{proof}

\begin{definition}\label{defn-full-set}
 Let $D$ be a weak divisor of constant degree $r$.  A \emph{full set
 of points} for $D$ is a list $u_1,\ldots,u_r$ of sections of $S$ such
 that $D=\sum_i[u_i]$.  If there exists a full set of points, it is
 clear that $D$ is actually a genuine divisor.  (This concept is due
 to Drinfeld, and is explained and used extensively
 in~\cite{kama:ame}.) 
\end{definition}

\begin{proposition}\label{prop-full-set-norm}
 If $u_1,\ldots,u_r$ is a full set of points for $D$, then
 $N_{D/S}(g)=\prod_ig(u_i)$ for any function $g$ on $D$.  Moreover, we
 have $f_D(a)=\prod_id(a,u_i)$, and so $\OD=\OC/f_D$.
\end{proposition}
\begin{proof}
 As the projection $[u_i]\xra{}S$ is an isomorphism, we see that
 $N_{[u_i]/S}(g)=g(u_i)$.  The first claim follows easily using from
 Lemma~\ref{lem-norm-product} by induction on $r$.  It
 follows similarly from Corollary~\ref{cor-fD-product} that
 $f_D(a)=\prod_id(a,u_i)$.  As $d$ is a difference function we have
 $\O_{[u_i]}=\OC/d(a,u_i)$ and so $\OD=\OC/f_D$ as claimed.
\end{proof}

\begin{lemma}\label{lem-full-set}
 Let $D$ be a weak divisor of constant degree $r$.  Then there is a
 finite, very flat scheme $T$ over $S$ such that the weak divisor
 $T\tm_SD$ on $T\tm_SC$ has a full set of points (and so is genuine).
\end{lemma}
\begin{proof}
 By an evident induction, it suffices to show that after very flat
 base change we can split $D$ as $[u]+D''$ for some section $u$ and
 some weak divisor $D''$.  It is enough to find a section
 $u\:S\xra{}D$, for then $[u]\sse D$ and we can apply the previous
 lemma.  For this we can simply pull back along the projection map
 $D\xra{}S$ (which is very flat by assumption) and then the diagonal
 map $D\xra{}D\tm_SD$ gives the required ``tautological'' section.
\end{proof}

\begin{proposition}\label{prop-weak-genuine}
 Every weak divisor is a genuine divisor.
\end{proposition}
\begin{proof}
 Let $D=\spf(R/J)$ be a weak divisor.  We may assume without loss that
 it has constant degree $r$.  We know from
 Corollary~\ref{cor-fD-regular} and Corollary~\ref{cor-fD-zero} that
 $f_D$ is regular in $R$ and lies in $J$; we need only show that it
 generates $J$.  It is enough to do this after faithfully flat base
 change, so by Lemma~\ref{lem-full-set} we may assume that we have a
 full set of points.  Proposition~\ref{prop-full-set-norm} completes
 the proof.
\end{proof}

\section{Embeddings}
\label{sec-embeddings}

Let $C$ be a nonempty formal multicurve
over a scheme $S$.  In this section we study embeddings of $S$ in the
affine line $\aff^1_S=\aff^1\tm S$.  If $q$ is the given map
$C\xra{}S$, then any map $C\xra{}\aff^1_S$ of schemes over $S$ has the
form $(x,q)$ for some $x\:C\xra{}\aff^1$, or equivalently $x\in\OC$.

Now choose a difference function $d$ on $C$.  Given $x\in\OC$, we can
define $x'\:C\tm_SC\xra{}\aff^1$ by $x'(a,b)=x(b)-x(a)$.
Equivalently, $x'$ is the element $1\ot x-x\ot 1$ in
$\O_{C\tm_SC}=\OC\hot_{\OS}\OC$.  It is clear that $x'$ vanishes on
the diagonal, and thus is divisible by $d$, say $x'=\tht(x)d$ for some
$\tht(x)\in\O_{C\tm_SC}$.  This element $\tht(x)$ is unique, because
$d$ is not a zero-divisor.

\begin{proposition}
 Let $C\xra{q}S$ be a nonempty formal multicurve.  A map
 $(x,q)\:C\xra{}\aff^1_S$ is injective if and only if $\tht(x)$ is
 invertible.  If so, then $(x,q)$ induces an isomorphism
 $C\xra{}\colim_kV(f^k)\subset\aff^1_S$ for some monic polynomial
 $f\in\OS[t]$, showing that $C$ is embeddable.
\end{proposition}
\begin{proof}
 Put $X=\{(a,b)\in C\tm_SC\st x(a)=x(b)\}=V(x')=V(\tht(x)d)$.  We see
 that $x$ is injective if and only if $V(x')=\Dl=V(d)$, if and only if
 $d=ux'=u\tht(x)d$ for some $u\in\O_{C\tm_SC}$.  As $d$ is not a zero
 divisor, this holds if and only if $\tht(x)$ is invertible.

 If so, we may assume without loss that $d=x'$.  Choose a good
 parameter $y$, so $\OC/y$ has constant rank $r$ over $\OS$ for some
 $r$.  Put $D=\spec(R/y)$, let $p\:C\tm_SD\xra{}C$ be the projection,
 and put $z=N_p(x')$.  The proof of
 Proposition~\ref{prop-weak-genuine} shows that $z$ is a unit multiple
 of $y$.

 We next claim that $\{1,x,\ldots,x^{r-1}\}$ is a basis for $R/y=R/z$
 over $k$, and that $z=f(x)$ for a unique monic polynomial $f$ of
 degree $r$.  It is enough to check this after faithfully flat base
 change, so we may assume that $D=[u_0]+\ldots+[u_{r-1}]$ for some list of
 sections $u_i$ of $D$.  If we put $a_i=x(u_i)$ we see that
 $z=\prod_i(x-a_i)$.  If we put $e_i=\prod_{j<i}(x-a_j)$ we also find
 that $\{e_0,\ldots,e_{n-1}\}$ is a basis for $R/z$.  As
 $e_i=x^i+\text{ lower terms }$, we also find that
 $\{1,\ldots,x^{n-1}\}$ is a basis as claimed.

 The rest of the proposition follows easily from this.  
\end{proof}

Now suppose we have an arbitrary element $x\in\OC$.  Given a map
$u\:S'\xra{}S$ we get a multicurve $C':=S'\tm_SC$ over $S'$ and a
function $x'=(C'\xra{}C\xra{x}\aff^1)\in\O_{C'}$.
\begin{lemma}
 There is a basic open subscheme $U\subset S$ such that
 $(x',q')\:C'\xra{}\aff^1_{S'}$ is an embedding if and only if
 $u\:S'\xra{}S$ factors through $U$.
\end{lemma}
\begin{proof}
 Choose a good parameter $y$ on $C$ and put $D=\spec(\OC/y)$.  Put
 $w=N_{D\tm_SD/S}(\tht(x))\in\OS$.  We see that $w$ is invertible in
 $\OS$ if and only if $\tht(x)$ is invertible in $\O_{D\tm_SD}$.  As
 $\OC$ is complete at $y$, we see that $\tht(x)$ is invertible in
 $\O_{C\tm_SC}$ if and only if it is invertible in $\O_{D\tm_SD}$.
 Given this, it is clear that the scheme $U=\spec(\OS[1/w])$ has the
 stated property.
\end{proof}

\begin{corollary}\label{cor-locally-embeddable}
 Let $C$ be a formal multicurve over $S$.  Then there is a faithfully
 flat map $S'\xra{}S$ such that the pullback $C':=S'\tm_SC$ is
 embeddable.
\end{corollary}
\begin{proof}
 Put $R=\OC$ and $k=\OS$, and let $y$ be a very good parameter on $C$.
 Let $P$ be the continuous dual of $R$, which is a projective module
 of countable rank over $k$.  We have
 $P\hot_kR\simeq\Hom^{\text{cts}}_k(R,R)$, so there is an element
 $x\in P\hot_kR$ corresponding to the identity map $1_R$.  The scheme
 $M:=\Map_S(C,\aff^1)$ is the spectrum of the symmetric algebra
 $k[P]$, with the tautological map $M\tm_SC\xra{}\aff^1$ corresponding
 to the element $x\in P\hot_kR\subset k[P]\hot_kR=\O_{M\tm_SC}$.  As
 in the lemma, there is a largest open subscheme $S'\sse M$ where this
 tautological map gives an embedding $S'\tm_SC\xra{}\aff^1_{S'}$.
 Note that $M=\spec(k[P])$ is flat over $S$ and $S'$ is open in $M$,
 it is again flat over $S$.  It is clear by construction that
 $S'\tm_SC$ has a canonical embedding in $\aff^1_{S'}$.  All that is
 left is to check that the map $u\:S'\xra{}S$ is \emph{faithfully}
 flat.  It will suffice to show that $u$ is surjective on geometric
 points, and this follows easily from Lemma~\ref{lem-alg-closed}.
\end{proof}

\section{Symmetric powers of multicurves}
\label{sec-symmetric}

In this section, we study the formal schemes $C^r/\Sg_r$, or in other
words the symmetric powers of $C$.  As usual, we write $R=\OC$ and
$k=\OS$.  We choose a good parameter $y$ on $C$, and a basis
$\{e_0,\ldots,e_{n-1}\}$ for $R/y$.  We then put $e_{ni+j}=y^ie_j$,
which gives a topological basis $\{e_i\st i\geq 0\}$ for $R$ over $k$
and thus an isomorphism $R\simeq\prod_ik$ of topological $k$-modules.
We write 
\begin{align*}
 R_r   &= R \hot_k \ldots \hot_k R \\
 S_r   &= R_r^{\Sg_r} \\
 \Rb   &= k\psb{y} \\
 \Rb_r &= \Rb \hot_k \ldots \hot_k \Rb = k\psb{y_1,\ldots,y_r} \\
 u_i   &= \text{ $i$'th elementary symmetric function of 
                 $y_1,\ldots,y_r$ } \\
 \Sb_r &= \Rb_r^{\Sg_r} = k\psb{u_1,\ldots,u_r} \\
 C^r   &= C\tm_S\ldots\tm_SC = \spf(R_r) \\
 C^r/\Sg_r &= \spf(S_r) \\
 \Cb   &= \spf(\Rb) \\
 \Cb^r &= \Cb\tm_S\ldots\tm_S\Cb = \spf(\Rb_r) \\
 \Cb^r/\Sg_r &= \spf(\Sb_r). 
\end{align*}
Here we have topologized $\Rb_r$, $S_r$ and $\Sb_r$ as closed subrings
of $R_r$.  We clearly have a commutative square of topological rings
as shown on the left below, and thus a commutative square of formal
schemes as shown on the right.
\[ \xymatrix{
 {R_r} & 
 {S_r} \ar@{ >->}[l] &
 {C^r} \ar@{->>}[r] \ar@{->>}[d] & 
 {C^r/\Sg_r} \ar@{->>}[d] \\
 {\Rb_r} \ar@{ >->}[u] & 
 {\Sb_r} \ar@{ >->}[l] \ar@{ >->}[u] &
 {\Cb^r} \ar@{->>}[r] &
 {\Cb^r/\Sg_r}
 } \]
We next exhibit topological bases for the above rings.  Put
\begin{align*}
 A   &= \N^r \\
 \Ab &= (n\N)^r = \{\al\in A\st \al_i=0\pmod{n}\text{ for all } i\} \\
 B   &= \{\bt\in\N^{\infty}\st \sum_{i=0}^\infty\bt_i=r\} \\
 \Bb &= \{\bt\in B\st \bt_i=0 \text{ whenever } i\neq 0\pmod{n}\}.
\end{align*}
Next, for $\al\in A$ we put 
\[ e_\al = e_{\al_1} \ot\ldots\ot e_{\al_r}\in R_r. \]
Note that $e_{n\al}=\prod_{i=1}^ry_i^{\al_i}\in\Rb_r$, and
$e_{n\al+\al'}=e_{n\al}e_{\al'}$.  

Now define $\tau\:A\xra{}B$ by $\tau(\al)_j=|\{i\st\al_i=j\}|$.  This
gives bijections $A/\Sg_r=B$ and $\Ab/\Sg_r=\Bb$.  For $\bt\in B$, we
put 
\[ e'_\bt = \sum_{\tau(\al)=\bt} e_\al \in S_r. \]
It is clear that
\begin{itemize}
 \item $\{e_\al\st\al\in A\}$ is a topological basis for $R_r$ over
  $k$, giving an isomorphism $R_r=\prod_Ak$.
 \item $\{e_\al\st\al\in\Ab\}$ is a topological basis for $\Rb_r$.
 \item $\{e'_\bt\st\bt\in B\}$ is a topological basis for $S_r$.
 \item $\{e'_\bt\st\bt\in \Bb\}$ is a topological basis for $\Sb_r$.
\end{itemize}
Of course, the monomials in the symmetric functions $u_i$ give another
topological basis for $\Sb_r$ over $k$.

\begin{proposition}
 If $S'=\spec(k')$ is any scheme over $S$, and $C'=S'\tm_SC$
 (considered as a multicurve over $S'$) then
 $(C')^r/\Sg_r=S'\tm_S(C^r/\Sg_r)$.  The schemes $\Cb$, $\Cb^r$ and
 $\Cb^r/\Sg_r$ are also compatible with base change in the same sense.  
\end{proposition}
\begin{proof}
 Put $R'=\O_{C'}=k'\hot_kR=\prod_{i\in\N}k'$, and
 $R'_r=\O_{(C')^r}=R'\hot_{k'}\ldots\hot_{k'}R'=\prod_{\al\in A}k'$,
 and $S'_r=\O_{(C')^r/\Sg_r}=\prod_{\bt\in B}k'$.  This is clearly the
 same as $k'\hot_kS_r$, so $(C')^r/\Sg_r=(C^r/\Sg_r)\tm_SS'$.  The
 same argument works for the other claims.
\end{proof}

We next need to formulate and prove various compatibility statements
for the topologies on the rings considered above.
\begin{definition}\label{defn-top-free}
 Let $A$ be a linearly topologised ring, and let $M$ be a topological
 module over $A$.  We say that $M$ is \emph{topologically free of rank
   $r$} if it is isomorphic to $A^r$ (with the product topology) as a
 topological module.
\end{definition}
\begin{definition}\label{defn-neat}
 Let $A$ be a linearly topologised ring, and let $B$ be a closed
 subring (with the subspace topology).  We write $I\leq_OA$ to
 indicate that $I$ is an open ideal in $A$.  We say that $B$ is
 \emph{neat} if for every open ideal $J\leq_OB$, the ideal $JA$ is
 open in $A$.
\end{definition}
\begin{remark}\label{rem-neat}
 As $B$ has the subspace topology, we see that
 $\{I\cap B\st I\leq_OA\}$ is a basis of open ideals in $B$.  It
 follows that $B$ is neat iff $(I\cap B)A$ is open in $A$ whenever
 $I\leq_OA$.  If so, then (using the inclusion $(I\cap B)A\leq I$) we
 see that $\{(I\cap B)A\st I\leq_OA\}$ is a basis of open ideals in
 $A$. 
\end{remark}
\begin{remark}\label{rem-free-neat}
 Suppose that we start with a linear topology on $B$.  We can then
 give $A$ a linear topology by declaring $\{AJ\st J\leq_OB\}$ to be a
 basis of open ideals in $A$.  By regarding $B$ as a subspace of $A$,
 we obtain a new linear topology on $B$, which may or may not be the
 same as the old one.  Now suppose that $A$ is faithfully flat over
 $B$.  It follows that $A/JA=A\ot_BB/J$ is faithfully flat over $B/J$,
 and in particular that the map $B/J\xra{}A/JA$ is injective, so
 $J=(JA)\cap B$.  Using this we see that the two topologies on $B$ are
 the same, and that $B$ is neat in $A$.

 In particular, if $A$ is topologically free of finite rank over $B$,
 then $B$ is neat in $A$.  Conversely, if $A$ is free of finite rank
 over $B$ and $B$ is neat, it is easy to see that $A$ is topologically
 free. 
\end{remark}

\begin{proposition}\label{prop-topologies}
 \begin{itemize}
  \item[(a)] $R_r$ is topologically free of rank $n^r$ over $\Rb_r$
  \item[(b)] $\Rb_r$ is topologically free of rank $r!$ over $\Sb_r$
  \item[(c)] $R_r$ is topologically free of rank $n^rr!$ over $\Sb_r$
  \item[(d)] $S_r$ is topologically free of rank $n^r$ over $\Sb_r$
  \item[(e)] $R_r$ is a finitely generated module over $S_r$, and
   $S_r$ is neat in $R_r$.
 \end{itemize}
 Moreover, in each of the four rings there is a finitely generated
 ideal $J$ such that $\{J^m\st m\geq 0\}$ is a basis of open ideals.
\end{proposition}
The proof will follow after a number of lemmas.  In
Corollary~\ref{cor-Rr-projective}, we will extend part~(e) by proving
that $R_r$ is a projective module of rank $r!$ over $R_r$.
\begin{lemma}\label{lem-ideal-powers}
 Suppose we have a ring $A$ and elements $a_1,\ldots,a_n\in A$, and we
 put $I_m=(a_1^m,\ldots,a_n^m)$.  Then 
 $I_1^{n(m-1)+1}\leq I_m\leq I_1^m$, so the topology defined by the
 ideals $I_m$ is the same as that defined by the ideals $I_1^m$.
\end{lemma}
\begin{proof}
 The ideal $I_1^m$ is generated by the monomials $\prod_ia_i^{v_i}$
 for which $\sum_iv_i=m$.  It is clear from this that $I_m\leq I_1^m$,
 and thus that $I_1^m$ is open in the topology defined by the ideals
 $I_k$.  Now suppose we have a monomial $\prod_ia_i^{v_i}$ that is not
 contained in $I_m$.  This means that $v_i\leq m-1$ for all $i$, and
 thus $\sum_{i=1}^nv_i<n(m-1)+1$.  By the contrapositive, we see that
 $I_1^{n(m-1)+1}\leq I_m$, so $I_m$ is open in the topology defined by
 the ideals $I_1^k$.
\end{proof}
\begin{lemma}\label{lem-invariants-neat}
 Let $A$ be a linearly topologised ring with a continuous action of a
 finite group $G$.  Suppose that there exists a finitely generated
 $G$-invariant ideal $I=(a_1,\ldots,a_n)$ such that
 $\{I^m\st m\geq 0\}$ is a basis of open ideals.  Then $A^G$ is neat
 in $A$.  Moreover, if $A$ is faithfully flat over $A^G$ then
 $\{(I^G)^m\st m\geq 0\}$ is a basis of open ideals in $A^G$.
\end{lemma}
\begin{proof}
 Put $r=|G|$.  For any $a\in A$, put
 $\phi_a(t)=\prod_{g\in G}(t-g.a)$, so $\phi_a(a)=0$.  If $J$ is any
 $G$-invariant ideal containing $a$ we see that
 $\phi_a(t)\in t^r+J^G[t]$, so the equation $\phi_a(a)=0$ gives
 $a^r\in A.J^G$.  Thus, all elements of $J$ are nilpotent modulo
 $A.J^G$.  If $J$ is finitely generated we deduce that there exists
 $s>0$ with $J^s\leq A.J^G$.  
 
 Now apply this with $J=I^m$; we see that $A.(I^m)^G$ contains
 $I^{ms}$ for some $s$, and thus is open.  This shows that $A^G$ is
 neat in $A$.

 Now suppose that $A$ is faithfully flat over $A^G$.  We claim that
 $(I^G)^m$ is open in $A^G$.  Indeed, the above shows that for large
 $j$ we have $I^j\leq A.I^G$.  It follows that
 $I^{jm}\leq(I^GA)^m=A.(I^G)^m$.  It is also clear that
 $A.(I^{jm})^G\leq I^{jm}$, so $A.(I^{jm})^G\leq A.(I^G)^m$.  By
 faithful flatness, for any ideals $J,J'\leq A^G$ we have $A.J\leq
 A.J'$ iff $J\leq J'$.  We deduce that $(I^{jm})^G\leq(I^G)^m$.  The
 ideal $(I^{jm})^G=I^{jm}\cap A^G$ is open in the subspace topology,
 so the same is true of $(I^G)^m$.  We also have $(I^G)^m\leq(I^m)^G$
 and the ideals $(I^m)^G$ form a basis of neighbourhoods of $0$; it
 follows that the same is true of the ideals $(I^G)^m$.
\end{proof}
\begin{corollary}\label{cor-invariants-neat}
 Let $A$, $I$ and $G$ be as in the lemma, and let $H$ be a subgroup of
 $G$.  Suppose that the inclusion $A^H\xra{}A$ is faithfully flat,
 and that $I^H$ is finitely generated.  Then $A^G$ is neat in $A^H$.
\end{corollary}
\begin{proof}
 The lemma (with $G$ replaced by $H$) tells us that
 $\{(I^H)^m\st m\geq 0\}$ is a basis of open ideals in $A^H$.  As
 $I^H$ is finitely generated, the same is true of $(I^H)^m$, say
 $(I^H)^m=(b_1,\ldots,b_n)$.  Consider the polynomial
 $\phi_{b_i}(t)=\prod_g(t-g.b_i)$ as in the proof of the lemma.  As
 $b_i\in(I^H)^m\sse I^m$, we see that
 $\phi_{b_i}(t)\in t^r+(I^m)^G[t]$.  Using the relation
 $\phi_{b_i}(b_i)=0$ we see that $b_i^r\in(I^m)^GA^H$, so
 \[ (I^H)^{m(n(r-1)+1)}\leq (b_1^r,\ldots,b_n^r)\leq (I^m)^GA^H, \]
 so $(I^m)^GA^H$ is open in $A^H$.  As the ideals $(I^m)^G$ are a
 basis of open ideals in $A^G$, we deduce that $A^G$ is neat as
 claimed. 
\end{proof}
\begin{lemma}\label{lem-topologies}
 Suppose that $A=k[y_1,\ldots,y_r]$, with the evident action of
 $G=\Sg_r$, and with topology determined by the powers of the ideal
 $I=(y_1,\ldots,y_r)$.  Let $H$ be a subgroup of $G$ of the form
 $\Sg_{r_1}\tm\ldots\tm\Sg_{r_k}$, with $r=r_1+\ldots+r_k$.  Then
 \begin{itemize}
  \item[(a)] $A$ is topologically free of rank $|G|=r!$ over $A^G$
  \item[(b)] $A$ is topologically free of rank $|H|=\prod_ir_i!$ over $A^H$
  \item[(c)] $A^H$ is topologically free of rank $|G/H|$ over $A^G$
  \item[(d)] The topology on $A^H$ (resp. $A^G$) is determined by
   powers of the ideal $I^H$ (resp. $I^G$), which is finitely
   generated.
 \end{itemize}
\end{lemma}
\begin{proof}
 It is well-known that $A^G=k[u_1,\ldots,u_r]$, where $u_i$ is the
 $i$'th elementary symmetric function in the variables $y_i$.
 Similarly, we have $A^H=k[v_1,\ldots,v_r]$, where
 $v_1,\ldots,v_{r_1}$ are the elementary symmetric functions of
 $y_1,\ldots,y_{r_1}$, and $v_{r_1+1},\ldots,v_{r_1+r_2}$ are the
 elementary symmetric functions of $y_{r_1+1},\ldots,y_{r_1+r_2}$ and
 so on.  By considering the maps 
 \[ A^G \xra{} A^H \xra{} A \xra{} A/I = k, \]
 we see that $I^G=(u_1,\ldots,u_r)$ and $I^H=(v_1,\ldots,v_r)$, so in
 particular these ideals are finitely generated.  

 We next claim that $A$ is algebraically free of rank $|H|$ over
 $A^H$.  Everything is compatible with base change, so it will be
 enough to prove this when $k=\Z$.  In this case, all the rings
 involved are Noetherian domains with unique factorisation and the
 claim is a standard piece of invariant theory.  Similarly, we see
 that $A$ and $A^H$ are algebraically free of the indicated ranks over
 $A^G$, and so the inclusions $A^G\xra{}A^H\xra{}A$ are faithfully
 flat.  

 Using Lemma~\ref{lem-invariants-neat} and
 Corollary~\ref{cor-invariants-neat} we deduce that the inclusions
 $A^G\xra{}A^H\xra{}A$ are neat.  A neat extension that is an
 algebraically free module is always topologically free, which
 proves~(a), (b) and~(c).  We have seen that $I^G$ and $I^H$ are
 finitely generated, and the rest of~(d) follows from
 Lemma~\ref{lem-invariants-neat}.
\end{proof}

\begin{proof}[Proof of Proposition~\ref{prop-topologies}]
 Claim~(a) is clear.  Claim~(b) follows from part~(a) of
 Lemma~\ref{lem-topologies} by passing to completions, and part~(c)
 follows immediately from~(a) and~(b).  

 For claim~(d), put
 \begin{align*}
  A' &= \{\al\in A\st \al_i<n \text{ for all } i\} \\
  B' &= A'/\Sg_r = 
        \{\bt\in B\st \bt_j=0 \text{ for all } j\geq n\}.
 \end{align*}
 For $\bt\in B'$ we put $H_\bt=\prod_i\Sg_{\bt_i}\leq\Sg_r$, so
 $\tau^{-1}\{\bt\}\simeq\Sg_r/H_\bt$.  As $A'$ is a basis for $R_r$
 over $\Rb_r$, we deduce that $S_r=R_r^{\Sg_r}$ is isomorphic to
 $\bigoplus_\bt\Rb_r^{H_\bt}$ as a module over $\Sb_r$.  It will thus
 suffice to show that $k\psb{y_1,\ldots,y_r}^{H_\bt}$ is topologically
 free of rank $|\Sg_r/H_\bt|$ over $k\psb{y_1,\ldots,y_r}^{\Sg_r}$,
 and this follows from part~(c) of Lemma~\ref{lem-topologies} by
 passing to completions.

 For part~(e), note that $R_r$ is finitely generated over $\Sb_r$ and
 thus is certainly finitely generated over the larger ring $S_r$.
 Neatness follows from Lemma~\ref{lem-invariants-neat}.

 Finally, we must show that for each of our rings there is a finitely
 generated ideal $J$ whose powers determine the topology.  For
 $\Rb_r$, we can obviously take $J$ to be the ideal
 $\Ib_r:=(y_1,\ldots,y_r)$.  Lemma~\ref{lem-invariants-neat} tells us
 that for $\Sb_r$ we can use the ideal
 $\Jb_r:=\Ib_r^{\Sg_r}=(u_1,\ldots,u_r)$.  For $S_r$ (which is
 topologically free over $\Sb_r$) we can therefore use the ideal
 $J_r:=\Jb_rS_r$.  Similarly, for $R_r$ we can use the ideal
 $I_r=\Ib_rR_r$.
\end{proof}

\begin{lemma}\label{lem-quotient-embeddable}
 If the curve $C$ is embeddable, then $R_r$ is topologically free of
 rank $r!$ over $S_r$.
\end{lemma}
\begin{proof}
 We may assume that 
 \[ C = \spf(k[x]^\wedge_{f(x)}) = \colim_m\spec(k[x]/f(x)^m) \]
 for some monic polynomial $f(x)$.  Put $A=k[x_1,\ldots,x_r]$, and
 give this the topology determined by the powers of the ideal
 $I=(f(x_1),\ldots,f(x_r))$, so $C^r=\spf(A^\wedge_I)$.  The evident
 action of $G:=\Sg_r$ on $A$ is continuous, and $A$ is free of rank
 $r!$ over $A^G$.  We see from Lemma~\ref{lem-invariants-neat} that
 $A^G$ is neat in $A$, so $A$ is topologically free over $A^G$ of rank
 $r!$, and the claim follows by passing to completions.
\end{proof}
\begin{lemma}\label{lem-proj-local}
 Let $A$ be a ring, $M$ a finitely generated $A$-module, and $B$ a
 faithfully flat $A$-algebra.  Suppose that $B\ot_AM$ is a free
 $B$-module of rank $s$.  Then $M$ is a projective $A$-module of the
 same rank.
\end{lemma}
\begin{proof}
 First, we claim that if $\mxi$ is a maximal ideal in $A$ with residue
 field $K=A/\mxi$, then $\dim_K(K\ot_AM)=s$.  Indeed, by faithful
 flatness there exists a prime ideal $\nxi\leq B$ with
 $\nxi\cap A=\mxi$.  Using Zorn's lemma we can find a maximal element
 of the set of all such ideals $\nxi$, and this is easily seen to be a
 maximal ideal in $B$.  It follows that the residue field $L=B/\nxi$
 is a field extension of $K$, so
 \[ \dim_K(K\ot_AM) = \dim_L(L\ot_KK\ot_AM) = \dim_L(L\ot_B(B\ot_AM)),
 \] 
 which is evidently equal to $s$.

 We now choose a finite generating set $\{m_1,\ldots,m_t\}$ for $M$.
 For each subset $S\sse\{1,\ldots,t\}$ with $|S|=s$, we let
 $f_S\:A^S\xra{}M$ be the map $\un{a}\mapsto\sum_sa_sm_s$, and we let
 $P_S$ and $Q_S$ be the kernel and cokernel of $f_S$.

 Next, we put $I_S=\ann(Q_S)\leq A$.  If $\mxi$ is maximal as before,
 we claim that there exists $S$ such that $I_S\not\leq\mxi$.  Indeed,
 as $\dim_K(K\ot M)=s$, we can certainly choose $S$ such that
 $\{m_i\st i\in S\}$ gives a basis for $K\ot_AM$.  It follows that
 $K\ot_AQ_S=0$, or equivalently that $\mxi Q_S=Q_S$.  The module $Q_S$
 is generated by the elements $m_j$ for $j\not\in S$, so we can find
 elements $u_{jk}\in\mxi$ for each $j,k\not\in S$ such that
 $m_j=\sum_ku_{jk}m_k$.  Let $U$ be the square matrix with entries
 $u_{jk}$ and put $u=\det(I-U)$.  As in~\cite{ma:crt}*{Theorem 2.1}, 
 we see that $u=1\pmod{\mxi}$ and $u\in I_S$, so $I_S\not\leq\mxi$ as
 claimed.

 It follows from this claim that $\sum_SI_S$ is not contained in any
 maximal ideal, so $\sum_SI_S=A$.  We can thus choose $a_S\in I_S$
 with $\sum_Sa_S=1$.  It follows that $\spec(A)$ is the union of the
 basic open subschemes $D(a_S)=\spec(A[a_S^{-1}])$.

 We have $a_SQ_S=0$ and so $Q_S[a_S^{-1}]=0$, so the map $f_S$ becomes
 surjective after inverting $a_S$.  It follows that the resulting map
 $1\ot f_S\:B[a_S^{-1}]^S\xra{}B[a_S^{-1}]\ot_AM$ is also surjective.
 Here both source and target are free modules of the same finite rank
 over $B[a_S^{-1}]$, so our map must in fact be an isomorphism.  As
 $B[a_S^{-1}]$ is faithfully flat over $A[a_S^{-1}]$, we deduce that
 $f_S$ actually gives an isomorphism
 $A[a_S^{-1}]^s\xra{}M[a_S^{-1}]$.  This shows that $M$ is locally
 free of rank $s$, and thus is projective.
\end{proof}

\begin{corollary}\label{cor-proj-local}
 Let $k$ be a ring, and let $A$ be a formal $k$-algebra whose topology
 is defined by the powers of a single open ideal $J$ (so
 $A=\invlim_mA/J^m$).  Let $M$ be a finitely generated $A$-module such
 that $M=\invlim_mM/J^mM$.  Let $k'$ be a faithfully flat
 $k$-algebra, and put $A'=k'\hot_kA$ and $M'=k'\hot_kM=A'\hot_AM$.
 Suppose that $M'$ is a free module of rank $s$ over $A'$; then $M$ is
 a projective module of rank $s$ over $A$.
\end{corollary}
\begin{proof}
 First, note that the map $A/J^m\xra{}A'/J^mA'=k'\ot_kA/J^m$ is a
 faithfully flat extension of discrete rings.  We can thus apply the
 lemma and deduce that $M/J^mM$ is a finitely generated projective
 module of rank $s$ over $A/J^m$.

 Next, as $M$ is finitely generated, we can choose an epimorphism
 $f\:A^t\xra{}M$ for some $t$.  Let $X_m$ be the set of $A$-module
 maps $g\:M/J^mM\xra{}(A/J^m)^t$ such that the induced map 
 \[ M/J^mM \xra{g} (A/J^m)^t \xra{f} M/J^mM \]
 is the identity.  As $M/J^mM$ is projective over $A/J^m$, we see that
 $X_m$ is nonempty.  There is an evident projection
 $\pi_m\:X_m\xra{}X_{m-1}$, which we claim is surjective.  Indeed,
 given $g\in X_{m-1}$ we can use the projectivity of $M/J^mM$ again to
 see that there exists a map $h\:M/J^mM\xra{}(A/J^m)^t$ lifting $g$.
 Let $\dl$ be the determinant of the resulting map
 $fh\:M/J^mM\xra{}M/J^mM$, so $\dl\in A/J^m$.  Because $g\in X_{m-1}$,
 we see that $\dl$ maps to $1$ in $A/J^{m-1}$.  As the kernel of the
 projection $A/J^{m-1}\xra{}A/J^m$ is nilpotent, it follows that $\dl$
 is a unit, so $fh$ is an isomorphism.  After replacing $h$ by
 $h(fh)^{-1}$ we may assume that $fh=1$, so $h\in X_n$ and
 $\pi(h)=g$.  It follows that $\invlim_mX_m\neq\emptyset$, and this
 gives a map $g\:M\xra{}A^t$ with $fg=1$.  Thus, $M$ is a retract of a
 free module, and hence is projective.  
\end{proof}

\begin{corollary}\label{cor-Rr-projective}
 $R_r$ is a projective module of rank $r!$ over $S_r$, so the
 projection $C^r\xra{}C^r/\Sg_r$ is a finite, faithfully flat map of
 degree $r!$.
\end{corollary}
\begin{proof}
 In Corollary~\ref{cor-proj-local}, we take $A=S_r$ and $M=R_r$.  We
 know from Proposition~\ref{prop-topologies} that the topology on
 $S_r$ is determined by powers of the ideal $J_r=(u_1,\ldots,u_r)$,
 and that $S_r$ is neat in $R_r$.  This means that the given topology
 on $R_r$ is determined by the ideals $J_r^mR_r$.  As $R_r$ is
 complete, we deduce that $R_r=\invlim_r R_r/J_r^mR_r$.  We next take
 $k'=\O_{S'}$ to be a faithfully flat extension of $k$ such that the
 curve $C'=S'\tm_SC$ is embeddable; this is possible by
 Corollary~\ref{cor-locally-embeddable}.  Using
 Lemma~\ref{lem-quotient-embeddable}, we see that $M'$ is
 topologically free of rank $r!$ over $A'$, so we can apply
 Corollary~\ref{cor-proj-local} and deduce that $R_r$ is projective
 over $S_r$.
\end{proof}

\section{Classification of divisors}
\label{sec-classify}

Our main task in this section is to prove the following result.
\begin{theorem}\label{thm-classify}
 Let $C$ be a formal multicurve over a scheme $S$.  Then for formal
 schemes $S'$ over $S$, there is a natural bijection between divisors
 of degree $r$ on $S'\tm_SC$ and maps $S'\xra{}C^r/\Sg_r$ over $S$.
\end{theorem}

\begin{construction}\label{cons-univ-div}
 We must first construct a universal example.  We start by putting
 $\Dl_i=\{(a_1,\ldots,a_r,b)\in C^{r+1}\st b=a_i\}$, which is a
 divisor of degree one on $C$ over $C^r$.  If we define
 $d_i(\un{a},b)=d(a_i,b)$ then $\O_{\Dl_i}=R_{r+1}/d_i$.  Now put
 $\dl_r=\prod_id_i\in R_{r+1}$ and
 $\tD_r=\sum_i\Dl_i=\spf(R_{r+1}/\dl_r)$, which is a divisor of degree
 $r$ on $C$ over $C^r$.  On the other hand, we note that
 $\dl_r\in R_{r+1}^{\Sg_r}=S_r\hot_kR$, so we can define
 $D_r=\spf((S_r\hot R)/\dl_r)$, which is a closed formal subscheme of
 $C^r/\Sg_r\tm_SC$.  It is clear that
 $R_r\ot_{S_r}\O_{D_r}=\O_{\tD_r}$, which is free of rank $r$ over
 $R_r$.  We know from Corollary~\ref{cor-Rr-projective} that $R_r$ is
 faithfully flat over $S_r$, and it follows from
 Lemma~\ref{lem-proj-local} that $\O_{D_r}$ is a projective module of
 rank $r$ over $S_r$.  Moreover, the relevant ideal is generated by
 the regular element $\dl_r$, so $D_r$ is a divisor on $C$ over
 $C^r/\Sg_r$. 
\end{construction}

Now put $Q_r=C^r/\Sg_r$ for brevity.  As in
Section~\ref{sec-symmetric}, we choose a topological basis $\{e_i\}$
for $\OC$, and use it to construct a topological basis
$\{e'_\bt\st\bt\in B\}$ for $\O_{Q_r}$.  We then put
$M=\spec(\Z[t_0,t_1,\ldots])$, and put $g=\sum_it_ie_i$, regarded as a
function on $M\tm Q_r\tm C$.  We then put
\[ h = N_{M\tm D_r/M\tm Q_r}(g) \in
     \O_{M\tm Q_r} = \prod_\bt \OS[t_i\st i\geq 0]e'_\bt.
\]
We claim that $h$ is actually equal to $\sum_\bt t^\bt e'_\bt$.
Indeed, although this is an equation in $\O_{M\tm Q_r}$, it will
suffice to prove it in the larger ring $\O_{M\tm C^r}$.  In that
context, we can describe $h$ as $N_{M\tm\tD_r/M\tm C^r}(g)$.  Now let
$\pi_j\:C^r\xra{}C$ be the $j$'th projection.  Using
Proposition~\ref{prop-full-set-norm} we see that
$h=\prod_j\pi_j^*g=\prod_j\sum_it_i\pi_j^*e_i$.  Expanding this out
gives 
\[ h = \sum_{\al\in A}t^{\tau(\al)} e_\al =
       \sum_{\bt\in B} \left(t^\bt \sum_{\tau(\al)=\bt} e_\al\right) = 
       \sum_{\bt\in B} t^\bt e'_\bt
\]
as claimed.

Now suppose we have a map $c\:S'\xra{}Q_r$ over $S$, and $D=c^*D_r$
over $S'$.  We deduce easily that 
\[ N_{M\tm D/M\tm S'}(g) = \sum_\bt t^\bt c^*(e'_\bt). \]
This shows that $c^*(e'_\bt)$ depends only on $D$, and
$\{e'_\bt\st\bt\in B\}$ is a topological basis for $S_r$, so the ring
map $c^*\:S_r\xra{}\O_{S'}$ depends only on $D$, so the map
$c\:S'\xra{}C^r/\Sg_r$ depends only on $D$.  We record this formally
as follows:

\begin{proposition}\label{prop-classify-inj}
 Suppose we have two maps $c_0,c_1\:S'\xra{}C^r/\Sg_r$ over $S$, and
 that $c_0^*D_r=c_1^*D_r$ as divisors over $S'$.  Then $c_0=c_1$. \qed
\end{proposition}

\begin{proof}[Proof of Theorem~\ref{thm-classify}]
 Let $S'$ be a scheme over $S$, and let $A$ be the set of maps
 $S'\xra{}C^r/\Sg_r$ over $S$, and let $B$ be the set of divisors of
 degree $r$ on $C$ over $S'$.  The construction $c\mapsto c^*D_r$
 gives a map $\phi\:A\xra{}B$, which is injective by
 Proposition~\ref{prop-classify-inj}.  To show that $\phi$ is
 surjective, suppose we have a divisor $D\in B$.  We can choose a
 faithfully flat map $q\:T\xra{}S'$ such that $q^*D$ has a full set of
 points, say $\un{u}=(u_1,\ldots,u_r)$.  We deduce that $q^*D$ is the
 pullback of $\tD_r$ along the map $T\xra{\un{u}}C^r$, and thus is the
 pullback of $D_r$ along the composite
 $c=(T\xra{\un{u}}C^r\xra{}C^r/\Sg_r)$.  Now let
 $q_0,q_1\:T\tm_{S'}T\xra{}T$ be the two projections, so $qq_0=qq_1$.
 Note that
 \[ (cq_0)^*D_r = q_0^*c^*D_r = q_0^*q^*D = (qq_0)^*D, \]
 and similarly
 \[ (cq_1)^*D_r = q_1^*c^*D_r = q_1^*q^*D = (qq_1)^*D. \]
 As $qq_0=qq_1$ we see that $(cq_0)^*D=(cq_1)^*D$, and so (by
 Proposition~\ref{prop-classify-inj}) we have $cq_0=cq_1$.  By
 faithfully flat descent, we have $c=\ov{c}q$ for a unique map
 $\ov{c}\:S'\xra{}C^r/\Sg_r$.  We then have
 $q^*\ov{c}^*D_r=c^*D_r=q^*D$, and using the faithful flatness of $q$,
 we deduce that $D=\ov{c}^*D_r=\phi(\ov{c})$.  This shows that $\phi$
 is also surjective, and thus a natural bijection.
\end{proof}

\begin{definition}\label{defn-divisor-sum}
 In the light of Theorem~\ref{thm-classify}, it makes sense to write
 $\Div_r^+(C)$ for $C^r/\Sg_r$.  The evident projection 
 \[ C^r/\Sg_r \tm_S C^s/\Sg_s = C^{r+s}/(\Sg_r\tm\Sg_s) 
    \xra{} C^{r+s}/\Sg_{r+s} 
 \]
 gives a map 
 \[ \sg_{r,s} \: \Div_r^+(C) \tm_S \Div_s^+(C) \xra{}
     \Div_{r+s}^+(C).
 \]
 It is easy to check that this classifies addition of divisors, in the
 following sense: if we have divisors $D=u^*D_r$ and $D'=v^*D_s$ on
 $C$ over $S'$, then $D+D'=w^*D_{r+s}$, where 
 \[ w = (S' \xra{(u,v)} \Div_r^+(C)\tm_S\Div_s^+(C)
         \xra{\sg} \Div_{r+s}^+(C)).
 \]
 We put $\Div^+(C)=\coprod_r\Div_r^+(C)$.  As one would expect, this
 is the free abelian monoid scheme generated by $C$; 
 see~\cite{st:fsfg}*{Section 6.2} for technical details.
\end{definition}
\begin{definition}\label{defn-convolution}
 Now suppose that $C$ has an abelian group structure, written
 additively.  We can then define
 $\tmu\:C^r\tm_SC^s\xra{}C^{rs}$ by
 \[ \tmu(a_0,\ldots,a_{r-1};b_0,\ldots,b_{s-1})_{i+rj} = a_i+b_j \]
 (for $0\leq i<r$ and $0\leq j<s$).  The composite 
 \[ C^r\tm_SC^s \xra{\tmu} C^{rs} \xra{q} C^{rs}/\Sg_{rs} \]
 is invariant under $\Sg_r\tm\Sg_s$, so we get an induced map 
 \[ \mu_{r,s} \: \Div_r(C) \tm_S \Div_s(C) \xra{} \Div_{rs}(C). \]
 If we have divisors $D=u^*D_r$ and $D'=v^*D_s$ on $C$ over $S'$, then
 we define $D*D'$ to be the divisor $w^*D_{rs}$, where
 \[ w = (S' \xra{(u,v)} \Div_r^+(C)\tm_S\Div_s^+(C)
         \xra{\mu} \Div_{rs}^+(C)).
 \]
 We call this the \emph{convolution} of $D$ and $D'$.  This operation
 makes $\Div^+(C)$ into a semiring.  If we have full sets of points,
 say $D=\sum_i[a_i]$ and $D'=\sum_j[b_j]$ then $D*D'$ is just
 $\sum_{i,j}[a_i+b_j]$.  
\end{definition}

\begin{proposition}\label{prop-containment-locus}
 Let $D$ and $D'$ be divisors on $C$ over $S$.  Then there exists a
 closed subscheme $T\sse S$ such that for any scheme $S'$ over $S$, we
 have $S'\tm_SD\leq S'\tm_SD'$ iff the map $S'\xra{}S$ factors through
 $T$.
\end{proposition}
\begin{proof}
 As $\OD$ is finitely generated and projective over $\OS$, we can
 choose an embedding $i\:\OD\xra{}\OS^N$ of $\OS$-modules, and a
 retraction $r\:\OS^N\xra{}\OD$.  We then have
 $i(f_{D'})=(a_1,\ldots,a_N)$ for some elements $a_j\in\OS$, and we
 put $J=(a_1,\ldots,a_N)$ and $T=\spec(\OS/J)$.  We find that a map
 $S'\xra{}S$ factors through $T$ iff $J$ maps to $0$ in $\O_{S'}$, iff
 $f_{D'}$ maps to $0$ in $\O_{S'}\ot_{\OS}\OD$, iff
 $S'\tm_SD\sse S'\tm_SD'$.
\end{proof}

\begin{proposition}\label{prop-SubrD}
 Let $D$ be a divisor on $C$ over $S$, and suppose that $r\geq 0$.
 Then there is a scheme $\Sub_r(D)$ over $S$ such that maps
 $S'\xra{}\Sub_r(D)$ over $S$ biject with divisors $D'\leq S'\tm_SD$
 of degree $r$.
\end{proposition}
\begin{proof}
 Over the formal scheme $\Div_r^+(C)$ we have both the originally
 given divisor $D$ and the universal divisor $D_r$.  We let
 $\Sub_r(D)$ denote the largest closed subscheme of $\Div_r^+(C)$ where
 $D_r$ is contained in $D$ (which makes sense by
 Proposition~\ref{prop-containment-locus}).  It is formal to check
 that this has the required property.
\end{proof}

\begin{proposition}\label{prop-PrD}
 Let $D$ be a divisor on $C$ over $S$, and suppose that $r\geq 0$.
 Then there is a scheme $P_r(D)$ over $S$ such that maps
 $S'\xra{}P_r(D)$ over $S$ biject with lists $(u_1,\ldots,u_r)$ of
 sections of $C$ over $S'$ such that $\sum_i[u_i]\leq S'\tm_SD$.
\end{proposition}
\begin{proof}
 Over the formal scheme $C^r$ we have both the originally given
 divisor $D$ and the divisor $\tD^r$.  We let $P_r(D)$ denote the
 largest closed subscheme of $C^r$ where $\tD_r$ is contained in $D$
 (which makes sense by Proposition~\ref{prop-containment-locus}).  It
 is formal to check that this has the required property.
\end{proof}
\begin{remark}\label{rem-Pr-top}
 Suppose that $D$ has degree $r$.  Then $P_r(D)$ classifies $r$-tuples
 for which $\sum_i[u_i]\leq D$, but by comparing degrees we see that
 this means that $\sum_i[u_i]=D$.  Thus, $P_r(D)$ classifies full sets
 of points for $D$.
\end{remark}

\begin{lemma}\label{lem-flat-flat}
 Suppose we have ring maps $A\xra{}B\xra{}C$, and $C$ is a projective
 module of degree $m>0$ over $B$, and also a projective module of
 degree $nm>0$ over $A$.  Then $B$ is a projective module of degree
 $n$ over $A$.
\end{lemma}
\begin{proof}
 We can use the second copy of $B$ to make $\Hom_A(B,B)$ into a
 $B$-module.  For any $B$-module $N$ there is an evident map
 $\Hom_A(B,B)\ot_BN\xra{}\Hom_A(B,N)$.  This is evidently an
 isomorphism if $N$ is a free module of finite rank, and thus (by
 taking retracts) also when $N$ is projective of finite rank over $B$.
 In particular, we have $\Hom_A(B,B)\ot_BC=\Hom_A(B,C)$.  As $C$ is
 also projective over $A$, the same kind of argument shows that 
 \[ \Hom_A(B,C)=\Hom_A(B,A)\ot_AC=(\Hom_A(B,A)\ot_AB)\ot_BC. \]
 It follows that $(\Hom_A(B,A)\ot_AB)\ot_BC=\Hom_A(B,B)\ot_BC$.  
 More precisely, there is a natural map
 \[ \al \: \Hom_A(B,A)\ot_A B\xra{} \Hom_A(B,B), \]
 given by $\al(\phi\ot b)(b')=\phi(b')b$.  By working through the
 above argument more carefully, we see that $\al\ot_B1_C$ is an
 isomorphism.  However, $C$ is faithfully flat over $B$ so $\al$
 itself must be an isomorphism.  In particular, we see that $1_B$ lies
 in the image of $\al$, so $1_B=\sum_{i=1}^N\al(\phi_i\ot b_i)$ for
 some maps $\phi_i\:B\xra{}A$ and elements $b_i\in B$.  This means
 that for all $b\in B$ we have $b=\sum_i\phi_i(b)b_i$.  We can use the
 elements $\phi_i$ to give a map $\phi\:B\xra{}A^N$, and the elements
 $b_i$ to give a map $\bt\:A^N\xra{}B$.  We find that $\bt\phi=1$,
 which proves that $B$ is projective.  It is now clear that the rank
 must be $n$.
\end{proof}

\begin{proposition}\label{prop-Pr-degree}
 Let $D$ be a divisor of degree $s$ on $C$ over $S$, and suppose that
 $0\leq r\leq s$.  Then there are natural maps
 $P_r(D)\xra{p}\Sub_r(D)\xra{q}S$ which are finite and very flat, with
 $\deg(p)=r!$ and $\deg(q)=s!/(r!(s-r)!)$ (so $\deg(qp)=s!/(s-r)!$).
\end{proposition}
\begin{proof}
 Over $P_r(D)$ we have tautological sections $u_1,\ldots,u_r$ of $C$
 giving a divisor $D'_r:=\sum_i[u_i]$ on $C$.  This is contained in
 $(qp)^*D$, so we can form the divisor $D''_r:=(qp)^*D-D'_r$, which
 has degree $s-r$ over $P_r(D)$.  It is easy to identify $P_{r+1}(D)$
 with $D''_r$, so $\deg(P_{r+1}(D)\xra{}S)=(s-r)\deg(P_r(D)\xra{}S)$.
 By an evident induction, we see that the map $pq$ is finite and very
 flat, with degree $s!/(s-r)!$, as claimed.

 Next, let $\Db$ be the tautological divisor of degree $r$ on $C$ over
 $\Sub_r(D)$.  We can then form the scheme $P_r(\Db)$, which
 classifies full sets of points on $\Db$.  As above, we see that the
 map $P_r(\Db)\xra{}\Sub_r(D)$ is finite and very flat, with degree
 $r!$.  We claim that $P_r(\Db)=P_r(D)$.  Indeed, a map
 $S'\xra{}P_r(\Db)$ over $S$ corresponds to a map $S'\xra{}\Sub_r(D)$,
 together with a lifting to $P_r(\Db)$.  Equivalently, it corresponds
 to a divisor of degree $r$ contained in $S'\tm_SD$, together with
 sections $u_1,\ldots,u_r\:S'\xra{}C$ giving a full set of points for
 that divisor.  The full set of points determines the divisor, so it
 is equivalent to just give sections $u_i$ with
 $\sum_i[u_i]\leq S'\tm_SD$, or equivalently, a map $S'\xra{}P_rD$
 over $S$.  The claim follows by Yoneda's lemma.  It follows that the
 map $p$ is finite and very flat, with degree $r!$.  We can now apply
 Lemma~\ref{lem-flat-flat} to see that $q$ is finite and very flat,
 with degree $s!/(r!(s-r)!)$.
\end{proof}

\begin{proposition}\label{prop-sub-univ}
 For the universal divisor $D_s$ over $\Div_s^+(C)$ we have
 \begin{align*} 
  \Sub_r(D_s) &= \Div_r^+(C) \tm_S \Div_{s-r}^+(C) \\
  P_r(D_s)    &= C^r \tm_S \Div_{s-r}^+(C).
 \end{align*}
\end{proposition}
\begin{proof}
 Let $S'$ be a scheme over $S$.  Then a map $S'\xra{}\Sub_r(D_s)$ over
 $S$ corresponds to a map $S'\xra{}\Div_r^+(C)$, together with a
 lifting to $\Sub_r(D_s)$.  Equivalently, it corresponds to a divisor
 $D$ of degree $s$ on $C$ over $S'$, together with a subdivisor
 $D'\leq D$ of degree $r$.  Given such a pair $(D,D')$, we have
 another divisor $D''=D-D'$, which has degree $s-r$.  There is
 evidently a bijection between pairs $(D,D')$ as above, and pairs
 $(D',D'')$ where $D'$ and $D''$ are arbitrary divisors of degrees $r$
 and $s-r$.  These pairs correspond in turn to maps
 $S'\xra{}\Div_r^+(C) \tm_S \Div_{s-r}^+(C)$ over $S$.  The first
 claim follows by Yoneda's lemma, and the second claim can be proved
 in the same way.
\end{proof}
\begin{corollary}\label{cor-sub-univ}
 We have $D_s=\Div_{s-1}^+(C)\tm_SC=C^s/\Sg_{s-1}$.
\end{corollary}
\begin{proof}
 Take $r=1$, and observe that $P_1(D_s)=\Sub_1(D_s)=D_s$ and
 $\Div_1^+(C)=C$.  
\end{proof}

%\section{\protect {Local structure of {$\Div_d^+(C)$}}}
\section{\protect {Local structure of the scheme of divisors}}
\label{sec-local-structure}

Let $C$ be a formal multicurve over a base $S$.  In the nonequivariant
case, we know that
$\Div_n^+(C)\simeq\spf(\OS\psb{c_1,\ldots,c_n})=\haf^n_S$, so
$\Div_n^+(C)$ is a formal affine space of dimension $n$ over $S$.
Equivariantly, this is not even true when $n=1$.  However, we will
show in this section that $\Div_n^+(C)$ is still a ``formal
manifold'', in the sense that the formal neighbourhood of any point is
isomorphic to $\haf^n_S$, at least up to a slight twisting.  Later we
will apply this to calculate $E^0BU(V)$, where $BU(V)$ is the
simplicial classifying space of the unitary group of a representation
$V$ of $A$.

We state the result more formally as follows.
\begin{theorem}
 Let $C=\spf(R)$ be a formal multicurve over $S=\spec(k)$, with a
 difference function $d$.  Let $s\:S\xra{}\Div_n^+(C)$ be a section,
 classifying a divisor $D=\spf(R/J)\subset C$.  Then the formal
 neighbourhood of $sS$ in $\Div_n^+(C)$ is isomorphic to the formal
 neighbourhood of zero in $\Map_S(D,\aff^1_S)$ (by an isomorphism that
 depends on the choice of $d$).
\end{theorem}
The rest of this section constitutes a more detailed explanation and a
proof.

We first examine the two formal schemes that are claimed to be
isomorphic.  We put $A_0=(R^{\hot n})^{\Sg_n}$ and
$X_0=\spf(A_0)=\Div_n^+(C)$.  The section $s$ corresponds to a
$k$-algebra map $A_0\xra{}k$, with kernel $K$ say.  We put
$A=(A_0)^\wedge_K$ and $X=\spf(A)$.  This is the formal neighbourhood
of $sS$ in $\Div_n^+(C)$.

Now consider the scheme $Y_0=\Map_S(D,\aff^1_S)$.  For any scheme $T$
over $S$, the maps $T\xra{}Y_0$ over $S$ are (essentially by
definition) the maps $D\tm_ST\xra{}\aff^1$ of schemes, or equivalently
the elements in the ring $\OD\ot_k\OT$.  These biject with the maps
$\OD^\vee=\Hom_k(\OD,k)\xra{}\OT$ of $k$-modules, or with the maps
$B_0=\Sym_k[\OD^\vee]\xra{}\OT$ of $k$-algebras.  Thus, we have
$Y_0=\spec(B_0)$.  We let $B$ be the completion of $B_0$ at the
augmentation ideal, and put $Y=\spf(B)$, which is the formal
neighbourhood of the zero section in $Y_0$.  Of course $B_0$ is just
the direct sum of all the symmetric tensor powers of $\OD^\vee$, and
$B$ is the direct product of the same terms.  If $\OD$ is free over
$k$ (rather than just projective) then $B$ is isomorphic to
$k\psb{c_1,\ldots,c_n}$; in the general case, it should be regarded
as a slight twist of this.  Note that maps $T\xra{}Y$ over $S$ biject
with $k$-linear maps $\OD^\vee\xra{}\Nil(\OT)$, or equivalently
elements of $\OD\ot_k\Nil(\OT)$.  Note also that a choice of
generators $x_1,\ldots,x_r$ for $\OD^\vee$ gives a split surjection 
$k\psb{x_1,\ldots,x_r}\xra{}B$.

There is an evident map 
\[ \al \: \Div_n^+(C) = X_0 \xra{} Y_0 = \Map_S(D,\aff^1_S), \]
sending the section $s'$ classifying a divisor $D'$ to the function
$(f_{D'})|_D\:D\xra{}\aff^1$.  This clearly sends $s$ itself to zero,
so it sends the formal neighbourhood of $s$ to the formal
neighbourhood of zero, so it gives a map $\al\:X\xra{}Y$.  We shall
show that this is an isomorphism.

Note that because $D$ and $D'$ have the same degree, we have $s'=s$ if
and only if $f_{D'}$ is divisible by $f_D$, if and only if
$\al(s')=0$.  This shows that the kernel $K$ of the map
$s^*\:A_0\xra{}k$ is generated by the image under $\al^*$ of the
augmentation ideal in $B_0$.  In particular, we see that $K$ is
finitely generated.

Because $\OD=R/J$ is projective over $k$, we can choose a
$k$-submodule $P\leq R$ such that $R=P\oplus J$.  It follows that the
map $P\xra{}R\xra{}\OD$ is an isomorphism, with inverse $\xi$ say.
\begin{lemma}
 Let $I\leq k$ be a finitely generated ideal with $I^m=0$, and let
 $g\in R$ be such that $g=f_D\pmod{IR}$.  Then $g$ is a regular
 element, the ideal $Rg$ is open, and we have $R=Rg\oplus P$.
\end{lemma}
\begin{proof}
 A standard topological basis for $R$ gives an isomorphism
 $R=\prod_ik$, and using the fact that $I$ is finitely generated we
 see that $I^jR=(IR)^j=\prod_iI^j$.  We thus have a finite filtration
 of $R$ with quotients $\prod_iI^j/I^{j+1}$.

 Now consider the $k$-linear self-map of $R$ given by
 $\lm(qf_D+r)=qg+r$ for $q\in R$ and $r\in P$.  This is easily seen
 to induce the identity map on the quotients of the above filtration,
 so it is an isomorphism.  It follows easily that $g$ is regular and
 $R=Rg\oplus P$.

 As $D$ is a divisor, we know that $Rf_D$ is open.  Thus, for any good
 parameter $y$ we have $y^l\in Rf_D$ for large $l$, say $y^l=uf_D$.
 We also know that $f_D=g+h$ for some $h\in IR$, so $y^l=uh\pmod{g}$.
 As $I^m=0$ we have $y^{lm}=u^mh^m=0\pmod{g}$, so $Rg$ is also open.
\end{proof}

We now define a map $\bt$ from sections of $Y$ to sections of $X$.  A
section of $Y$ is an element $r\in\Nil(k)\OD$.  As $\OD$ is finitely
generated we have $r\in I\OD$ for some finitely generated ideal
$I\leq\Nil(k)$, and by finite generation this satisfies $I^m=0$ for
some $m$.  We can thus apply the lemma to the function
$g=f_D+\xi(r)\in R$ and conclude that the subscheme $D'=\spec(R/g)$ is
a divisor of degree $n$, classified by a section $s'$ of $X_0$ say.
Over the subscheme $\spec(k/I)\subset S$ it clearly coincides with
$s$, so $(s')^*(K)\leq I$, so $(s')^*(K^m)=0$.  This shows that $s'$
is actually a section of $X$, as required.  We can thus define
$\bt(r)=s'$.  

In order to define a map $\bt\:Y\xra{}X$ of formal schemes over $S$,
we need to define maps $\bt_T$ from sections of $Y$ over $T$ to
sections of $X$ over $T$, naturally for all schemes $T$ over $S$.  For
this we just replace $C$ by $T\tm_SC$, $P$ by $\OT\ot_kP$ and follow
the same procedure.

We now define another map $\al'\:X\xra{}Y$.  It will again be
sufficient to do this for sections defined over $S$.  Let $s'$ be a
section of $X$, classifying a divisor $D'$.  Put $I=(s')^*K\leq k$;
this is finitely generated because $K$ is, and nilpotent because $s'$
lands in $X$.  Over $\spec(k/I)$ we have $D'=D$, so
$f_{D'}=f_D\pmod{IR}$.  The lemma tells us that $R=Rf_{D'}\oplus P$,
so there are unique elements $h\in R$ and $p\in P$ such that
$f_D=hf_{D'}-p$.  By reducing modulo $I$ we see that $h=1\pmod{IR}$
and $p\in IP$.  We let $r$ be the image of $p$ in $R/f_D=\OD$, so
$r\in I\OD$ and $\xi(r)=p$.  The map $\al'\:X\xra{}Y$ is defined by
$\al'(s')=r$.  Note that $h$ is invertible so $f_{D'}$ is a unit
multiple of $f_D+p=f_D+\xi(r)$, so $D'=\spec(R/(f_D+\xi(r)))=\bt(r)$.
This shows that $\bt\al'=1$.  

In the other direction, suppose we start with $r\in I\OD$ and put
$D'=\spec(R/(f_D+\xi(r)))$ (corresponding to $\bt(r)$).  There is then
a unique element $p\in P$ congruent to $-f_D$ modulo $f_{D'}$, and
$\al'\bt(r)$ is by definition the image of $p$ in $\OD$.  It is clear
that $-f_D$ is congruent to $\xi(r)$ modulo $f_D+\xi(r)$, which is a
unit multiple of $f_{D'}$, so $p=\xi(r)$ and $\al'\bt(r)=r$.  This
shows that $\al'\bt=1$, so $\al'$ and $\bt$ are isomorphisms.

We actually started by claiming that the (slightly more canonical) map
$\al$ is an isomorphism.  As $\bt$ is an isomorphism, it suffices to
check that the map $\al\bt\:Y\xra{}Y$ is an isomorphism, or that
$(\al\bt)^*$ is an automorphism of $\OY=B$.  As $B$ is the completed
symmetric algebra of a finitely generated projective module, it will
suffice to show that $(\al\bt)^*$ is the identity modulo the square of
the augmentation ideal.  By base-change to the universal case, it will
suffice to show that $\al\bt(r)=r$ whenever $r\in I\OD$ with $I^2=0$.
Given such an $r$, we form the divisor $D'=\spec(R/(f_D+\xi(r)))$
corresponding to $\bt(r)$, and observe that $f_{D'}=u(f_D+\xi(r))$ for
some $u\in R^\tm$.  As $f_{D'}=f_D\pmod{IR}$ we must have
$u=1\pmod{IR}$.  As $\xi(r)\in IR$ and $I^2=0$ we have
$u\xi(r)=\xi(r)$ and so $f_{D'}=\xi(r)\pmod{f_D}$, so $\al\bt(r)=r$ as
claimed. 

\section{Generalised homology of Grassmannians}
\label{sec-grassmann}

Consider a periodically orientable theory $E$ with associated
equivariant formal group $C=\spf(E^0P\CU)$ over $S=\spec(E^0)$.  Let
$G_r\CU$ be the space of $r$-dimensional subspaces of $\CU$, and put
$G\CU=\coprod_{r=0}^\infty G_r\CU$.  Here we reprove the following
result from~\cite{cogrkr:uec}.

\begin{theorem}[Cole, Greenlees, Kriz]\label{thm-grassmann}
 There are natural isomorphisms
 \begin{align*}
  E_*G_r\CU &= (E_*P\CU)^{\ot r}_{\Sg_r} \\
  E^*G_r\CU &= ((E^*P\CU)^{\hot r})^{\Sg_r} \\
  \spf(E^0G_r\CU) &= C^r/\Sg_r.
 \end{align*}
\end{theorem}

We first introduce some additional structure.  Rather than working
with spaces, it will be convenient to use (pre)spectra that happen to
be homotopy equivalent to suspension spectra of spaces.  Recall
from~\cite{mama:eos} that an \emph{orthogonal prespectrum} $X$
consists of a space $X(V)$ for each finite-dimensional inner product
space $V$, together with maps $\al_*\:X(V)\to X(W)$ for each isometric
isomorphism $\al\:V\to W$ and maps $\sg\:S^U\Smash X(V)\to X(U\op V)$
satisfying certain axioms.  The category of orthogonal prespectra has
a smash product, defined so that pairings 
$X(U)\Smash Y(V)\to Z(U\op V)$ satisfying some obvious axioms biject
with maps $X\Smash Y\to Z$.  There is a parallel theory of equivariant
unitary prespectra, where we use unitary representations of $A$ rather
than inner product spaces.  Our work in this section will be based on
that theory.

For any complex inner product space $V$, we put
\[ R_0(V)=\Sg^VGV_+=\bigWedge_r\Sg^VG_rV_+. \] 
Using the evident maps $G_rU\tm G_sV\xra{}G_{r+s}(U\op V)$ we get maps
$\mu_{U,V}\:R_0(U)\Smash R_0(V)\xra{}R_0(U\op V)$.  We also have
inclusions $\eta_U\:S^U=\Sg^UG_0U_+\xra{}R_0(U)$.  These maps make
$R_0$ into a commutative and associative ring in the category of
unitary prespectra.  All this works equivariantly in an obvious way.
The weak homotopy type of $R_0$ is
\[ R_0 \simeq \colim_{U\ll\CU} \Sg^{-U}\Sgi R_0(U) 
       = \colim_{U\ll\CU} \Sgip GU = \Sgip G\CU = 
        \Sgip\coprod_r G_r\CU.
\]
We write $Q_rR_0$ for the subfunctor $V\mapsto\Sg^VG_rV_+$, so that
$R_0=\bigWedge_rQ_rR_0$ and $Q_rR_0\simeq\Sgip G_r\CU$.  In particular,
we have $Q_0R_0=S^0$ and $Q_1R_0=\Sgip P\CU$.  This gives a map
$E_*P\CU\xra{}E_*R_0$ and thus a ring map
$\Sym_{E_*}E_*P\CU\xra{}E_*R_0$.  The theorem says that this is an
isomorphism.  For the proof, we need some intermediate spectra.  Let
$T$ denote the tautological bundle over $G_rV$.  For any
representation $W$, we put
\begin{align*}
 Q_rR_W(V) &= \Sg^V (G_rV)^{\Hom(T,W)} \\
 R_W(V)    &= \bigWedge_r Q_rR_W(V) = \Sg^V GV^{\Hom(T,W)}.
\end{align*}
This again gives a commutative unitary ring spectrum, with weak
homotopy type $R_W\simeq G\CU^{\Hom(T,W)}$.  In the case $W=0$ we
recover $R_0$ as before.  An inclusion $W\xra{}W'$ gives a ring map
$i\:R_W\xra{}R_{W'}$.  In particular, we have a ring map
$R_0\xra{}R_W$, whose fibre we denote by $J_W$.  This is weakly
equivalent to the stable fibre of the zero section
$G\CU_+\xra{}G\CU^{\Hom(T,W)}$, and thus is the sphere bundle of the
bundle $\Hom(T,W)$ over $G\CU$.  

Next, recall that there is an isometric embedding $\CU\op W\xra{}\CU$,
and that the space of such embeddings is connected.  We have
\[ Q_1R_W =  P\CU^{\Hom(T,W)} = P(\CU\op W)/PW \simeq P\CU/PW. \]
Using this, we have a diagram as follows, in which the rows are
cofibrations: 
\[ \xymatrix{
 {PW_+}    \dto \rto & 
 {P\CU_+}  \dto \rto &
 {P\CU/PW} \dto \\
 {J_W}          \rto &
 {R_0}          \rto &
 {R_W.}
 } \]
The map $P\CU/PW\xra{}R_W$ gives a ring map
\[ \tht_W\:\Sym_{E_*}E_*(P\CU,PW) \xra{} E_*R_W. \]

\begin{theorem}\label{thm-ERW}
 The above maps $\tht_W$ are isomorphisms.
\end{theorem}
The proof will be given after some preparatory results.

First, suppose we have representations $W$ and $L$ with $\dim(L)=1$.  We
put $W'=W\op L$ and investigate the fibre of the map
$R_W\xra{}R_{W'}$.  We may assume that $W'\leq\CU$, and then we have a
map
\[ S^{\Hom(L,W)} = PL^{\Hom(T,W)} \sse
   P\CU^{\Hom(T,W)} \xra{} R_W,
\]
which we denote by $b_{W,L}$.  Multiplication by $b_{W,L}$ gives a map
$\Sg^{\Hom(L,W)}R_W\xra{}R_W$, which we again denote by $b_{W,L}$.
(Note that this sends $\Sg^{\Hom(L,W)}Q_rR_W$ into $Q_{r+1}R_W$, or in
other words, it increases internal degrees by one.)
\begin{proposition}
 The sequence
 \[ \Sg^{\Hom(L,W)} R_W \xra{b_{W,L}} R_W \xra{} R_{W\op L} = R_{W'}
 \] 
 is a cofibration.
\end{proposition}
\begin{proof}
 This is a special case of the following fact.  Suppose we have a
 space $X$ with vector bundles $U$ and $V$.  Let $S(U)$ for the unit
 sphere bundle in $U$, and $D(U)$ for the unit disc bundle, so $X^U$
 is homeomorphic to $D(U)/S(U)$.  We can pull back $V$ along the
 projection $q\:S(U)\xra{}X$ and thus form the Thom space
 $S(U)^{q^*V}=(S(U)\tm_XD(V))/(S(U)\tm_XS(V))$.  It is not hard to see
 that there is a cofibration 
 \[ S(U)^{q^*V} \xra{} D(U)^{q^*V}\simeq X^V \xra{} X^{U\op V}. \]
 We will apply this with $X=G\CU$ and $V=\Hom(T,W)$ and $U=\Hom(T,L)$,
 so that $X^V=R_W$ and $X^{U\op V}=R_{W'}$.  To prove the proposition,
 we need to identify $S(U)^{q^*V}$ with $\Sg^{\Hom(L,W)}R_W$.  

 To do this, observe that $S(U)$ is the space of pairs $(M,\al)$
 where $M$ is a finite-dimensional subspace of $\CU$ and
 $\al\:M\xra{}L$ is a linear map of norm one.  As $L$ has dimension
 one, we find that $\al$ can be written as the composite of the
 orthogonal projection $M\xra{}M\ominus\ker(\al)$ with an isometric
 isomorphism $M\ominus\ker(\al)\xra{}L$.  Using this, we identify
 $S(U)$ with the space of pairs $(N,\bt)$, where $N$ is a
 finite-dimensional subspace of $\CU$ and $\bt\:L\xra{}\CU\ominus N$
 is an isometric embedding; the correspondence is that $M=N\op\bt(L)$
 and 
 \[ \al = (N\op\bt(L) \xra{\text{proj}} \bt(L) \xra{\bt^{-1}} L). \]
 We can thus define a map $k\:S(U)\xra{}G\CU$ by $k(N,\bt)=N$ (or
 equivalently, $k(M,\al)=\ker(\al)$).  This makes  $S(U)$ into an
 equivariant fibre bundle over $G\CU$.  The fibre over a point
 $N\in G\CU$ is the space $\CL(L,\CU\ominus N)$, which is well-known
 to be contractible, and the contraction can be chosen to be
 equivariant with respect to the stabilizer of $N$ in $A$.  It follows
 that $k$ is an equivariant equivalence.  We next analyse the inverse
 of $k$ (which will help to make the above argument more explicit).

 Recall that $L\leq W'\leq \CU$, so we can put
 \[ Y = \{N \in G\CU \st N \text{ is orthogonal to } L\}
     = G(\CU\ominus L).
 \]
 Define $j\:Y\xra{}X$ by $j(N)=N\op L$, and then define
 $\tj\:Y\xra{}S(U)$ by  
 \[ \tj(N) = (N\op L,\text{proj}\:N\op L\xra{}L), \]
 so $q\tj=j$.  Clearly $k\tj\:Y=G(\CU\ominus L)\xra{}G\CU=X$ is just
 the map induced by the inclusion $\CU\ominus L\xra{}\CU$.  As the
 space of linear isometries between any two complete $A$-universes is
 equivariantly contractible, we see that this inclusion is an
 equivariant equivalence.  As the same is true of $k$, we deduce that
 $\tj$ is also an equivariant equivalence.  We can thus identify
 $S(U)$ with $Y$ and $q\:S(U)\xra{}X$ with $j\:Y\xra{}X$.  It follows
 that we can identify $q^*V$ with $j^*V$, but the fibre of $j^*V$ over
 a point $N\in Y$ is $\Hom(j(N),W)=\Hom(L,W)\op\Hom(N,W)$, so 
 \[ Y^{j^*V}=\Sg^{\Hom(L,W)}Y^{\Hom(T,W)}\simeq\Sg^{\Hom(L,W)}R_W. \]
 This gives a cofibration $\Sg^{\Hom(L,W)}R_W\xra{}R_W\xra{}R_{W'}$,
 and one can check from the definitions that the first map is just
 multiplication by $b_{W,L}$.
\end{proof}

Now choose a complete flag 
\[ 0 = W_0 < W_1 < \ldots < \CU, \]
where $\dim_\C W_i=i$ and $\CU=\colim_iW_i$.  Put $R(i)=R_{W_i}$, so
we have maps
\[ R_0 = R(0) \xra{} R(1) \xra{} R(2) \xra{} \ldots. \]
Put $L_i=W_{i+1}\ominus W_i$ and $U_i=\Hom(L_i,W_i)$ and
$b_i=b_{W_i,L_i}$, so we have a cofibration
\[ \Sg^{U_i} R(i) \xra{b_i} R(i) \xra{} R(i+1). \]

\begin{lemma}
 Suppose that $B\leq A$, and split $W$ as
 $\bigoplus_{\bt\in B^*}W[\bt]$ in the usual way.  Then
 \[ \phb^B R_W = \bigSmash_{\bt\in B^*} R_{W[\bt]}, \]
 where as before
 \[ W[\bt] =
    \{ w \in W \st bw = e^{2\pi i\bt(b)}w \text{ for all } b \in B\},
 \]
 and so the connectivity of $(\phb^BR_W)/S^0$ is at least
 $\min_\bt(2\dim_\C(W[\bt])-1)$.
\end{lemma}
\begin{proof}
 We have 
 \[ \phb^BG\CU = (G\CU)^B =
     \{\text{ $B$-invariant subspaces of $\CU$ }\}.
 \]
 Any complex representation $U$ of $B$ splits as
 $\bigoplus_\bt U[\bt]$, so a subspace $U\leq\CU$ is invariant iff it
 is the direct sum of its intersections with the subspaces
 $\CU[\bt]$.  It follows that 
 \[ (G\CU)^B = \prod_\bt G\CU[\bt] \simeq \prod_\bt G\C^\infty. \]
 We have a tautological bundle $T[\bt]$ over $G\CU[\bt]$, and the
 bundle $\Hom_{\C B}(T,W)$ over $(G\CU)^B$ is the external direct sum
 of the bundles $\Hom_\C(T[\bt],W[\bt])$.  The Thom complex
 $G\CU[\bt]^{\Hom(T[\bt],W[\bt])}$ is just $R_{W[\bt]}$, and it follows
 that $\left[(G\CU)^{\Hom(T,W)}\right]^B$ is just the smash product of
 these factors, as claimed.

 For the last statement, note that if $X$ is a space and $U$ is a
 vector bundle of real dimension $d$ over $X$, then $X^U$ is always
 $(d-1)$-connected.  Now let $\CV$ be a complex universe, and $V$ a
 complex vector space of finite dimension $d$.  The bundle $\Hom(T,V)$
 over $G_r\CV$ has real dimension $2rd$, so $\conn(Q_rR_V)\geq 2rd-1$,
 and 
 \[ \conn(R_V/S^0)=\conn(\bigWedge_{r>0}Q_rR_V)\geq 2d-1. \]
 The claim follows easily.
\end{proof}

\begin{corollary}\label{cor-colimit}
 $\colim_i R(i)=S^0$.
\end{corollary}
\begin{proof}
 The unit map $S^0\xra{}Q_0R(i)$ is an isomorphism for all $i$, so
 $\colim_iQ_0R(i)=S^0$.  It will thus suffice to show that
 $\colim_iR(i)/S^0=0$, or equivalently that the spectrum
 $\phb^B(\colim_iR(i)/S^0)=\colim_i((\phb^BR(i))/S^0)$ is
 nonequivariantly contractible for all $B$.  As $\CU$ is a complete
 universe, we have $\dim W_i[\bt]\xra{}\infty$ as $i\xra{}\infty$ for
 all $\bt$, so $\conn(\phb^BR(i)/S^0)\xra{}\infty$, and the claim
 follows. 
\end{proof}

We now let $E$ be a periodically oriented theory, with orientation $x$
say.  This gives a universal generator $u_i$ for $\tE_0S^{U_i}$, and a
basis $\{c_i\st i\geq 0\}$ for $\tE_0P\CU$.  Put
\[ ER(i) = \Sym_{E_*} E_*(P\CU,PW_i) = E_*[c_j\st j\geq i] =
    ER(0)/(c_k\st k<i),
\]
and let $Q_rER(i)$ be the submodule generated by monomials of weight
$r$ (where each generator $c_j$ is considered to have weight one).
We then have maps
\[ \tht_i = \tht_{W_i}\: ER(i) \xra{} E_*R(i), \]
which restrict to give maps
\[ \tht_{ir} \: Q_rER(i) \xra{} E_*Q_rR(i). \]
The elements $u_i$ and $c_i$ are related as follows: the inclusion
$PL_i\xra{}P\CU$ gives an inclusion
$S^{U_i}\xra{}P\CU^{\Hom(T,W_i)}\simeq P\CU/PW_i$, and the image of
$u_i$ under this map is the same as the image of $c_i$ under the
evident quotient map $P\CU\xra{}P\CU/PW_i$.  It follows that the
cofibration $\Sg^{U_i}R(i)\xra{b_i}R(i)\xra{}R(i+1)$ gives rise to a
cofibration
\[ E\Smash R(i) \xra{c_i} E\Smash R(i) \xra{} E\Smash R(i+1), \]
which restricts to give a cofibration 
\[ E\Smash Q_{r-1}R(i) \xra{c_i}
   E\Smash Q_rR(i)     \xra{}
   E\Smash Q_rR(i+1).
\]

\begin{proposition}\label{prop-ERW}
 The maps $\tht_{ir}$ are isomorphisms for all $i$ and $r$.
\end{proposition}
\begin{proof}
 The maps $\tht_{j0}$ and $\tht_{j1}$ are visibly isomorphisms, so we
 may assume inductively that $\tht_{j,r-1}$ is an isomorphism for all
 $j$.  The cofibration displayed above gives a diagram $D(i)$ as
 follows, with exact rows:
 \[ \xymatrix{
  {Q_{r-1}ER(i)}
   \ar@{ >->}[r]^{c_i}
   \dto_{\tht_{i,r-1}}^{\simeq} &
  {Q_rER(i)} \ar@{->>}[r] \dto_{\tht_{i,r}} &
  {Q_rER(i+1)} \dto^{\tht_{i+1,r}} \\
  {E_*Q_{r-1}R(i)} \rto_{c_i} &
  {E_*Q_rR(i)} \rto_{q_{ir}} &
  {E_*Q_rR(i+1)}
  } \]
 We first prove that $\tht_{ir}$ is surjective for all $i$.  Let
 $\Tht(i)$ be the image of $\tht_{ir}$, so the claim is that
 $\Tht(i)=E_*Q_rR(i)$.  For $j\geq i$ we write $K(j)$ for the kernel
 of the map $E_*Q_rR(i)\xra{}E_*Q_rR(j)$.  Clearly
 $K(i)=0\leq\Tht(i)$.  Suppose that $K(j)\leq\Tht(i)$, and consider
 the diagram $D(j)$.  Suppose that $a\in K(j+1)$.  Let $b$ be the
 image of $a$ in $E_*Q_rR(j)$, so $q_{jr}(b)=0$.  As the bottom row of
 $D(j)$ is exact and $\tht_{j,r-1}$ is an isomorphism we see that
 there exists $c\in Q_rER(j)$ with $\tht_{jr}(c)=b$.  Moreover, the
 map $Q_rER(i)\to Q_rER(j)$ is epi, so we can lift $c$ to get an
 element $d\in Q_rER(i)$.  Now $a-\tht_{ir}(d)\in K(j)\leq\Tht(i)$ and
 visibly $\tht_{ir}(d)\in\Tht(i)$ so $a\in\Tht(i)$.  It follows by
 induction that $K(j)\leq\Tht(i)$ for all $j$.  Moreover,
 Corollary~\ref{cor-colimit} implies that $E_*Q_rR(i)$ is the union of
 the subgroups $K(j)$, so $\Tht(i)=E_*Q_rR(i)$ as required.

 We now see that in $D(i)$, the vertical maps are surjective, so
 $q_{ir}$ is surjective.  As the bottom row is part of a long exact
 sequence and the right hand map is surjective, we conclude that the
 bottom row is actually a short exact sequence.  Using the snake
 lemma, we conclude that the induced map
 $\ker(\tht_{ir})\xra{}\ker(\tht_{i+1,r})$ is an isomorphism.  It
 follows that for any $m>i$, the map
 $\ker(\tht_{ir})\xra{}\ker(\tht_{mr})$ is an isomorphism.  However,
 we have $\ker(\tht_{ir})\leq Q_rER(i)$, and it is also clear that
 when $r>0$, any element of $Q_rER(i)$ maps to zero in $Q_rER(m)$ for
 $m\gg 0$.  It follows that $\ker(\tht_{ir})$ must be zero, so
 $\tht_{ir}$ is an isomorphism as claimed.
\end{proof}

\begin{proof}[Proof of Theorem~\ref{thm-ERW}]
 Given any subrepresentation $W<\CU$, we can choose our flag $\{W_i\}$
 such that $W=W_i$ for some $i$.  The theorem then follows from
 Proposition~\ref{prop-ERW}.  
\end{proof}

\section{Thom isomorphisms and the projective bundle theorem}
\label{sec-thom}

Let $E$ be a periodically orientable $A$-equivariant cohomology
theory, with associated equivariant formal group $(C,\phi)$ over $S$.
For any $A$-space $X$, we will write $X_E=\spf(E^0X)$.  

Now let $V$ be an equivariant complex vector bundle over $X$.  We
write $PV$ for the associated bundle of projective spaces, and $X^V$
for the Thom space (so $X^V=P(V\op\C)/PV$).  In this section, we will
give a Thom isomorphism and a projective bundle theorem to calculate
$\tE^*X^V$ and $E^*PV$.

First, it is well-known that equivariant bundles of dimension $r$ over
$X$ are classified by homotopy classes of $A$-maps $X\xra{}G_r\CU$.
(See for example~\cite{at:kt}*{Section 1.6}.)  We saw above that
$E^0G_r\CU=S_r$, and moreover the standard topological basis
$\{e'_\bt\}$ for $S_r$ is dual to a universal basis for $E_0G_r\CU$.
It follows that $(G_r\CU)_E=C^r/\Sg_r=\Div_r^+(C)$.

Now let $T$ denote the tautological bundle over $G_r\CU$.  It is not
hard to identify the projective bundle $PT\xra{}G_r\CU$ with the
addition map 
\[ G_{r-1}\CU\tm P\CU=G_{r-1}\CU\tm G_1\CU \xra{} G_r\CU, \]
and thus to identify $E^0PT$ with $S_{r-1}\hot R=\O_{C^r/\Sg_{r-1}}$,
so $PT_E=C^r/\Sg_{r-1}$.  On the other hand, we can use
Corollary~\ref{cor-sub-univ} to identify $C^r/\Sg_{r-1}$ with the
universal divisor $D_r$ over $C^r/\Sg_r$.

Now suppose we have a vector bundle $V$ over $X$, classified by a map
$c\:X\xra{}G_r\CU$, so $c^*T\simeq V$.  The map $c$ is then covered by
a map $\tc\:PV\xra{}PT$, which gives a map
$\tc^*\:\O_{D_r}=E^0PT\xra{}E^0PV$.  We can combine this with the
evident map $E^*X\xra{}E^*PV$ to get a map
\[ \tht_{X,V} \: \O_{D_r}\ot_{S_r} E^*X \xra{} E^*PV. \]

\begin{theorem}\label{thm-proj-bundle}
 For any $X$ and $V$ as above, the map $\tht_{X,V}$ is an isomorphism
 (and so $E^*PV$ is a projective module of rank $r$ over $E^*X$).
\end{theorem}

\begin{proof}
 We first examine the simplifications that occur when $V$ admits a
 splitting $V=L_1\op\ldots\op L_r$, where each $L_i$ is a line bundle.
 In this case, the classifying map $X\xra{}G_r\CU$ factors through
 $P\CU^r$, so the map $S_r\xra{}E^0X$ factors through $R_r$.  As
 $\O_{D_r}\ot_{S_r}R_r=\O_{\tD_r}$, we see that $\tht_{X,V}$ is the
 composite of an isomorphism with a map
 \[ \tht'_{X,V} \: \O_{\tD_r}\ot_{R_r} E^*X \xra{} E^*PV. \]
 Next, choose a coordinate $x$ on $C$ and define a difference function
 $d(a,b)=x(b-a)$ as usual.  Define a function $d_i$ on $C^{r+1}$ by 
 \[ d_i(a_1,\ldots,a_r,b) = d(a_i,b) = x(b-a_i), \]
 as in Construction~\ref{cons-univ-div}.  We then put
 $c_i=\prod_{j<i}d_j$.  By the method of
 Proposition~\ref{prop-top-basis} we see that $\{c_1,\ldots,c_r\}$ is
 a basis for $\O_{\tD}$ over $R_r$, so $\tht'_{X,V}$ is just the map
 $(E^*X)^r\xra{}E^*PV$ given by $(t_1,\ldots,t_r)\mapsto\sum_it_ic_i$.  

 Now consider the case where $X$ is a point, so $V$ is just a
 representation of $A$.  In this case there is always a splitting
 $V=L_1\op\ldots\op L_r$ as above, where $L_i=L_{\al_i}$ for some
 $\al_i\in A^*$.  In this case the image of $c_i$ in $E^0PV$ is just
 the element $x_{U_{i-1}}$ from Corollary~\ref{cor-flag-basis}, so the
 map $\tht_{X,V}$ is an isomorphism.

 More generally, suppose that $X$ is arbitrary but $V$ is a constant
 bundle, with fibre given by a representation $W=\sum_iL_{\al_i}$ say.
 As the elements $x_{U_i}$ form a universal basis for $E^*PW$, we see
 that $E^*PV=E^*X\ot_{E^*}E^*PW=\bigoplus_iE^*.x_{U_i}$ and it follows
 easily that $\tht_{X,V}$ is again an isomorphism.

 Now consider the case $X=A/B$ for some $B\leq A$.  It is easy to see
 that any bundle over $X$ has the form $A\tm_BW_0$ for some
 representation $W_0$ of $B$.  However, as $A$ is a finite abelian
 group, we can find a representation $W$ of $A$ such that $W|_B=W_0$,
 and it follows that $A\tm_BW_0$ is isomorphic to the constant bundle
 $A/B\tm W$.  It follows that $\tht_{X,V}$ is again an isomorphism.

 Now let $X$ and $V$ be arbitrary, and suppose we can decompose $X$ as
 the union of two open sets $X_0$ and $X_1$, with intersection $X_2$.
 Suppose we know that the maps $\tht_{X_i,V}$ are isomorphisms for
 $i=0,1$ and $2$; we claim that $\tht_{X,V}$ is also an isomorphism.
 Indeed, the decomposition gives a Mayer-Vietoris sequence involving
 $E^*X$.  We can tensor this by the projective module $\O_{D_r}$ over
 $S_r$, and it will remain exact.  Alternatively, we can pull back the
 decomposition to get a decomposition of $PV$, and obtain another
 Mayer-Vietoris sequence.  The $\tht$ maps are easily seen to be
 compatible with these sequences, so the claim follows by the five
 lemma.

 We can now now prove that $\tht_{X,V}$ is an isomorphism for all $X$
 and $V$, by induction on the number of cells and passage to colimits.
\end{proof}
\begin{corollary}\label{cor-proj-bundle}
 If $V$ is an equivariant vector bundle over $X$, then the formal
 scheme $D(V):=PV_E$ is a divisor on $C$ over $X_E$, of degree equal
 to the dimension of $V$.  Moreover, we have $D(V\op W)=D(V)+D(W)$ and
 $D(V\ot W)=D(V)*D(W)$.
\end{corollary}
\begin{proof}
 The first statement is clear from the theorem.  We need only check
 the equation $D(V\op W)=D(V)+D(W)$ in the universal case, where
 $X=G_r\CU\tm G_s\CU$.  As the map $P\CU^{r+s}\xra{}G_r\CU\tm G_s\CU$
 induces a faithfully flat map $C^{r+s}\xra{}C^r/\Sg_r\tm_SC^s/\Sg_s$,
 it suffices to check that equation for the obvious bundles
 $V=L_1\op\dotsb\op L_r$ and $W=L_{r+1}\op\dotsb\op L_{r+s}$ over
 $P\CU^{r+s}$, in which case it is clear.  A similar approach works
 for convolution of divisors.
\end{proof}

We next consider the Thom isomorphism. 
\begin{definition}\label{defn-thom-module}
 Let $C$ be a formal multicurve group over $S$, with zero section
 $\zt\:S\xra{}C$.  Given a divisor $D$ on $C$ over $S$, we let $J_D$
 denote the kernel of the restriction map $\OC\xra{}\OD$, which is a
 free module of rank one over $\OC$.  We also use the map
 $\zt^*\:\OC\xra{}\OS$ to make $\OS$ into a module over $\OC$, and we
 define $L(D)=\OS\ot_{\OC}J_D$, which is a free module of rank one
 over $\OS$, or equivalently a trivialisable line bundle over $S$.  We
 call this the \emph{Thom module} for $D$.  More generally, given a
 scheme $S'$ over $S$ and a divisor $D$ on $C$ over $S'$, we obtain a
 trivialisable line bundle $L(D)$ over $S'$.
\end{definition}
\begin{remark}\label{rem-thom-kernel}
 Note that $\ker(\zt^*)=J_{[0]}$ and that $J_DJ_{[0]}=J_{D+[0]}$.  It
 follows that 
 \[ L(D) = J_D/J_{D+[0]} = \ker(\O_{D+[0]} \xra{} \O_D). \]
\end{remark}
\begin{remark}\label{rem-thom-generator}
 If we fix a coordinate $x$ and put $d(a,b)=x(b-a)$, we get a
 generator $f_D$ for $J_D$ as in Definition~\ref{defn-fD}, and thus a
 generator $u_D=1\ot f_D$ for $L(D)$, which we call the \emph{Thom
   class}.  However, these generators are not completely canonical
 because of the choice of coordinate.

 We also define the \emph{Euler class} $e_D$ to be the element
 $f_D(0)=\zt^*f_D\in\OS$.  Note that if $D=[u]$ for some section $u$,
 then $f_D(a)=x(a-u)$ and so $e_D=x(-u)=\xb(u)$.
\end{remark}
\begin{remark}\label{rem-thom-products}
 For any two divisors $D$ and $D'$, we have $J_{D+D'}=J_DJ_{D'}$,
 which can be identified with $J_D\ot_{\OC}J_{D'}$ (because $J_D$ and
 $J_{D'}$ are each generated by a single regular element).  It follows
 that $L(D+D')=L(D)\ot_{\OS}L(D')$.  In terms of a coordinate, we have
 $u_{D+D'}=u_D\ot u_{D'}$ and $e_{D+D'}=e_De_{D'}$.
\end{remark}

\begin{theorem}\label{thm-thom}
 Let $V$ be an equivariant complex bundle over a space $X$, giving a
 divisor $D(V)=PV_E$ on $C$ over $X_E$ as in
 Corollary~\ref{cor-proj-bundle} and thus a free rank one module
 $L(D(V))$ over $E^0X$.  Then there is a natural isomorphism
 $\tE^0X^V=L(D(V))$ (and $\tE^*X^V=L(D(V))\ot_{E^0X}E^*X$).  Moreover,
 if we choose a coordinate and thus obtain a Thom class
 $u_{D(V)}$ as in Remark~\ref{rem-thom-generator}, then this gives a
 \emph{universal} generator for $\tE^0X$.
\end{theorem}
\begin{proof}
 Consider the cofibration $P(V)\xra{}P(V\op\C)\xra{}X^V$.  Using
 Theorem~\ref{thm-proj-bundle} we see that
 \begin{align*}
  E^*P(V)      &= E^*X\ot_{E^0X}\O_{D(V)} \\
  E^*P(V\op\C) &= E^*X\ot_{E^0X}\O_{D(V\op\C)}
                = E^*X\ot_{E^0X}\O_{D(V)+[0]}.
 \end{align*}
 As the map $\rho\:\O_{D(V)+[0]}\xra{}\O_{D(V)}$ is a split surjection
 of $E^0X$-modules, we see that the long exact sequence of the
 cofibration splits into short exact sequences.  As
 $\ker(\rho)=L(D(V))$, we see that $\tE^*X^V=L(D(V))\ot_{E^0X}E^*X$.
 By looking in degree zero, we see that $\tE^0X^V=L(D(V))$.  As this
 isomorphism is natural in $X$, it is easy to see that the generator
 is universal.
\end{proof}

\begin{remark}\label{rem-thom-product}
 If we have two bundles $V$ and $V'$, the above results give
 \begin{align*}
  \tE^0X^{V\op V'} &= L(D(V\op V')) = L(D(V) + D(V')) \\ 
  &= L(D(V)) \ot_{E^0X} L(D(V')) = \tE^0X^V \ot_{E^0} \tE^0X^{V'}.
 \end{align*}
 One can check that this isomorphism
 $\tE^0X^V\ot\tE^0X^{V'}=\tE^0X^{V\op V'}$ is induced by the usual
 diagonal map $X^{V\op V'}\xra{}X^V\Smash X^{V'}$.
\end{remark}

\begin{definition}\label{defn-thom-class}
 We write $u_V$ for $u_{D(V)}$, and call this the Thom class of $V$.
 We also write $e_V$ for $e_{D(V)}$, and call this the Euler class of
 $V$.  (Using Remark~\ref{rem-thom-generator}, we see that this is
 consistent with the definition for line bundles given in
 Section~\ref{sec-orientability}.) 
\end{definition}
It is easy to see that the Euler class is the pullback of the Thom
class along the zero section $X\xra{}X^V$, and that 
$e_{V\op W}=e_Ve_W$.

Now suppose that $r\leq\dim(V)$, and consider the space $P_r(V)$
consisting of all tuples $(x;L_1,\ldots,L_r)$ where $x\in X$ and the
$L_i$ are orthogonal lines in $V_x$.  Recall also that $P_r(D(V))$ is
the classifying scheme for $r$-tuples $(u_1,\ldots,u_r)$ of sections
of $C$ such that $\sum_i[u_i]\leq D(V)$, as in
Proposition~\ref{prop-PrD}.

\begin{proposition}\label{prop-Pr}
 There is a natural isomorphism $P_r(V)_E=P_r(D(V))$.
\end{proposition}
\begin{proof}
 For each $i$ we have a line bundle over $P_r(V)$ whose fibre over
 $(x,L_1,\ldots,L_r)$ is $L_i$.  This is classified by a map
 $P_r(V)\xra{}P\CU$, which gives rise to a map $u_i\:P_r(V)_E\xra{}C$.
 The direct sum of these line bundles corresponds to the divisor
 $[u_1]+\ldots+[u_r]$.  This direct sum is a subbundle of $V$, so
 $[u_1]+\ldots+[u_r]\leq D(V)$.  This construction therefore gives us
 a map $P_r(V)_E\xra{}P_r(D(V))$.

 In the case $r=1$ we have $P_1(V)=PV$ and $P_1(D(V))=D(V)$ so the
 claim is that $(PV)_E=D(V)$, which is true by definition.  In
 general, suppose we know that $P_{r-1}(V)_E=P_{r-1}(D(V))$.  We can
 regard $P_r(V)$ as the projective space of the bundle over
 $P_{r-1}(V)$ whose fibre over a point $(x,L_1,\ldots,L_{r-1})$ is the
 space $V_x\ominus(L_1\op\ldots\op L_{r-1})$.  It follows that
 $P_r(V)_E$ is just the divisor $D(V)-([u_1]+\ldots+[u_{r-1}])$ over
 $P_{r-1}(D(V))$, which is easily identified with $P_r(D(V))$.  The
 proposition follows by induction.
\end{proof}

We next consider the Grassmannian bundle 
\[ G_r(V)=\{(x,W)\st x\in X\;,\; W\leq V_x \text{ and }\dim(W)=r\}.
\]
\begin{proposition}\label{prop-Gr}
 There is a natural isomorphism $G_r(V)_E=\Sub_r(D(V))$.
\end{proposition}
\begin{proof}
 Let $T$ denote the tautological bundle over $G_r(V)$.  This is a rank
 $r$ subbundle of the pullback of $V$ so we have a degree $r$
 subdivisor $D(T)$ of the pullback of $D(V)$ over $G_r(V)_E$.  This
 gives rise to a map $f\:G_r(V)_E\xra{}\Sub_r(D(V))$, so if we put
 $A=\O_{\Sub_r(D(V))}$ we get a map $f^*\:A\xra{}E^0G_r(V)$, and we
 must show that this is an isomorphism.  Now consider the tautological
 divisor $\Db$ of degree $r$ over $\Sub_r(D(V))$.  As the module
 $B=\O_{P_r\Db}$ is faithfully flat over $A$, it will suffice to show
 that the map $f^*\:B\xra{}B\ot_AE^0G_r(V)$ is an isomorphism.
 However, we saw in the proof of Proposition~\ref{prop-Pr-degree} that
 $P_r\Db=P_rD(V)=(P_rV)_E$, so $B=E^0P_rV$.  If we let $T$ be the
 tautological bundle over $G_rV$, it is easy to see that $P_rT=P_rV$
 and so $B=E^0P_rT=\O_{P_rD(T)}$.  It is also easy to see that
 $D(T)=f^*\Db$, so $P_rD(T)=f^*P_r\Db$, and so
 \[ B=\O_{P_rD(T)}=E^0G_r(V)\ot_A\O_{P_r\Db}=E^0G_r(V)\ot_AB, \]
 as required. 
\end{proof}

We conclude this section with a consistency check that will be useful
later.  
\begin{definition}\label{defn-rho}
 Given a one-dimensional complex vector space $L$ and an arbitrary
 complex vector space $V$, we define $\rho\:PV\xra{}P(L\ot V)$ by
 $\rho(M)=L\ot M$.  This is evidently a homeomorphism.  If $V$ has the
 form $V=\Hom(L,W)=L^*\ot W$ then we identify $L\ot V$ with $W$ in the
 obvious way, and thus obtain a homeomorphism
 $\rho\:P(\Hom(L,W))\xra{}PW$.  All this clearly works equivariantly,
 and fibrewise for vector bundles.
\end{definition}

\begin{proposition}\label{prop-PVW-consistency}
 Let $X$ be a space equipped with two complex vector bundles $V$ and
 $W$.  Let $p\:PV\xra{}X$ be the projection, and let $T$ be the
 tautological bundle over $PV$, so we have a bundle $\Hom(T,p^*W)$
 over $PV$.  Then there is a natural homeomorphism 
 \[ P(V\op W)/PW = PV^{\Hom(T,p^*W)}. \]
\end{proposition}
\begin{proof}
 Put $U=\Hom(T,p^*W))$.  We will construct a diagram as follows:
 \[ \xymatrix{
  {PU} \ar@{ >->}[d]_{i_0} \rto^{\rho}_{\simeq} &
  {P(p^*W)} \ar@{ >->}[d]_{i_1} \ar@{->>}[rr]^{\text{proj}} & & 
  {PV} \ar@{ >->}[d]^{i_2} \\
  {P(\C\op U)} \rto^{\simeq}_{\rho} &
  {P(T\op p^*W)} \ar@{ >->}[r] &
  {P(p^*(V\op W))} \ar@{->>}[r]_{\text{proj}} &
  {P(V\op W)}
  } \]
 First note that the obvious map $\C\xra{}\Hom(T,T)$ is an isomorphism
 (because $T$ has dimension one), so
 \[ \C\op U = \Hom(T,T)\op \Hom(T,p^*W) = \Hom(T,T\op p^*W). \]
 Given this, it is clear that we have homeomorphisms $\rho$ as
 indicated; this gives the left hand half of the diagram, and shows
 that the cofibre of $i_0$ is homeomorphic to that of $i_1$.

 Next, observe that $T$ is a subbundle of $p^*V$ so
 $T\op p^*W\leq p^*(V\op W)$, so $P(T\op p^*W)\sse P(p^*(V\op W))$ as
 indicated.  There is also an obvious projection
 $P(p^*(V\op W))\xra{}P(V\op W)$, giving the right hand rectangle in
 the diagram.  Note also that $P(p^*W)=p^*PV=PV\tm_XPW$.  
 
 We next consider in more detail the map
 $P(T\op p^*W)\xra{}P(V\op W)$, which we shall call $\tau$.  A point
 in $P(T\op p^*W)$ consists of a triple $(x,L,M)$, where $x\in X$ and
 $L$ is a one-dimensional subspace of $V_x$ and $M$ is a
 one-dimensional subspace of $L\op W_x$.  We have
 $\tau(x,L,M)=(x,M)\in P(V\op W)$.  Suppose we start with a point
 $(x,M)\in P(V\op W)$.  If $M\in PW_x$ then it is clear that
 $\tau^{-1}\{(x,M)\}=PV_x\tm\{M\}$.  On the other hand, if
 $M\not\in PW_x$ then the image of $M$ under the projection
 $V_x\op W_x\xra{}V_x$ is a one-dimensional subspace $L\leq V_x$, and
 the point $(x,L,M)$ is the unique preimage of $(x,M)$ under $\tau$.
 This means that the rectangle is a pullback, in which the horizontal
 maps are surjective.  Using this, we see that $\tau$ induces
 a homeomorphism from the cofibre of $i_1$ to that of $i_2$.

 The cofibre of $i_0$ is $PV^U$, and the cofibre of $i_2$ is
 $P(V\op W)/PW$, so these are homeomorphic as claimed.
\end{proof}

As a corollary of the above, we have 
$E^0(P(V\op W),PW)=\tE^0PV^{\Hom(T,p^*W)}$.  We can use the projective
bundle theorem and the Thom isomorphism to calculate both sides in
terms of divisors, and they are not obviously the same.  Nonetheless,
there is an isomorphism between them that can be constructed by pure
algebra, as explained in the following result.

\begin{proposition}\label{prop-PQ-consistency}
 Let $C=\spf(R)$ be a formal multicurve group over $S=\spec(k)$,
 equipped with two divisors $P=\spec(R/K)$ and $Q=\spec(R/L)$.  Define
 an automorphism $\rho$ of $P\tm_SC$ by $\rho(a,b)=(a,b+a)$ (so
 $\rho^{-1}(a,b)=(a,b-a)$) and let $p\:P\xra{}S$ be the projection.
 Then there is a natural isomorphism 
 \[ L(\rho^{-1}(p^*Q)) = L/KL = \ker(\O_{P+Q}\xra{}\O_Q). \]
 Moreover, if we have a coordinate $x$ and use it to define a
 difference function and Thom classes, then the above isomorphism
 sends the Thom class in $L(\rho^{-1}(p^*Q))$ to the element
 $f_Q\in L/KL$.
\end{proposition}
\begin{remark}
 In the last part of the statement, it is important that we are using
 the generator $f_Q$ defined as a norm as in
 Definition~\ref{defn-fD}.  As explained in
 Remark~\ref{rem-ordinary-difference}, in the nonequivariant case,
 this is different from the Chern polynomial which is more usually
 used as a generator.  
\end{remark}
\begin{proof}
 Write $Z$ for the scheme $P\tm 0$ and $\Dl$ for the image of the
 diagonal map $P\xra{}P\tm_SP$.  Both of these can be regarded as
 divisors on the multicurve $P\tm_SC$ over $P$, and it is clear that
 $\rho^{-1}(\Dl)=Z$, and so $\rho^{-1}(\Dl+p^*Q)=Z+\rho^{-1}(p^*Q)$.
 It is also clear that $\Dl\leq p^*P$ and so $\Dl+p^*Q\leq p^*(P+Q)$.
 There are evident projection maps $p^*Q\xra{}Q$ and
 $p^*(P+Q)\xra{}P+Q$.  All this fits together into the following
 diagram. 
 \[ \xymatrix{
  {\rho^{-1}(p^*Q)} \rto^{\rho}_{\simeq} \ar@{ >->}[d]_{i_0} &
  {p^*Q} \ar@{->>}[rr] \ar@{ >->}[d]_{i_1} & & 
  {Q} \ar@{ >->}[d]^{i_2} \\
  {Z+\rho^{-1}(p^*Q)} \rto^{\simeq}_{\rho} &
  {\Dl+p^*Q} \ar@{ >->}[r] &
  {p^*(P+Q)} \ar@{->>}[r] &
  {P+Q}   
  } \]
 The kernel of the ring map $i_0^*$ is (essentially by definition) the
 Thom module $L(\rho^{-1}(p^*Q))$.  As $\rho$ is an isomorphism, it
 induces an isomorphism $\ker(i_1^*)\simeq\ker(i_0^*)$.  It will thus
 be enough to show that the map $\ker(i_2^*)\xra{}\ker(i_1^*)$ is also
 an isomorphism.  To be more explicit, write $A=\O_P=R/K$ and
 $B=\O_Q=R/L$.  Let $I$ be the ideal in $A\ot R=\O_{P\tm_SC}$ defining
 the closed subscheme $\Dl$, so $I$ is generated by the images of
 elements $1\ot a-a\ot 1$ for $a\in R$.  The right hand half of the
 above diagram then gives the following diagram of rings and ideals:
 \[ \xymatrix{
  {\frac{A\ot R}{A\ot L}} &
  {\frac{R}{L}} \ar@{ >->}[l]_{f_0} \\
  {\frac{A\ot R}{I.(A\ot L)}} \ar@{->>}[u]^{i_1^*} & 
  {\frac{R}{KL}} \lto_{f_1}  \ar@{->>}[u]_{i_2^*} \\
  {\frac{A\ot L}{I.(A\ot L)}} \ar@{ >->}[u] &
  {\frac{L}{KL}} \lto^{f_2} \ar@{ >->}[u]
  } \]
 The maps $f_i$ have the form $a\mapsto 1\ot a$, and we must show that
 $f_2$ is an isomorphism.  Now choose a generator $g$ for the ideal
 $L$, giving an isomorphism $L\simeq R$ of $R$-modules.  This gives
 isomorphisms $(A\ot L)/I.(A\ot L)\simeq(A\ot R)/I$ and
 $L/KL\simeq R/K=A$, in terms of which $f_2$ becomes the ring map
 $A\xra{}(A\ot R)/I$ given by $a\mapsto(a\ot 1+I)=(1\ot a+I)$.  This
 corresponds to the projection $\Dl\xra{}B$, which is evidently an
 isomorphism as required.
\end{proof}

\begin{remark}\label{rem-consistency}
 In the context of Proposition~\ref{prop-PVW-consistency}, we can take
 $P=D(V)$ and $Q=D(W)$.  We find that
 $D(\Hom(T,q^*V))=\rho^{-1}(q^*P)$, and the diagram in the proof of
 Proposition~\ref{prop-PQ-consistency} can be identified with that in
 the proof of Proposition~\ref{prop-PVW-consistency}.  It follows that
 the isomorphism $L(\rho^{-1}(q^*P))=\ker(\O_{P+Q}\xra{}\O_P)$
 obtained by applying $E^0(-)$ to
 Proposition~\ref{prop-PVW-consistency} is the same as that given by
 Proposition~\ref{prop-PQ-consistency}. 
\end{remark}

\section{Duality}
\label{sec-duality}

Let $D=\spf(R/I)$ be a divisor of degree $r$ on $C$.  In this section
we will prove that $\Hom_S(\OD,\OS)$ is a free module of rank one over
$\OD$, which means that $\OD$ is a Poincar\'e duality algebra over
$\OS$.  More precisely, we will identify $\Hom_S(\OD,\OS)$ with a
subquotient of the module of meromorphic differential forms on $C$.
In the case where $C$ is embeddable, the duality is given by a kind of
residue.  It is therefore reasonable to define the residue map in the
general case so that this continues to hold.

\subsection{Abstract duality}
\label{subsec-duality-general}

It will be convenient to start by considering a more abstract
situation.  Fix a ground ring $k$, and write $\Hom(M,N)$ for
$\Hom_k(M,N)$ and $M\ot N$ for $M\ot_kN$.  
Let $A$ be a $k$-algebra that is a finitely generated projective
module of rank $r$ over $k$, and write
\begin{align*}
 M^\vee &= \Hom(M,k) = \Hom_k(M,k) \\
 N^*    &= \Hom_A(N,A).
\end{align*}
If $M$ is an $A$-module, then we make $M^\vee$ an $A$-module by the
usual rule $(a\phi)(m)=\phi(am)$.

Now Let $I$ be the kernel of the multiplication map 
$\mu\:A\ot A\xra{}A$, and let $J$ be the annihilator of $I$ in 
$A\ot A$, and put $B=(A\ot A)/J$.  Assume that $I$ and $J$ are both
principal.

Given a $k$-linear map $\phi\:A\xra{}k$, we get an $A$-linear map
$1\ot\phi\:A\ot A\xra{}A$, so we can define
\[ \tph = \tht_0(\phi) = (1\ot\phi)|_J \: J \xra{} A. \]
This construction gives a map $\tht_0\:A^\vee\xra{}J^*$.  

\begin{theorem}\label{thm-duality-general}
 The $A$-modules $A^\vee$ and $J$ are both free of rank one (but without
 canonical generator) and the map $\tht_0\:A^\vee\xra{}J^*$ is an
 $A$-linear isomorphism.
\end{theorem}
The rest of this section constitutes the proof.

\begin{lemma}\label{lem-theta-one}
 The map $\tht_0$ is $A$-linear, and the adjoint map
 $\tht_1\:J\xra{}A^{\vee *}$ is an isomorphism.
\end{lemma}
\begin{proof}
 First suppose that $a\in A$ and $\phi\in A^\vee$ and $u\in J$; we must
 show that $(1\ot a\phi)(u)=a((1\ot\phi)(u))$.  From the definitions
 we have $(1\ot a\phi)(u)=(1\ot\phi)((1\ot a)u)$, and
 $(1\ot\phi)((a\ot 1)u)=a((1\ot\phi)(u))$, so it will suffice to show
 that $(1\ot a)u=(a\ot 1)u$.  This holds because $1\ot a-a\ot 1\in I$
 and $IJ=0$.  We now see that $\tht_0$ is $A$-linear, which allows us
 to define the adjoint map $\tht_1\:J\xra{}A^{\vee *}$ by
 $\tht_1(u)(\phi)=\tht_0(\phi)(u)=(1\ot\phi)(u)$.

 Next, as $A$ is $k$-projective, we have $A\ot A=\Hom(A^\vee,A)$.  If an
 element $u\in A\ot A$ corresponds to a map $\al\:A^\vee\xra{}A$, then
 $(1\ot a)u$ corresponds to the map $x\mapsto a\al(x)$, and
 $(a\ot 1)u$ corresponds to the map $x\mapsto\al(ax)$.  It follows
 that $(1\ot a-a\ot 1)u=0$ iff $\al(ax)=a\al(x)$ for all $x\in A$.
 As $I$ is generated by elements of the form $1\ot a-a\ot 1$, we find
 that 
 \[ A^{\vee *} = \Hom_A(A^\vee,A) = \ann(I,\Hom(A^\vee,A)) =
     \ann(I,A\ot A) = J.
 \]
 One can check that the isomorphism arising from this argument is just
 $\tht_1$.
\end{proof}
\begin{remark}\label{rem-duality-general}
 It follows immediately that if $A^\vee$ has an inverse as an
 $A$-module, then that inverse must be $J$, and $\tht_0$ must be an
 isomorphism.
\end{remark}

We now define $\eta_0,\eta_1\:A\xra{}A\ot A$ by
\begin{align*}
 \eta_0(a) &= a\ot 1 \\
 \eta_1(a) &= 1\ot a.
\end{align*}
We regard $A\ot A$ as an $A$-algebra (and thus $I$, $J$ and $B$ as
$A$-modules) via the map $\eta_0$.

\begin{lemma}
 The $A$-modules $I$ and $B$ are projective, both with rank $r-1$.
 Moreover, $J$ is free of rank one as an $A$-module.
\end{lemma}
\begin{proof}
 As $A$ is projective over $k$ with rank $r$, we see that $A\ot A$ is
 projective over $A$ with rank $r$.  There is a short exact sequence
 $I\xra{}A\ot A\xra{\mu}A$, that is $A$-linearly split by $\eta_0$.
 It follows that $I$ is projective over $A$, with rank $r-1$.  As $I$
 is principal with annihilator $J$, we see that $I\simeq(A\ot A)/J=B$
 as $A$-modules, so $B$ is also projective of rank $r-1$.  It follows
 that the short exact sequence $J\xra{}A\ot A\xra{}B$ is $A$-linearly
 split, and so $J$ is projective of rank one.  It is also a principal
 ideal and thus a cyclic module, so it must in fact be free of rank
 one.
\end{proof}

Next, write $\lm^t$ for the $t$'th exterior power functor, and observe
that $\eta_1$ induces a $k$-linear map
\[ \hat{\eta}_1 \: \lm^{r-1}A \xra{}
   A\ot\lm^{r-1}A=\lm^{r-1}_A(A\ot A),
\] 
and the projection $q\:A\ot A\xra{}B$ induces a map
$\lm^{r-1}(q)\:\lm^{r-1}_A(A\ot A)\xra{}\lm^{r-1}_AB$.  We define
$\psi\:\lm^{r-1}A\xra{}\lm^{r-1}_AB$ to be the composite
of these two maps, so
\[ \psi(a_1\wedge\dotsb\wedge a_{r-1}) = 
    q(1\ot a_1)\wedge_A\dotsb\wedge_A q(1\ot a_{r-1}).
\]

\begin{lemma}
 The map $\psi\:\lm^{r-1}A\xra{}\lm_A^{r-1}B$ is an isomorphism.
\end{lemma}
\begin{proof}
 One can see from the definitions that the image of $\psi$ generates
 $\lm^{r-1}_AB$ as an $A$-module.  We will start by showing that
 the image is itself an $A$-submodule, so $\psi$ must be surjective.

 First, we show that $\lm^{r-1}A$ has a natural structure as an
 $A$-module.  Indeed, there is an evident multiplication
 $\nu\:A\ot\lm^{r-1}A\xra{}\lm^rA$, which induces an isomorphism 
 $\nu^\#\:\lm^{r-1}A\simeq\Hom(A,\lm^rA)$.  The $A$-module structure
 on $A$ gives an $A$-module structure on $\Hom(A,\lm^rA)$ which can
 thus be transported to $\lm^{r-1}A$.  More explicitly, there is a
 unique bilinear operation $*\:A\ot\lm^{r-1}A\xra{}\lm^{r-1}A$ such
 that $a\Smash(b*u)=(ab)\Smash u$ for all $a,b\in A$ and
 $u\in\lm^{r-1}A$.  

 This in turn gives an $A\ot A$-module structure on the group
 $A\ot\lm^{r-1}A=\lm^{r-1}_A(A\ot A)$, by the formula
 $(a\ot b)*(c\ot u)=(ac)\ot(b*u)$.  It follows that
 \[ \hat{\eta}_1(b*u) = \eta_1(b) * \hat{\eta}_1(u). \]

 We next claim that $\lm^{r-1}_AB$ is a quotient $A\ot A$-module of
 $\lm^{r-1}_A(A\ot A)$, and is annihilated by $I$.  Indeed, we
 certainly have an $A$-module structure and $(A\ot A)/I=A$ so it will
 suffice to show that the map
 $\lm^{r-1}(q)\:\lm^{r-1}_A(A\ot A)\xra{}\lm^{r-1}_AB$ annihilates
 $I*\lm^{r-1}_A(A\ot A)$.  To see this, choose an element $e$ that
 generates $J$ (so $Ie=0$).  Using the splittable exact sequence
 $J\xra{}A\ot A\xra{q}B$, we see that there is a commutative diagram
 as follows, in which $\chi$ is an isomorphism.
 \[ \xymatrix{
  {\lm^{r-1}_A(A\ot A)} \ar@{->>}[r]^{e\Smash(-)}
   \ar@{->>}[d]_{\lm^{r-1}(q)} &
  {\lm^r_A(A\ot A)} \\
  {\lm^{r-1}_AB} \urto^{\chi}_{\simeq}
  } \]
 Now for $u\in\lm^{r-1}_A(A\ot A)$ we have
 $e\Smash(I*u)=(Ie)\Smash u=0$, so $\lm^{r-1}(q)(I*u)=0$ as required.

 We can now apply the map $\lm^{r-1}(q)$ to the equation
 $\hat{\eta}_1(b*u)=\eta_1(b)*\hat{\eta}_1(u)$ to see that our map
 $\psi\:\lm^{r-1}A\xra{}\lm^{r-1}_AB$ is $A$-linear.  In particular,
 the image of $\psi$ is an $A$-submodule, and thus $\psi$ is
 surjective as explained previously.

 Next, note that $A$ is a projective $k$-module of rank $r$, so the
 same is true of $\lm^{r-1}A$.  On the other hand, $B$ is a projective
 $A$-module of rank $r-1$, so $\lm^{r-1}_AB$ is a projective
 $A$-module of rank one, and thus also a projective $k$-module of rank
 $r$.  Thus $\psi$ is a surjective map between projective $k$-modules
 of the same finite rank, so it is necessarily an isomorphism.
\end{proof}

\begin{example}\label{eg-psi-embedded}
 It is illuminating to see how this works out in the case where
 $A=k[x]/f(x)$, where $f(x)=\sum_{i=0}^ra_ix^{r-i}$ is a monic
 polynomial of degree $r$.  We write $x_0$ for $x\ot 1$ and $x_1$ for
 $1\ot x$, so $A\ot A=k[x_0,x_1]/(f(x_0),f(x_1))$.  Put 
 \begin{align*}
  d(x_0,x_1) &= x_1-x_0 \\
  e(x_0,x_1) &= (f(x_1)-f(x_0))/(x_1-x_0)
              = \sum_{i+j<r} a_{r-i-j-1}x_0^ix_1^j.
 \end{align*}
 One checks that $I$ is generated by $d$ and $J$ is generated by $e$.
 Put 
 \[ u_i =
    (-1)^i x^0\Smash\ldots\Smash\widehat{x^i}
              \Smash\ldots\Smash x^{r-1},
 \]
 so $\{u_0,\ldots,u_{r-1}\}$ is a basis for $\lm^{r-1}A$ over $k$.  If
 we put $v=x^0\Smash\ldots\Smash x^{r-1}\in\lm^rA$, then
 $x^i\Smash u_j=\dl_{ij}v$, so
 $x_1^i\Smash\hat{\eta}_1(u_j)=\dl_{ij}\hat{\eta}_1(v)$.  Using this, we
 find that 
 \[ \chi\psi(u_j) = e\Smash\hat{\eta}_1(u_j) =
     (\sum_{i=k}^{r-1} a_{i-k}x^{r-1-i})\hat{\eta}_1(v).
 \]
 This means that the matrix of the map
 $\chi\psi\:\lm^{r-1}A\xra{}\lm^r_A(A\ot A)$ (with respect to the
 obvious bases) is triangular, with ones on the diagonal.  The map
 $\chi\psi$ is thus an isomorphism, and the same is true of $\chi$, so
 $\psi$ is an isomorphism as expected.
\end{example}

\begin{proof}[Proof of Theorem~\ref{thm-duality-general}]
 First, suppose we have a map $\phi\:A\xra{}k$.  It is well-known that
 there is a unique derivation $i_\phi$ of the exterior algebra
 $\lm^*A$ whose effect on $\lm^1A$ is just the map
 $\lm^1A=A\xra{\phi}k=\lm^0A$ (this is called \emph{interior
 multiplication} by $\phi$).  We write $\zt(\phi)$ for the map
 $i_\phi\:\lm^rA\xra{}\lm^{r-1}A$.  This construction gives a map
 $\zt\:A^\vee\xra{}\Hom(\lm^rA,\lm^{r-1}A)$.  If we have a basis for
 $A$ then we find that $\zt$ sends the obvious basis of $A^\vee$ to
 the obvious basis for $\Hom(\lm^rA,\lm^{r-1}A)$ (up to sign), so
 $\zt$ is an isomorphism.  In general, $A$ need only be projective
 over $k$ but we can still choose a basis Zariski-locally on
 $\spec(k)$ and the argument goes through.  It follows that $\zt$ is
 always an isomorphism.

 Next, as mentioned above, the short exact sequence 
 $J\xra{}A\ot A\xra{}B$ gives
 $J\ot_A\lm^{r-1}_AB=\lm^r_A(A\ot A)=A\ot\lm^rA$.  As the modules $J$,
 $\lm^{r-1}_AB$ and $A\ot\lm^rA$ are all dualisable, we deduce that
 \begin{align*}
  J^* &= \Hom_A(A\ot\lm^rA,\lm^{r-1}_AB) \\
      &= \Hom(\lm^rA,\lm^{r-1}_AB) \\
      &= \Hom(\lm^rA,\lm^{r-1}A) \\
      &= A^\vee.
 \end{align*}
 In particular, we see that $A^\vee$ is an invertible $A$-module, so
 Remark~\ref{rem-duality-general} tells us that the map $\tht_0$ must
 be an isomorphism.  (In fact, the above chain of identifications
 implicitly constructs an isomorphism $\tht'_0\:A^\vee\xra{}J^*$, and
 one could presumably check directly that $\tht'_0=\tht_0$, but we
 have not done so.)
\end{proof}

\begin{definition}
 The isomorphism $J^*=A^\vee$ gives $J^{*\vee}=A^{\vee\vee}=A$; we
 let $\ep\:J^*\xra{}k$ be the element of $J^{*\vee}$ corresponding to
 $1\in A$ under this isomorphism.  Equivalently, $\ep$ is the unique
 map such that $\ep(\tht_0(\phi))=\ep((1\ot\phi)|_J)=\phi(1)$ for all
 $\phi\in A^\vee$.
\end{definition}

We will prove later that in cases arising from topology, the map $\ep$
can be identified with a Gysin map.  We conclude this section with an
algebraic characterisation of $\ep$ that will be the basis of the
proof.  

\begin{construction}\label{cons-u}
 Given $\lm\in J^*$ and $a\in A$ we can define a map
 $m(a\ot\lm)\:J\xra{}A\ot A$ by $e\mapsto a\ot\lm(e)$.  This map is
 $A$-linear if we use the second copy of $A$ to make $A\ot A$ into an
 $A$-module.  Thus, this construction gives a map
 $m\:A\ot J^*\xra{}\Hom_A(J,A\ot A)$, and as $A$ 
 is projective over $k$, this is easily seen to be an isomorphism.
 Under the inverse of this isomorphism, the inclusion $J\xra{}A\ot A$
 corresponds to an element $u\in A\ot J^*$.  Lemma~\ref{lem-u} will
 give a more concrete description of this element.
\end{construction}

\begin{proposition}\label{prop-epsilon}
 The map $\ep\:J^*\xra{}k$ is such that $(1\ot\ep)(u)=1$ (where $u$ is
 as constructed above).  Moreover, $\ep$ is the unique map with this
 property. 
\end{proposition}
The proof will follow after some discussion.

It is convenient to choose a generator $e=\sum_ia_i\ot b_i$ for $J$,
and a dual generator $\eta$ for $J^*$, so $\eta(e)=1$.  We then put
$\psi=\tht_0^{-1}(\eta)\in A^\vee$; this is the unique map
$\psi\:A\xra{}k$ such that $(1\ot\psi)(e)=1$.  

\begin{lemma}\label{lem-psi}
 We have $\psi(a)=\ep(a\eta)$ for all $a\in A$.
\end{lemma}
\begin{proof}
 Using the $A$-linearity of $\tht_0$ and the fact that
 $\ep\tht_0(\phi)=\phi(1)$, we see that
 \[ \ep(a\eta) = \ep\tht_0(a\psi) = (a\psi)(1) = \psi(a). \]
\end{proof}

\begin{lemma}\label{lem-u}
 The element $u$ in Construction~\ref{cons-u} is given by
 \[ u = e.(1\ot\eta) = \sum_i a_i\ot(b_i\eta). \]
\end{lemma}
\begin{proof}
 The element $v:=e.(1\ot\eta)$ corresponds to the map
 $i=m(v)\:J\xra{}A\ot A$ given by $i(x)=\sum_i a_i\ot (b_i\eta(x))$.  In
 particular, we have $i(e)=\sum_ia_i\ot b_i=e$, so $i$ is the
 inclusion, so $v=u$.
\end{proof}

\begin{proof}[Proof of Proposition~\ref{prop-epsilon}]
 Recall that $\ep$ is the image of $1$ under an isomorphism
 $A\simeq A^{\vee\vee}\simeq J^{*\vee}$, so it certainly generates
 $J^{*\vee}$.  It will thus suffice to check that $(1\ot(t\ep))(u)=t$
 for all $t\in A$.  The calculation is as follows:
 \begin{align*}
  (1\ot(t\ep))(u) &= \sum_i a_i\ep(tb_i\eta) \\
                  &= \sum_i a_i\psi(t b_i) \\
                  &= (1\ot\psi)((1\ot t)e) \\
                  &= (1\ot\psi)((t\ot 1)e) \\
                  &= t.(1\ot\psi)(e) \\
                  &= t.
 \end{align*}
 The first equality is Lemma~\ref{lem-u}, and the second is
 Lemma~\ref{lem-psi}.  The fourth equality holds because
 $(1\ot t)e=(t\ot 1)e$, and the last equality is essentially the
 definition of $\psi$.
\end{proof}

\subsection{Duality for divisors}
\label{subsec-duality-divisors}

Now consider a divisor $D=\spec(A)$ on a multicurve $C=\spf(R)$ over a
scheme $S=\spec(k)$.  In this section, we explain and prove the
following theorem.

\begin{theorem}\label{thm-duality}
 For any divisor $D=\spec(\OC/I_D)$, there is a natural isomorphism 
 \[ \Hom_{\OS}(\OD,\OS) = (I_D^{-1}/\OC)\ot_{\OC}\Om. \]
 (The right hand side consists of meromorphic differential forms whose
 polar divisor is less than or equal to $D$, modulo holomorphic
 differential forms; it is easily seen to be free of rank one over
 $\OD$.)  Moreover, if a map $\phi\:\OD\xra{}\OS$ corresponds to a
 meromorphic form $\mu$, then $\phi(1)=\res(\mu)$.
\end{theorem}

The proof is postponed to the end of the section.  The last part of
the theorem is not yet meaningful, as we have not defined residues.
The definition will be such as to make the claim trivial, but we will
also check that the definition is compatible with the usual one in the
case of embeddable multicurves.

The first step in proving the theorem is to show that the theory in
the previous section is applicable.

\begin{definition}\label{defn-perfect}
 Let $X$ be a scheme over $S$, with closed subschemes $Y$ and $Z$.
 Suppose that $\OX$, $\OY$ and $\OZ$ are all finitely generated
 projective modules over $\OS$.  We say that $Y$ and $Z$ are
 \emph{perfectly complementary} if
 \begin{itemize}
  \item[(a)] the ideals $I_Y$ and $I_Z$ are principal.
  \item[(b)] $\ann(I_Y)=I_Z$ and $\ann(I_Z)=I_Y$.
 \end{itemize}
\end{definition}

\begin{lemma}\label{lem-sum-complementary}
 If $D_0$ and $D_1$ are divisors on a multicurve $C$, then $D_0$ and
 $D_1$ are perfectly complementary in $D_0+D_1$.
\end{lemma}
\begin{proof}
 Put $A_i=\O_{D_i}$; this is a finitely generated projective module
 over $\OS$, and has the form $\OC/f_i$ for some regular element
 $f_i\in\OC$.  We also put $B=\OC/(f_0f_1)=\O_{D_0+D_1}$, so
 $A_i=B/f_i$.   

 Suppose that $gf_1=0\pmod{f_0f_1}$.  Then $(g-hf_0)f_1=0$ for some
 $h\in\OC$, but $f_1$ is regular, so $g=hf_0=0\pmod{f_0}$.  This shows
 that the annihilator of $f_1$ in $B$ is generated by $f_0$, and by
 symmetry, the annihilator of $f_0$ is generated by $f_1$, and this
 proves the lemma.
\end{proof}

\begin{corollary}\label{cor-div-comp}
 Let $D$ be a divisor on a multicurve $C$.  Then the diagonal
 subscheme $\Dl\sse D\tm_SD$ and the subscheme $P_2D\sse D\tm_SD$ are
 perfectly complementary.
\end{corollary}
\begin{proof}
 Let $q\:D\xra{}S$ be the projection.  We can regard $\Dl$ and $P_2D$
 as divisors on the multicurve $q^*C$ over $D$, with
 $D\tm_SD=q^*D=\Dl+P_2D$, so the claim follows from the lemma.
\end{proof}

\begin{corollary}\label{cor-duality}
 Let $D$ be a divisor on a multicurve $C$, and put
 $J=\ker(\O_{D^2}\xra{}\O_{P_2D})$.  Then $J$ is a free $\OD$-module
 of rank one, and there is a natural isomorphism
 \[ \tht_0\:\Hom_{\OD}(J,\OD) \xra{} \Hom_{\OS}(\OD,\OS) \]
 given by $\tht_0(\phi)=(1\ot\phi)|_J$.
\end{corollary}
\begin{proof}
 Let $I$ be the kernel of the multiplication map
 $\O_{D^2}=\OD^{\ot 2}\xra{}\OD$, so that $\O_{\Dl}=\O_{D^2}/I$.  Note
 that $I$ is principal, generated by any difference function on $C$.
 We see from Corollary~\ref{cor-div-comp} that $J$ is the annihilator
 of $I$, and that $J$ is also principal.  We can thus apply
 Theorem~\ref{thm-duality-general} to get the claimed isomorphism.
\end{proof}

To proceed further, we need a better understanding of the ideal $J$.

\begin{definition}
 For the rest of this section, we will use the following notation.
 \begin{align*}
  k &= \OS \\
  R &= \OC \\
  A &= \OD \\
  R_2 &= R\hot R \\
  A_2 &= A \ot A \\
  \tI &= \ker(R_2 \xra{\mu} R) \\
  I &= \ker(A_2 \xra{\mu} A) \\
  J &= \ann_{A_2}(I) \\
  K &= I_D = \ker(R \xra{} A) \\
  K_2 &= K\hot R + R\hot K = \ker(R_2 \xra{} A_2) \\
  \tJ &= \{u \in R_2 \st u\tI\sse K_2\}.
 \end{align*}
 We will also choose a difference function $d$ on $C$ (so $d\in\tI$),
 and let $\al$ denote the image of $d$ in $\Om$.  We choose a
 generator $f$ of $K$, and we let $e$ denote the unique element of
 $R_2$ such that $1\ot f-f\ot 1=de$.  We will check later that the
 image of $e$ in $A_2$ is a generator of $J$.  After that, we will
 write $\eta$ for the dual generator of $J^*$, and $\psi$ for the
 corresponding element of $A^\vee$.
\end{definition}

\begin{remark}\label{rem-tJ}
 It is clear that $\tJ$ is the preimage of $J$ in $R_2$, so
 $J=\tJ/K_2$.  As $\tI$ is generated by a single regular element, one
 can check that $\tI\cap K_2=\tI\tJ$.
\end{remark}

\begin{definition}\label{defn-xi}
 We define $\chi_0\:R\xra{}\tI$ by $\chi_0(a)=1\ot a-a\ot 1$, and
 observe that 
 \[ \chi_0(ab) = (a\ot 1) \chi_0(b) + (1\ot b) \chi_0(a). \]
 Given a difference function $d$ on $C$, we let $\xi(a)\in R_2$ denote
 the unique element such that $\chi_0(a)=\xi(a)d$.  We again have
 \[ \xi(ab) = (a\ot 1) \xi(b) + (1\ot b) \xi(a). \]
\end{definition}
\begin{lemma}\label{lem-nu}
 There is a unique map $\nu\:K_2\xra{}K/K^2$ such that 
 \begin{align*}
  \nu(a\ot b) &= ab \text{ whenever } a\in K \text{ and } b \in R \\
  \nu(a\ot b) &= 0  \text{ whenever } a\in R \text{ and } b \in K.
 \end{align*}
 Moreover, we have $\nu(K_2\tI)=0$.
\end{lemma}
\begin{proof}
 Using a $k$-linear splitting of the sequence $K\xra{}R\xra{}A$, we
 see that $(K\hot R)\cap(R\hot K)=K\hot K$, so we have
 \[ \frac{K_2}{R\hot K} = \frac{K\hot R}{K\hot K}. \]
 The multiplication map $\mu\:R\hot R\xra{}R$ evidently induces a map
 $(K\hot R)/(K\hot K)\xra{}K/K^2$.  Putting this together, we get a
 map
 \[ \nu = \left(
     K_2 \xra{} \frac{K_2}{R\hot K} = \frac{K\hot R}{K\hot K}
         \xra{\mu} \frac{K}{K^2}
    \right).
 \]
 It is clear that this is uniquely characterised by the stated
 properties.  As $\nu$ is essentially given by $\mu$ on $K\hot R$ and
 $\mu(\tI)=0$, we see that $\nu(\tI.(K\hot R))=0$.  We also have
 $\nu(R\hot K)=0$ and so $\nu(\tI.(R\hot K))=0$, so $\nu(\tI K_2)=0$
 as claimed. 
\end{proof}

\begin{proposition}\label{prop-xi}
 The map $\xi$ induces a $A$-linear isomorphism $K/K^2\xra{}J$, so $J$
 is freely generated over $A$ by the element $e=\xi(f)$.  (Note
 however that this isomorphism is not canonical, because it depends on
 the choice of $d$.)
\end{proposition}
\begin{proof}
 Suppose that $a\in K$, so $\chi_0(a)\in K_2$.  As $\tI=R_2d$ we see
 that $\xi(a)\tI=R_2\chi_0(a)\sse K_2$, so $\xi(K)\sse\tJ$, so we get an
 induced map $K\xra{}\tJ/K_2=J$.  Using the product formula for $\xi$,
 we deduce that this map is $R$-linear and induces a map
 $K/K^2\xra{}J$.  In the opposite direction, we define
 $\zt\:\tJ\xra{}K/K^2$ by $\zt(u)=\nu(ud)$.  If $u\in K_2$ then
 $ud\in \tI K_2$ so $\nu(ud)=0$.  Thus, $\zt$ induces a map
 $J=\tJ/K_2\xra{}K/K^2$.  It is easy to see that
 $\zt\xi=1\:K/K^2\xra{}K/K^2$, and both $J$ and $K/K^2$ are invertible
 $A$-modules, so $\zt$ and $\xi$ must be mutually inverse
 isomorphisms.  It follows immediately that $J$ is freely generated by
 $e$.
\end{proof}

Our next task is to reformulate the above isomorphism in a way that is
independent of any choices.  
\begin{proposition}\label{prop-chi}
 Put 
 \[ \bOm=\Om|_D=\Om\ot_RA=\tI\ot_{R_2}A. \]
 There is a \emph{natural} isomorphism
 $\chi\:K/K^2\xra{}J\ot_R\Om=J\ot_A\bOm$ of $A$-modules, satisfying
 $\chi(a)=\xi(a)\ot\al$, where $\al$ is the image of $d$ in $\Om$.  By
 adjunction, there is also a natural isomorphism
 $\chi'\:J^*\xra{}(K/K^2)^*\ot_A\bOm=(K/K^2)^*\ot_R\Om$.
\end{proposition}
\begin{proof}
 We have 
 \[ J\ot_A\bOm = J\ot_A(A\ot_{R_2}\tI) = (\tJ/K_2)\ot_{R_2}\tI =
    (\tI\tJ)/(\tI K_2) = (\tI\cap K_2)/(\tI K_2).
 \]
 We have seen that the map $\chi$ sends $K$ to $\tI\cap K_2$ and $K^2$
 to $\tI K_2$, so it induces a map $\chi\:K/K^2\xra{}J\ot_A\bOm$,
 which is obviously canonical.  One checks from the definitions that
 $\chi(a)=\xi(a)\ot\al$, where $\xi$ and $\al$ are defined in terms of
 a difference function $d$ as above.  As $\xi$ is an isomorphism and
 $\bOm$ is freely generated by $\al$ over $A$, we conclude that $\chi$
 is an isomorphism.
\end{proof}

Our next task is to interpret the module $(K/K^2)^*$.
\begin{definition}\label{defn-MC}
 Let $C$ be a formal multicurve over $S$.  We say that an element
 $f\in\OC$ is \emph{divisorial} if it is not a zero-divisor, and
 $\OC/f$ is a projective $\OS$-module of finite rank.  One can check
 that the set of divisorial elements is closed under multiplication,
 so we can invert it to get a new ring $\MC$, whose elements we call
 \emph{meromorphic functions}.  We say that a meromorphic function is
 \emph{divisorial} if it can be written as $f/g$, where $f$ and $g$
 are divisorial elements of $\OC$ (this can be seen to be compatible
 with the previous definition).  A \emph{fractional ideal} is an
 $\OC$-submodule $I\leq\MC$ that can be generated by a divisorial
 meromorphic function.  The set of fractional ideals forms a group
 under multiplication, with $I^{-1}=\{f\in\MC\st fI\sse\OC\}$.
\end{definition}

\begin{lemma}\label{lem-fract-id}
 There is a natural isomorphism $(K/K^2)^*=(K^{-1}/R)$.  With respect
 to this, the generator $f$ of $K/K^2$ is dual to the generator $1/f$
 of $K^{-1}/R$.
\end{lemma}
\begin{proof}
 The multiplication map $K\hot_RK^{-1}\xra{}R$ induces a map
 $(K/K^2)\ot_A(K^{-1}/R)\xra{}A$, and thus a map
 $K^{-1}/R\xra{}(K/K^2)^*$.  This is easily seen to be an
 isomorphism.  The statement about generators is clear.
\end{proof}

\begin{proof}[Proof of Theorem~\ref{thm-duality}]
 Everything except the last part now follows immediately from
 Corollary~\ref{cor-duality}, Proposition~\ref{prop-chi} and
 Lemma~\ref{lem-fract-id}.  The residue map will be defined in
 Definition~\ref{defn-res-multi}, and then the last part of the
 theorem will be true by definition.
\end{proof}

We can make the theorem more explicit as follows.
\begin{proposition}\label{prop-duality-explicit}
 The natural isomorphism 
 \[ A^\vee \xra{\tht_0} J^* \xra{\chi'} K^{-1}/R\ot_A\Om \]
 sends $\phi$ to the element $((1\ot\phi)(e)/f)\ot\al$.
\end{proposition}
\begin{proof}
 Using Proposition~\ref{prop-xi}, we see that $e$ generates $J$, so
 there is a unique element $\eta\in J^*$ with $\eta(e)=1$.  The
 natural isomorphism $K/K^2\xra{}J\ot_A\Om$ sends $f$ to $e\ot\al$, so
 the adjoint map $\chi'\:J^*\xra{}K^{-1}/R\ot\Om$ sends $\eta$ to
 $(1/f)\ot\al$, and thus sends $a\eta$ to $(a/f)\ot\al$.  Next, we
 certainly have $\tht_0(\phi)=a\eta$ for some $a\in A$, and by
 evaluating this equation on $e$ we find that $a=(1\ot\phi)(e)$.
 It follows that
 \[ \chi'\tht_0(\phi) = (a/f)\ot\al = ((1\ot\phi)(e)/f)\ot\al \]
 as claimed.
\end{proof}

We next examine how this works in the case where $C$ is embeddable,
say $C=\spf(k[x]^\wedge_{(g)})$ for some monic polynomial $g$.  We
then have $K=Rf$ for some monic polynomial $f$ that divides some power
of $g$, and $A=R/K=k[x]/f$.  

\begin{definition}\label{defn-res-embedded}
 Suppose we have a ring $k$ and an expression
 $\al=f(x)\bd x=p(x)\bd x/q(x)$, where $p$ and $q$ are polynomials
 with $q$ monic; we then define the residue $\res(\al)$ as follows.
 Let $R'$ denote the ring of series of the form
 $u(x)=\sum_{n=-\infty}^Na_nx^n$ for some finite $N$.  Clearly
 $k[x]\sse R'$.  Moreover, if $q(x)$ is a monic polynomial then we can
 write $q(x)=x^nr(1/x)$ for some polynomial $r(t)$ with $r(0)=1$.  It
 follows that $r(1/x)$ is invertible in $R'$, and thus the same is
 true of $q(x)$, so we can expand out $f(x)=p(x)/q(x)$ as an element
 of $R'$, say $f(x)=\sum_{n=-\infty}^Na_nx^n$.  We then put
 \[ \res(\al) = a_{-1}. \]
\end{definition}

\begin{remark}\label{rem-res-props}
 In the case $k=\C$, one can check that $\res(\al)$ is the sum of the
 residues of $\al$ at all its poles, so the integral of $\al$ around
 any sufficiently large circle is $2\pi i\res(\al)$.  By standard
 arguments, most formulae that hold when $k=\C$ will be valid for all
 $k$.  In particular, we have
 \begin{itemize}
  \item $\res(f(x)\bd x)=0$ if $f$ is a polynomial
  \item $\res(f'(x)\bd x)=0$ for any $f=p/q$ as above
  \item $\res(f'(x)\bd x/f(x))=\deg(p)-\deg(q)$ if $f=p/q$ for some
   monic polynomials $p$ and $q$.
 \end{itemize}
\end{remark}
\begin{lemma}\label{lem-res-calc}
 Suppose that $f(x)=p(x)/q(x)$ where $q$ is monic of degree $n$, and
 \[ p(x)=\sum_{i=0}^{n-1}b_ix^i\pmod{q(x)}. \]
 Then $\res(f(x)\bd x)=b_{n-1}$.
\end{lemma}
\begin{proof}
 First, put $m(x)=\sum_{i=0}^{n-1}b_ix^i$, so $p(x)=m(x)+l(x)q(x)$ for
 some polynomial $l(x)$, so $f(x)=m(x)/q(x)+l(x)$.  We have
 $\res(l(x)\bd x)=0$, so it will suffice to show that
 $\res(m(x)/q(x)\bd x)=b_{n-1}$.

 Next, write $q(x)=x^nr(1/x)$ as in the definition, and put
 $u(x)=1/r(1/x)$, so $u(x)=\sum_{i=0}^\infty a_ix^{-i}$ for some
 coefficients $a_i\in k$ with $a_0=1$.  We have
 \[ x^j/q(x)=x^{j-n}u(x)=\sum_{i\geq 0}a_ix^{j-n-i}, \]
 so 
 \[ \res(x^j\bd x/q(x)) = 
     \begin{cases}
      0 & \text{ if } 0 \leq j < n-1 \\
      1 & \text{ if } j = n-1.
     \end{cases}
 \]
 The claim follows immediately.
\end{proof}

\begin{proposition}\label{prop-res-embedded}
 If $A=k[x]/f(x)$, then the map 
 \[ (K^{-1}/R)\ot\Om \xra{(\chi')^{-1}} J^* \xra{\ep} k \]
 is just the residue.
\end{proposition}
\begin{proof}
 As in example~\ref{eg-psi-embedded}, we put
 \begin{align*}
  f(x)       &= \sum_{i=0}^ra_ix^{r-i} \\
  d(x_0,x_1) &= x_1-x_0 \\
  \al        &= \bd x\in \Om \\
  e(x_0,x_1) &= (f(x_1)-f(x_0))/(x_1-x_0)
              = \sum_{i+j<r} a_{r-i-j-1}x_0^ix_1^j.
 \end{align*}
 This is of course compatible with the notation in
 Proposition~\ref{prop-duality-explicit}, so
 $\chi'\tht_0(\phi)=(1\ot\phi)(e)\bd x/f(x)$.  
 The map $\ep$ is characterised by the fact that
 $\ep(\tht_0(\phi))=\phi(1)$, so it will suffice to check that 
 $\res((1\ot\phi)(e)\bd x/f(x))=\phi(1)$ for all $\phi\in A^\vee$.

 Now let $\{\zt_0,\ldots,\zt_{r-1}\}$ be the basis for $A^\vee$ dual to
 the basis $\{x^0,\ldots,x^{r-1}\}$ for $A$.  We then have
 $(1\ot\zt_j)(x_1^k)=\dl_{jk}$, and it follows that
 \[ (1\ot\zt_j)(e) = \sum_{i=0}^{r-1-j} a_{r-i-j-1}x^i. \]
 Using Lemma~\ref{lem-res-calc} we deduce that
 $\res((1\ot\zt_j)(e)\bd x/f(x))=0$ for $j>0$, whereas for $j=0$ we
 get $a_0$, which is $1$ because the polynomial
 $f(x)=\sum_{i+j=r}a_jx^j$ is monic.  On the other hand, we also have
 $\zt_j(1)=\dl_{0j}$ by definition, so
 $\res((1\ot\zt_j)(e)\bd x/f(x))=\zt_j(1)$ as required.
\end{proof}

Given this, it would be reasonable to define residues on multicurves
using the maps $\ep$.  To make this work properly, we need to check
that these maps are compatible for different divisors.

\begin{proposition}\label{prop-res-compatible}
 Suppose we have divisors $D_0\sse D_1$, corresponding to ideals
 $K_1\leq K_0\leq R$.  Let $j\:K_0^{-1}/R\xra{}K_1^{-1}/R$ be the
 evident inclusion, and let $q\:A_1=R/K_1\xra{}R/K_0=A_0$ be the
 projection.  Define $\dl_i\:A_i^\vee\xra{}k$ by
 $\dl_i(\phi)=\phi(1)$.  Then the following diagram commutes:
 \[ \xymatrix{
  {k} \ar@{=}[d] &
  {A_0^\vee}
   \lto_{\dl_0} \ar@{ >->}[d]^{q^\vee}
   \rto^(0.3){\chi'\tht_0}_(0.3){\simeq} & 
  {K_0^{-1}/R \ot\Om} \ar@{ >->}[d]^j \\
  k & 
  A_1^\vee \lto^{\dl_1} \rto^(0.3){\simeq}_(0.3){\chi'\tht_0} &
  {K_1^{-1}/R \ot\Om}
  } \]
\end{proposition}
\begin{proof}
 As $q(1)=1$, it is clear that the left hand square commutes.  For the
 right hand square, choose generators $f_i$ for $K_i$ and put
 $e_i=\xi(f_i)$, so that
 \[ \chi'\tht_0(\phi)=((1\ot\phi)(e_i)/f_i)\ot\al \]
 for $\phi\in A_i^\vee$.  

 As $D_0\sse D_1$, we have $f_1=gf_0$ for some $g$, and so 
 $\xi(f_1)=(g\ot 1)\xi(f_0)+(1\ot f_0)\xi(g)$, or in other words
 $e_1=(g\ot 1)e_0+(1\ot f_0)\xi(g)$.

 Now suppose we have $\phi\in A_0^\vee$, so $(q^*\phi)(f_0)=0$, so 
 $(1\ot q^*\phi)((1\ot f_0)\xi(g))=0$.  We then
 have
 \begin{align*}
  \chi'\tht_0q^*(\phi) &=
      ((1\ot q^*\phi)(e_1)/f_1) \ot\al \\
   &= (g (1\ot\phi)(e_0))/(g f_0) \ot \al \\
   &= (1\ot\phi)(e_0)/f_0 \ot\al \\
   &= j\chi'\tht_0(\phi)
 \end{align*}
 as claimed.
\end{proof}

\begin{definition}\label{defn-res-multi}
 Define $\dl\:\Hom_{\OS}(\OD,\OS)\xra{}\OS$ by $\dl(\phi)=\phi(1)$.
 We let $\res\:\MC\ot_{\OC}\Om\xra{}\OS$ be the unique map whose
 restriction to $I_D^{-1}\ot_{\OC}\Om$ is the composite
 \[ I_D^{-1}\ot_{\OC}\Om \xra{}
    (I_D^{-1}/\OC)\ot_{\OC}\Om \simeq
     \Hom_{\OS}(\OD,\OS) \xra{\dl} \OS.
 \]
 (This is well-defined, by Proposition~\ref{prop-res-compatible}, and
 compatible with the classical definition, by
 Proposition~\ref{prop-res-embedded}.)
\end{definition}

\begin{proposition}\label{prop-res-props}
 For any $g,f\in\OC$, if $f$ is divisorial than 
 \[ \res((g/f)\bd f)=\trace_{(\OC/f)/\OS}(g). \]
 In particular, we have $\res((1/f)\bd f)=\dim_{\OS}(\OC/f)$.
 Moreover, we also have
 \[ \res(\bd(g/f)) =
    \res\left(\frac{f \bd(g) - g \bd(f)}{f^2}\right) = 0.
 \]
\end{proposition}
\begin{proof}
 Both facts are well-known for residues in the classical sense, so
 they hold whenever $C$ is embeddable.  Using
 Corollary~\ref{cor-locally-embeddable}, we deduce that they hold for
 a general multicurve $C$.  We will give a more direct and
 illuminating proof for the first fact; we have not been able to
 find one for the second fact.

 We use abbreviated notation as before, with $K=Rf$ so that $A=R/f$.
 The multiplication map $\mu\:A_2\xra{}A$ restricts to give an
 $A$-linear map $\mu\:J\xra{}A$, or in other words an element of
 $J^*$.  The trace map $\tau\:A\xra{}k$ can be regarded as an element
 of $A^\vee$.  We claim that the elements $\tau$, $\mu$ and $(\bd f)/f$
 correspond to each other under our standard isomorphisms
 \[ A^\vee\simeq J^*\simeq (K^{-1}/R)\ot\Om. \]
 To see this, note that $(1\ot\tau)(u)=\trace_{A_2/A}(u)$ for all
 $u\in A_2$.  Using the splittable short exact sequence
 \[ I \xra{} A_2 \xra{\mu} A, \]
 we see that
 \[ (1\ot\tau)(u) = \trace(I \xra{\tm u} I) + \mu(u). \]
 If $u\in J$ then multiplication by $u$ is zero on $I$ and we deduce
 that $(1\ot\tau)(u)=\mu(u)$.  This shows that $\tht_0(\tau)=\mu$ as
 claimed.  
 
 Next, let $e=\xi(f)$ be the standard generator of $J$, and let $\eta$
 be the dual generator of $J^*$, so $\eta(e)=1$.  Using
 Proposition~\ref{prop-duality-explicit}, we see that $\mu$
 corresponds to the element $(\mu(e)/f)\ot\al=(1/f)\ot(\mu(e)\al)$ in
 $(K^{-1}/R)\ot\Om$.  Now, the module $\Om=\tI/\tI^2$ is originally a
 module over $R_2$ that happens to be annihilated by $\ker(\mu)=\tI$,
 and so is regarded as a module over $R$ via $\mu$.  Thus, $\mu(e)\al$
 is just the same as $e\al$.  Moreover, $\al$ is just the image of $d$
 in $\Om$, so $e\al$ is the image of $ed=\xi(f)d=1\ot f-f\ot 1$, and
 this image is by definition just $\bd f$.  Thus, $\mu\in J^*$
 corresponds to $(1/f)\ot\bd f$ as claimed.

 As the isomorphism $A^*\xra{}(K^{-1}/R)\ot\Om$ is $A$-linear, we see
 that $g\tau$ maps to $(g/f)\bd f$, so 
 \[ \res((g/f) \bd f) = (g\tau)(1) = \tau(g), \]
 as claimed.
\end{proof}

\begin{remark}\label{rem-kld}
 It is useful and interesting to reconcile this result
 with~\cite{st:kld}*{Proposition 9.2}.  There we have a $p$-divisible
 formal group $\hC=\spf(R)$ of height $n$ over a formal scheme
 $S=\spf(k)$, where $k$ is a complete local Noetherian ring of residue
 characteristic $p$, and we will assume for simplicity that $k$ is
 torsion-free.  Fix $m\geq 1$ and let $\psi\:\hC\xra{}\hC$ be $p^m$
 times the identity map.  In this context the subgroup scheme
 $D:=\hC[p^m]=\ker(\psi)$ is a divisor of degree $p^{nm}$, so the ring
 $\OD$ is self-dual (with a twist) as before.  Given any coordinate
 $x$, we note that $\OD=R/\psi^*x$, so the meromorphic form
 $\al=\bD(x)/(\psi^*x)$ is a generator of the twisting module
 $(K^{-1}/R)\ot\Om$.  We claim that $\al$ is actually independent of
 $x$.  Indeed, any other coordinate $x'$ has the form $x'=(x+x^2q)u$
 for some $u\in k^\tm$ and $q\in R$.  It follows that
 $\bd_0(x')=u\bd_0(x)$, so that $\bD(x')=u\bD(x)$.  We also have
 $\psi^*(x')=u\psi^*(x)\pmod{\psi^*(x)^2}$, and it follows that
 $\psi^*(x')^{-1}=u^{-1}\psi^*(x)^{-1}\pmod{R}$, so
 $\bD(x')/\psi^*(x')=\bD(x)/\psi^*(x)$ mod holomorphic forms, as
 claimed.  Thus, we have a canonical generator for $(K^{-1}/R)\ot\Om$
 and thus a canonical generator for $A^\vee$, giving a Frobenius form
 on $\OD$.  The cited proposition says that this is the same as the
 Frobenius form coming from a transfer construction.  As discussed in
 the preamble to that proposition, $p^m$ times the transfer form is
 the same as the trace form.  In view of
 Proposition~\ref{prop-res-props}, this means that
 $p^m\al=\bd(\psi^*x)/(\psi^*x)$.  In fact, this is easy to see
 directly.  We know that $\bD(x)$ generates $\Om$ and agrees with
 $\bd(x)$ at zero, so $\bd(x)=(1+xr)\bD(x)$ for some function $r\in
 R$.  It follows that
 \[ \bd(\psi^*x) = \psi^*(\bd(x)) =
     (1+\psi^*(x)\psi^*(r))\psi^*(\bD(x)).
 \]
 As $\bD(x)$ is invariant we have $\psi^*\bD(x)=p^m\bD(x)$.  It
 follows that $\bd(\psi^*x)/(\psi^*x)=p^m\bD(x)/(\psi^*x)=p^m\al$ in
 $(K^{-1}/R)\ot\Om$, as claimed.
\end{remark}

\begin{remark}\label{rem-tate}
 It should be possible to connect our treatment of residues with that
 of Tate~\cite{ta:rdc}.  However, Tate assumes that the ground ring
 $k$ is a field, and it seems technically awkward to remove this
 hypothesis. 
\end{remark}

\subsection{Topological duality}
\label{subsec-duality-top}

Consider a periodically orientable theory $E$, an $A$-space $X$, and
an equivariant complex bundle $V$ over $X$.  To avoid some minor
technicalities, we will assume that $X$ is a finite $A$-CW complex;
everything can be generalised to the infinite case by passage to
(co)limits.  Let $C$ be the multicurve $\spf(E^0(P\CU\tm X))$ over
$S:=\spec(E^0X)$.  We then have a divisor $D=D(V)$ on $C$, to which we
can apply all the machinery in the previous section.  In particular,
this gives us a residue map 
\[ \res\:(I_D^{-1}/\OC)\ot_{\OC}\Om\xra{} \OS. \] 
On the other hand, if we let $\tau$ denote the tangent bundle along
the fibres of $PV$, then there is a stable Pontrjagin-Thom collapse
map $X_+\xra{}PV^{-\tau}$, giving a Gysin map 
\[ p_! \: \tE^0PV^{-\tau} \xra{} E^0X = \OS. \]

\begin{theorem}\label{thm-duality-top}
 There is a natural isomorphism
 $\tE^0PV^{-\tau}=(I_D^{-1}/\OC)\ot_{\OC}\Om$, which identifies the
 Gysin map with the residue map.
\end{theorem}

This is an equivariant generalisation of a result stated by Quillen
in~\cite{qu:fgl}.  Even in the nonequivariant case, we believe that
there is no published proof.  The rest of this section constitutes the
proof of our generalisation.  (The case of nonequivariant
\emph{ordinary} cohomology is easy, and is a special case of the
result proved in~\cite{da:ghf}.)

We retain our previous notation for rings, and write $P^2V=PV\tm_XPV$,
so
\begin{align*}
 k &= \OS = E^0X \\
 R &= \OC = E^0(P\CU\tm X) \\
 A &= \OD = E^0(PV) \\
 R_2 &= R\hot R = E^0(P\CU\tm P\CU\tm X) \\
 A_2 &= A\ot A = E^0(P^2V).
\end{align*}
Next, observe that $P_2V$ is a subspace of $P^2V$, and
Proposition~\ref{prop-Pr} tells us that $E^0P_2V=\O_{P_2D}=A_2/J$, so
$J=E^0(P^2V,P_2V)$.  On the other hand, there is another natural
description of $E^0(P^2V,P_2V)$, which we now discuss.  Let $T$ be the
tautological bundle on $PV$, consider the vector bundles
$T^\perp=V\ominus T$ and $U=\Hom(T,T^\perp)$, and let $\oB U$ denote
the open unit ball bundle in $U$.  A point in $\oB U$ is a triple
$(x,L,\al)$ where $x\in X$ and $L\in PV_x$ and 
$\al\:L\xra{}V_x\ominus L$, such that $\|\al(u)\|<\|u\|$ for all 
$u\in L\sm\{0\}$.  We can thus consider $\graph(\al)$ and
$\graph(-\al)$ as one-dimensional subspaces of 
$L\tm(V_x\ominus L)=V_x$, or in other words points of $PV_x$, so we
have a map $\dl'\:\oB U\xra{}P^2V$ given by 
\[ \dl'(x,L,\al) = (\graph(\al),\graph(-\al)). \]
\begin{proposition}
 The map $\dl'$ is a diffeomorphism $\oB U\xra{}P^2V\sm P_2V$.
\end{proposition}
\begin{proof}
 First, we must show that $\dl'(x,L,\al)\not\in P_2V$, or in other
 words that $\graph(\al)$ is not perpendicular to $\graph(-\al)$.  For
 this, we choose a nonzero element $u\in L$, so
 $v_0=u+\al(u)\in\graph(\al)$ and $v_1=u-\al(u)\in\graph(-\al)$.  It
 follows that $\ip{v_0,v_1}=\|u\|^2-\|\al(u)\|^2$; this is strictly
 positive because $\|\al\|<1$, so the lines are not orthogonal, as
 required.  We therefore have a map $\dl'\:\oB U\xra{}P^2V\sm P_2V$.

 Any element of $P^2V\sm P_2V$ has the form $(x,M_0,M_1)$ where
 $x\in X$ and $M_0,M_1\in PV_x$ and $M_0$ is not orthogonal to $M_1$.
 This means that we can choose elements $u_i\in M_i$ with $\|u_i\|=1$
 and such that $t:=\ip{u_0,u_1}$ is a positive real number.  One
 checks that the pair $(u_0,u_1)$ is unique up to the diagonal action
 of $S^1$.  Put $v=u_0+u_1$ and $w=u_0-u_1$.  By Cauchy-Schwartz we
 have $t\leq 1$, and by direct expansion we have
 \begin{align*}
  \ip{v,w} &= 0 \\
  \ip{v,v} &= 2(1+t) > 0 \\
  \ip{w,w} &= 2(1-t) < \ip{v,v}.
 \end{align*}
 We can thus put $L=\C v\in PV_x$ and define $\al\:L\xra{}L^\perp$ by
 $\al(zv)=zw$; these are clearly independent of the choice of pair
 $(u_0,u_1)$.  As $\|w\|<\|v\|$ we have $\|\al\|<1$.  As
 $v+\al(v)=2u_0$ we have $\graph(\al)=M_0$, and as $v-\al(v)=2u_1$ we
 have $\graph(-\al)=M_1$.  It follows that the construction
 $(x,M_0,M_1)\mapsto(x,L,\al)$ gives a well-defined map
 $\zt\:P^2V\sm P_2V\xra{}\oB U$, with $\dl'\zt=1$.  One can check
 directly that $\zt\dl'$ is also the identity, so $\dl'$ is a
 diffeomorphism as claimed.
\end{proof}

\begin{corollary}
 The bundle $U$ is the normal bundle to the diagonal embedding
 $\dl\:PV\xra{}P^2V$, there is a homeomorphism $P^2V/P_2V=PV^U$, and
 the quotient map $P^2V\xra{}P^2V/P_2V$ can be thought of as the
 Pontrjagin-Thom collapse associated to $\dl$. \qed
\end{corollary}
\begin{remark}
 There are easier proofs of this corollary if one is willing to be
 less symmetrical.
\end{remark}

Now, the above corollary together with Proposition~\ref{prop-chi} and
Section~\ref{sec-differential} gives
\[ \tE^0PV^U = E^0(P^2V,P_2V) = J = K/K^2 \ot_A \bOm^* = 
    K/K^2 \ot_k \om^\vee. 
\]
On the other hand, we have
\[ U = \Hom(T,T^\perp) = \Hom(T,p^*V) \ominus \C, \]
so (using the case $W=V$ of Proposition~\ref{prop-PVW-consistency})
\[ PV^U = \Sg^{-2} PV^{\Hom(T,p^*V)} = \Sg^{-2} P(V\op V)/PV. \]
Remark~\ref{rem-om-top} tells us that $\tE^0(S^{-2})=\om^\vee$, and it
is clear that $E^0(P(V\op V),PV)=K/K^2$.  We thus obtain
\[ \tE^0PV^U = K/K^2 \ot_k \om^\vee \]
again.  These two arguments apparently give two different isomorphisms
$\tE^0PV^U\xra{}K/K^2\ot_k\om^\vee$, but one can show (using
Remark~\ref{rem-consistency}) that they are actually the same.
% {\bf More details}.

We next recall some ideas about Gysin maps.  We discuss the situation
for manifolds, and leave it to the reader to check that everything
works fibrewise for bundles of manifolds, at least in sufficient
generality for the arguments below.  Let $f\:M\xra{}N$ be an analytic
map between compact complex manifolds.  (It is possible to work with
much less rigid data, but we shall not need to do so.)  Let $\tau_M$
and $\tau_N$ be the tangent bundles of $M$ and $N$, and let $\nu_f$ be
the virtual bundle $f^*\tau_N-\tau_M$ over $M$.  Then for any virtual
bundle $U$ over $N$, a well-known variant of the Pontrjagin-Thom
construction gives a stable map $T(f,U)\:N^U\xra{}M^{f^*U+\nu_f}$, and
thus a map $f_!=T(f,U)^*\:\tE^0M^{f^*U+\nu_f}\xra{}\tE^0N^U$.  Using
the ring map $f^*\:E^0N\xra{}E^0M$ we regard the source and target of
$T(f,U)^*$ as $E^0N$-modules, and we find that $T(f,U)^*$ is
$E^0N$-linear.  We also find that $T(f,U)^*$ can be obtained from
$T(f,0)^*$ by tensoring over $E^0N$ with $\tE^0N^U$.  Finally, we have
a composition formula: given maps $M\xra{f}N\xra{g}P$, we have
$\nu_{gf}=\nu_f+f^*\nu_g$ and 
\[ T(f,\nu_g) \circ T(g,0) = T(gf,0) \:
     P^V \xra{} M^{(gf)^*V+\nu_{gf}}.
\]

Now consider the maps $M\xra{\dl}M^2\xra{1\tm\pi}M$, where $\pi$ is
the constant map from $M$ to a point.  We have $\nu_\dl=\tau_M$ and
$\nu_\pi=-\tau_M$, so the transitivity formula says that the composite 
\[ M_+ \xra{1\Smash T(\pi,0)}
   M_+\Smash M^{-\tau} \xra{T(\dl,\nu_{1\tm\pi})} M_+ 
\]
is the identity.  Assuming a K\"unneth isomorphism, we get maps
\[ E^0M \xra{\dl_!} E^0M\ot\tE^0M^{-\tau} \xra{1\ot\pi_!} E^0M, \]
whose composite is again the identity.

Now specialise to the case $M=PV$.  As before we put $A=E^0PV$ and
identify $\tE^0M^\tau$ with $J$, and the map 
$\dl_!=T(\dl,0)^*\:\tE^0M^\tau\xra{}E^0(M^2)$ with the inclusion 
$J\xra{}A\ot A$.  We know that the map 
\[ \dl_! = T(\dl,\nu_{1\tm\pi})^*\:
    A = E^0M \xra{} E^0M\ot\tE^0M^{-\tau} = A\ot J^*
\]
is obtained from $T(\dl,0)^*$ by tensoring over $A\ot A$ with 
$A\ot J^*$.  It follows easily that $\dl_!(1)=u\in A\ot J^*$, where
$u$ is as in Construction~\ref{cons-u}.  The equation 
$(1\ot\pi_!)\dl_!=1$ now tells us that $(1\ot\pi_!)(u)=1$.
Proposition~\ref{prop-epsilon} now tells us that
$\pi_!=\ep\:J^*\xra{}k$.  This proves Theorem~\ref{thm-duality-top}.

\section{Further theory of infinite Grassmannians}
\label{sec-grassmann-more}

Recall from Section~\ref{sec-grassmann} that $G\CU$ denotes the space
of finite-dimensional subspaces of $\CU$, which is the natural
equivariant generalization of the space $G\C^\infty=\coprod_{d\geq
  0}BU(d)$.  We know from Theorem~\ref{thm-grassmann} that $E_0G\CU$
is the symmetric algebra over $E_0$ generated by $E_0P\CU=\O^\vee_C$.
It follows that
\begin{align*}
 \spec(E_0G\CU)  &= \Map(C,\aff^1) \\
 \spf(E^0G_d\CU) &= \Div_d^+(C).
 \spf(E^0G\CU)   &= \Div^+(C).
\end{align*}
In this section, we obtain similar results for spaces analogous to
$\Z\tm BU$, $BU$ and $BSU$.

\begin{definition}
 For any finite-dimensional $A$-universe $U$, we put $\two U=U\op U$.
 We write $U_+$ for $U\op 0$ and $U_-$ for $0\op U$ so
 $\two U=U_++U_-$.  We put
 \[ \tG(U) = G(\two U) = \coprod_{d=0}^{2\dim(U)}G_d(\two U); \]
 a point $X\in\tG(U)$ should be thought of as a representative of the
 virtual vector space $X-U_-$.  We embed $G(U)$ in $\tG(U)$ by
 $X\mapsto X\op U=X_++U_-$.  We define $\tdim\:\tG(U)\xra{}\Z$ by
 $\tdim(X)=\dim(X)-\dim(U)$, and $\tG_d(U)=\{X\st\tdim(X)=d\}$.  Given
 an isometric embedding $j\:U\xra{}V$, we define
 $j_*\:\tG(U)\xra{}\tG(V)$ by $j_*(X)=(j\op j)(X)+W_-$, where
 $W=V\ominus j(U)$.  There is an evident map
 $\sg\:\tG(U)\tm\tG(V)\xra{}\tG(U\op V)$ sending $(X,Y)$ to $X\op Y$;
 one checks that $\tdim(X\op Y)=\tdim(X)+\tdim(Y)$ and that the map
 $\sg$ is compatible in an obvious sense with the maps $j_*$.  

 If $\CU$ is an infinite $A$-universe, we define $\two\CU=\CU\op\CU$
 as before, and put $\tG(\CU)=\colim_U\tG(U)$, where $U$ runs over
 finite-dimensional subspaces.  Equivalently, $\tG(\CU)$ is the space
 of subuniverses $\CV<\two\CU$ such that the space $\CV\cap\CU_-$ has
 finite codimension in $\CV$ and also has finite codimension in
 $\CU_-$.  This is a natural analogue of the space $\Z\tm BU$.
\end{definition}

\begin{proposition}
 For any $B\leq A$ we have
 \begin{align*}
  (G\CU)^B   &= \prod_{\bt\in B^*} G(\CU[\bt])
              = \Map(B^*,\coprod_dBU(d)) \\
  (\tG\CU)^B &= \prod_{\bt\in B^*} \tG(\CU[\bt])
              = \Map(B^*,\Z\tm BU)
 \end{align*}
 where 
 \[ \CU[\bt] =
     \{u\st b.u=\exp(2\pi i\bt(b))u \text{ for all } b\in B\} 
 \]
 is the $\bt$-isotypical part of $\CU$.  In each case, the first
 equivalence is $A/B$-equivariant, but the second is not.
\end{proposition}
\begin{proof}
 For the first isomorphism, just note that $\CU$ splits
 $A$-equivariantly as $\bigoplus_\bt\CU[\bt]$, and a subspace $V<\CU$
 is $B$-invariant iff it is the direct sum of its intersections with
 the subspaces $\CU[\bt]$.  This gives an $A/B$-equivariant
 isomorphism $(G\CU)^B=\prod_\bt G(\CU[\bt])$, and it is clear that
 $G(\CU[\bt])$ is nonequivariantly a copy of $\coprod_dBU(d)$ so
 $(G\CU)^B=\Map(B^*,\coprod_dBU(d))$.  The argment for $(\tG\CU)^B$ is
 essentially the same.
\end{proof}

We next write $R^+A=\N[A^*]=\pi^A_0(G\CU)$ for the additive semigroup
of honest representations of $A$, and $RA=\Z[A^*]=\pi^A_0(\tG\CU)$ for
the additive group of virtual representations.  It is clear that the
semigroup ring $E_0[R^+A]$ is a polynomial algebra over $E_0$ with one
generator $u_\al$ for each character $\al$, and the group ring
$E_0[RA]$ is the Laurent series ring with the same generators.  In
other words, we have
\begin{align*}
 E_0[R^+A] &= E_0[u_\al\st \al\in A^*] \\
 E_0[RA]   &= E_0[u_\al,u_\al^{-1}\st\al\in A^*] = E_0[R^+A][v^{-1}]
\end{align*}
where $v=\prod_\al u_\al$.  Note that 
\begin{align*}
 \spec(E_0[R^+A]) &= \Map(A^*,\aff^1) \\
 \spec(E_0[RA])   &= \Map(A^*,\MG),
\end{align*}
and the isomorphisms $R^+A=\pi^A_0G\CU$ and $RA=\pi^A_0\tG\CU$ give
maps $E_0[R^+A]\xra{}E_0G\CU$ and $E_0[RA]\xra{}E_0\tG\CU$.

Now let $\phi$ be the obvious isomorphism 
\[ \C[A]\op\CU = \C[A]\op\C[A]^\infty \xra{}\C[A]^\infty = \CU, \]
and define $\phi'\:G_d\CU\xra{}G_{d+|A|}\CU$ by
$\phi'(X)=\phi(\C[A]\op X)$.

\begin{proposition}
 The space $\tG\CU$ is the telescope of the self-map $\phi'$ of
 $G\CU$.  We thus have 
 \[ E_0\tG\CU = v^{-1}E_0G\CU = E_0[RA]\ot_{E_0[R^+A]} E_0G\CU, \]
 and so $\spec(E_0\tG\CU)=\Map(C,\MG)$.
\end{proposition}
\begin{proof}
 Put $\CU'=\C[A][z,z^{-1}]$, and identify this with $\two\CU$ by
 sending $(e_k,0)$ to $z^k$ and $(0,e_k)$ to $z^{-k-1}$.  The standard
 embedding $G\CU\xra{}\tG\CU$ now sends $X$ to $X\oplus\CU_-$.  It is
 easy to check that $\tG\CU=\colim_k z^{-k}G\CU$ on the nose, and that
 the inclusion $z^{-k}G\CU\xra{}z^{-k-1}G\CU$ is isomorphic to the map
 $\phi'$.  The first claim follows, and the second claim is just
 the obvious consequence in homology.  The tensor product description
 of $E_0\tG\CU$ gives us a pullback square
 \[ \xymatrix{
  {\spec(E_0\tG\CU)}
   \rto
   \dto &
  {\spec(E_0G\CU)=\Map(C,\aff^1)}
   \dto \\
  {\Map(A^*,\MG)=\spec(E_0[RA])}
   \rto &
  {\spec(E_0[R^+A])=\Map(A^*,\aff^1).}
  } \]
 As $C$ is a formal neighbourhood of the image of $\phi$, we see that
 a map $C\xra{f}\aff^1$ lands in $\MG$ if and only if the composite
 $A^*\tm S\xra{\phi}C\xra{f}\aff^1$ lands in $\MG$.  Given this, we
 see that the pullback is just $\Map(C,\MG)$ as required.
\end{proof}

We next introduce the analogue of $BU$.
\begin{proposition}\label{prop-Z-factor}
 There is a natural splitting $\tG\CU=\Z\tm\tG_0\CU$, and we have 
 $\spec(E_0\tG_0\CU)=\Map_0(C,\MG)$ (the scheme of maps
 $f\:C\xra{}\MG$ with $f(0)=1$).
\end{proposition}
\begin{proof}
 We have already described an equivariant map $\tdim\:\tG\CU\xra{}\Z$,
 and defined $\tG_0\CU=\ker(\tdim)$.  We also have
 $(\tG\CU)^A=\Map(A^*,\Z\tm BU)$ so $\pi^A_0(\tG\CU)=\Map(A^*,\Z)=RA$,
 which gives an equivariant map $i\:RA\xra{}\tG\CU$ (where $RA$ has
 trivial action).  The composite $\tdim\circ i\:RA\xra{}\Z$ is just
 the usual augmentation map $\ep$ sending a virtual representation to
 its dimension.  Thus, if we let $\eta\:\Z\xra{}RA$ be the unit map,
 then $i\circ\eta$ is a section of $\tdim$.  As $\tG\CU$ is a
 commutative equivariant $H$-space, we can define a map
 $\dl\:\tG\CU\xra{}\tG_0\CU$ by $x\mapsto(i(\eta(\tdim(x))) - x)$, and we
 find that the resulting map $(\tdim,\dl)\:\tG\CU\xra{}\Z\tm\tG_0\CU$
 is an equivalence.  This is easily seen to be parallel to the
 splitting $\Map(C,\MG)=\MG\tm\Map_0(C,\MG)$ given by
 $f\mapsto(f(0),f(0)/f)$, which gives the claimed description of
 $\spec(E_0G\CU)$.
\end{proof}

\begin{remark}\label{rem-BU-versions}
 There are two other possible analogues of $BU$.  Firstly, one could
 take the colimit of the spaces $G_d\CU$ using the maps
 $V\mapsto V\op\C$; the scheme associated to the corresponding space
 is then $\Map_0(C,\aff^1)$, which classifies maps $f\:C\xra{}\aff^1$
 with $f(0)=1$.  Alternatively, one could take the colimit of the
 spaces $G_{d|A|}\CU$ using the maps $V\mapsto V\op\C[A]$.  This gives
 the scheme of maps $f\:C\xra{}\MG$ for which
 $\prod_{\al\in A^*}f(\phi(\al))=1$.  However, the space $\tG_0\CU$
 described above is the one that occurs in Greenlees's definition of
 the spectrum $kU_A$, and is also the one whose Thom spectrum is
 $MU_A$. 
\end{remark}

We next introduce the analogue of $BSU$.  For this, we need an
analogue of the map $B\det\:BU\xra{}\CPi$.
\begin{definition}\label{defn-tdet}
 Given a universe $U$ of finite dimension $d$, we put 
 $\tF U=\Hom(\lm^dU_-,\lm^d(\two U))$.  We make this a functor as
 follows.  Given an isometric embedding $j\:U\xra{}V$, we put
 $W=V\ominus jU$ and $e=\dim(W)$.  As $j\:U\xra{}jU$ is an isomorphism
 and $\lm^eW$ is one-dimensional, we have an evident isomorphism 
 \[ \tF U=\Hom(\lm^djU_-\ot\lm^eW_-,\lm^d(\two jU)\ot\lm^eW_-). \]
 The isomorphism $V=jU\op W$ gives an isomorphism
 $\lm^djU_-\ot\lm^eW=\lm^{d+e}V$ and an embedding
 $\lm^d(\two jU)\ot\lm^eW_-\xra{}\lm^{d+e}(\two V)$.  Putting this
 together gives the required map $j_*\:\tF U\xra{}\tF V$.

 There are also obvious maps $\tF(U)\ot\tF(V)\xra{}\tF(U\oplus V)$,
 compatible with the above functorality.  This gives maps
 $P\tF(U)\tm P\tF(V)\xra{}P\tF(U\op V)$ of the associated projective
 spaces. 

 Next, recall that a point of $\tG_0U$ is a $d$-dimensional subspace
 $X\leq\two U$.  We define 
 \[ \tdet(X) = \Hom(\lm^dU_-,\lm^dX) \in P\tF U. \]
 One can check that this gives a natural map
 $\tdet\:\tG_0\xra{}P\tF$, with $\tdet(X\op Y)=\tdet(X)\ot\tdet(Y)$.
 
 Finally, for our complete universe $\CU$ we put
 $\tF\CU=\colim_U\tF U$, where $U$ runs over the finite-dimensional
 subuniverses.  It is easy to check that this is again a complete
 $A$-universe, and thus is unnaturally isomorphic to $\CU$.  The maps
 $\tdet$ pass to the colimit to give an $H$-map
 $\tdet\:\tG_0\CU\xra{}P\tF\CU\simeq P\CU$.  We write $S\tG_0\CU$ for
 the pullback of the projection $S(\tF\CU)\xra{}P\tF\CU$ along the map
 $\tdet$, or equivalently the space of pairs $(\CV,u)$ where
 $\CV\in\tG_0\CU$ and $u$ is a unit vector in the one-dimensional
 space $\tdet(\CV)$.  As $S(\tF\CU)$ is equivariantly contractible,
 this is just the homotopy fibre of $\tdet$.
\end{definition}

\begin{proposition}
 There is a natural splitting $\tG_0\CU=S\tG_0\CU\tm P\CU$ (which does
 not respect the $H$-space structure).
\end{proposition}
\begin{proof}
 It is enough to give a section of the $H$-map
 $\tdet\:\tG_0\CU\xra{}P\CU$.  We can include $P\CU=G_1\CU$ in
 $\tG_1\CU\subset\tG\CU$ in the usual way, then use the projection
 $\tG\CU\xra{}\tG_0\CU$ from Proposition~\ref{prop-Z-factor}.  We find
 that the resulting composite $P\CU\xra{}P\CU$ is actually minus the
 identity, but we can precompose by minus the identity to get the
 required section.
\end{proof}

\begin{remark}
 Cartier duality identifies $\spec(E_0P\CU)$ with $\Hom(C,\MG)$, and
 the proposition suggests that $\spec(E_0S\tG_0\CU)$ should be the
 quotient $\Map_0(C,\MG)/\Hom(C,\MG)$.  However, there are
 difficulties in interpreting this quotient, and it is in fact more
 useful to take a slightly different approach as
 in~\cites{hoanst:esw,anst:wpm}.  We will not give details here.
\end{remark}

Next, recall that Greenlees has defined an equivariant analogue of
connective $K$-theory (denoted by $kU_A$) by the homotopy pullback
square 
\[ \xymatrix{
 {kU_A} \rto \dto & {F(EA_+,kU)} \dto \\
 {KU_A}  \rto & {F(EA_+,KU).}
 } \]
If $v\in\pi_2kU$ is the Bott element then $kU[v^{-1}]=KU$ and
$kU/v=H$.  It is not hard to see that there is a corresponding element
in $\pi^A_2kU_A$ with $kU_A[v^{-1}]=KU_A$ and $kU_A/v=F(EA_+,H)$.  

\begin{proposition}
 The zeroth, second and fourth spaces of $kU_A$ are $\tG\CU$,
 $\tG_0\CU$ and $S\tG_0\CU$ respectively.
\end{proposition}
\begin{proof}
 We take it as well-known that the zeroth space of $KU_A$ is $\tG\CU$,
 and $KU_A$ is two-periodic so this is also the $2k$'th space for all
 $k$.  Let $X_k$ denote the $2k$'th space of $kU_A$, so we have a
 homotopy pullback square
 \[ \xymatrix{
  {X_k} \rto \dto^{i} & {F(EA_+,\BU{2k})} \dto \\
  {\tG\CU} \rto_(0.3){j} & {F(EA_+,\Z\tm BU)}
  } \]
 (where $\BU{0}$ is interpreted as $\Z\tm BU$).  In the case $k=0$,
 the map $i$ is the identity and so $X_0=\tG\CU$.  In the case $k=1$,
 the map $i$ is just the inclusion
 \[ F(EA_+,BU) \xra{} \Z \tm F(EA_+,BU) = F(EA_+,\Z\tm BU) \]
 and the map $j$ sends $\tG_k\CU$ into $\{k\}\tm F(EA_+,BU)$.  It
 follows easily that $X_1=\tG_0\CU$.  In the case $k=2$, we note that
 the cofibration $\Sg^2kU_A\xra{v}kU_A\xra{}F(EA_+,H)$ gives a
 fibration $X_2\xra{}X_1\xra{}F(EA_+,K(\Z,2))$.  We know that
 $X_1=\tG_0\CU$ and Proposition~\ref{prop-PU-FEACP} that
 $F(EA_+,K(\Z,2))=P\CU$.  One can check that the resulting map
 $\tG_0\CU\xra{}P\CU$ is just $\pm\tdet$, and so $X_2=S\tG_0\CU$ as
 claimed. 
\end{proof}

\section{Transfers and the Burnside ring}
\label{sec-transfer}

There is a well-known relationship between transfers, the Burnside
ring, and equivariant stable homotopy.  In this section we will show
how this relationship is encoded in the theory of equivariant formal
groups. 

Let $\Om=\Om(A)$ be the Burnside ring of $A$, which is the group
completion of the monoid of isomorphism classes of finite $A$-sets.
Additively this is freely generated over $\Z$ by the elements $[A/B]$
for subgroups $B\leq A$, with the product rule
\[ [A/B].[A/B'] = \left|\frac{A}{BB'}\right|\;[A/B\cap B']. \]
There is a well-known isomorphism $\Om\to\pi^G_0$, sending $[A/B]$ to
$\trf_B^A(1)$.  

Now let $E$ be an $E$-equivariant evenly periodic ring spectrum, with
associated equivariant formal group $(C,\phi)$ over the scheme
$S=\spec(E_0)$.  The unit map $\eta\:S_G\to E$ gives us a ring map
$\Om\to\OS=E_0$.  We will show that this map is determined
by the EFG structure.

\subsection{The functions $v_n$}

Let $C$ be a formal multicurve group over a scheme $S$, and let $x$ be
a coordinate on $C$.  (Later we will show that our constructions are
independent of $x$, but we have not been able to make this visible
from the outset.)  It will be convenient to write $x_n(a)=x(na)$.

\begin{lemma}
 There is a unique function $f\:C\tm_SC\to\aff^1$ such that 
 \[ x(a+b) = x(a) + x(b) + x(a)x(b)f(a,b) \]
 for all $(a,b)\in C\tm_SC$.
\end{lemma}
\begin{proof}
 From the axioms for formal multicurves, we see that
 $\OC=\OS\oplus\OC.x$.  After tensoring this decomposition with
 itself, we see that and function $p$ on $C\tm_SC$ can be written
 uniquely in the form 
 \[ p(a,b) = q + x(a)r(b) + s(a)x(b) + f(a,b)x(a)x(b), \]
 for some $q\in\OS$ and some functions $r$, $s$ and $f$.  If
 $p(a,0)=p(0,b)=0$ then it follows that $q$, $r$ and $s$ vanish, so
 $p(a,b)=f(a,b)x(a)x(b)$.  We apply this to the function
 $p(a,b)=x(a+b)-x(a)-x(b)$ to prove the lemma.
\end{proof}

We can apply the lemma inductively to understand $x_n(a)$.  For any
finite set $I$ of integers we put $y_I(a)=\prod_{i\in I}f(a,ia)$.  We
also write $\min_n(I)$ for the smallest element of $I$, with the
convention $\min_n(\emptyset)=n$.  

\begin{lemma}
 $x_n=\sum_I\min_n(I)x^{|I|+1}y_I$, where $I$ runs over subsets of
 $\{1,2,\dotsc,n-1\}$.  
\end{lemma}
\begin{proof}
 Induction, using the identity
 $x_{n+1}=x+x_n+x x_n y_n=x+(1+x y_n)x_n$.
\end{proof}

\begin{definition}\label{defn-vw}
 We put 
 \begin{align*}
  v_n &= \sum_I \min_n(I) x^{|I|}y_I \\
  w_n &= \sum_{I\neq\emptyset} \min_n(I) x^{|I|-1}y_I,
 \end{align*}
 so that
 \[ x_n = x\,v_n = nx + x^2 w_n. \]
 We also write $v_{n,k}(a)=v_n(ka)$ and $w_{n,k}(a)=w_n(ka)$ and
 $y_{I,k}(a)=y_I(ka)$.  
\end{definition}

\begin{lemma}\label{lem-nm}
 $v_{nm}=v_n v_{m,n}=v_m v_{n,m}$.
\end{lemma}
\begin{proof}
 By evaluating the identity $x(ma)=x(a)v_m(a)$ at $a=nb$, we get
 $x_{nm}=x_nv_{m,n}=xv_nv_{m,n}$.  On the other hand, we have
 $x_{nm}=x\,v_{nm}$ and $x$ is not a zero divisor so
 $v_{nm}=v_nv_{m,n}$.  The other identity follows symmetrically. 
\end{proof}

\begin{corollary}
 $v_n^2=nv_n\pmod{x_n}$.
\end{corollary}
\begin{proof}
 We have $v_n^2-nv_n=xw_nv_n=w_nx_n$.
\end{proof}

We can push the lemma one step further, as follows.
\begin{lemma}\label{lem-asym}
 $v_n^2w_{m,n}-v_m^2w_{n,m}=n\,w_m-m\,w_n=v_nw_m-v_mw_n$.
\end{lemma}
\begin{proof}
 We have 
 \begin{align*}
  x_{nm} &= nx_m + x_m^2 w_{n,m} \\
   &= nmx + nx^2w_m + x^2 v_m^2 w_{n,m} \\
\intertext{
 By exchanging the roles of $n$ and $m$, we also have 
}
  x_{nm} &= nmx + mx^2w_n + x^2 v_n^2 w_{m,n}.
 \end{align*}
 Subtracting these expressions gives
 \[ x^2(nw_m-mw_n+v_m^2w_{n,m}-v_n^2w_{m,n}) = 0. \]
 As $x$ is not a zero-divisor in $\O_C$, it follows that 
 \[ v_n^2w_{m,n}-v_m^2w_{n,m}=n\,w_m-m\,w_n. \]
 Now write $v_n=n+xw_n$ and $v_m=m+xw_m$ to see that 
 \[ nw_m-mw_n = v_nw_m-v_mw_n. \]
\end{proof}

\begin{lemma}\label{lem-coprime}
 Suppose that $n$ and $m$ are coprime.  Then $v_n=v_{n,m}\pmod{x_n}$.
\end{lemma}
\begin{proof}
 Put 
 \[ s_{n,m} = v_mw_{n,m}-w_n. \]
 Note that 
 \[ v_{n,m}-v_n = 
   (n+x_mw_{n,m}) - (n+xw_n) = x_mw_{n,m} - xw_n = 
    x(v_mw_{n,m}-w_n) = x\,s_{n,m}.
 \]
 On the other hand, a simple rearrangement of Lemma~\ref{lem-asym}
 tells us that $v_ms_{n,m}=v_ns_{m,n}$, so $x_m\,s_{n,m}=x_ns_{m,n}$.

 Put $I=\{t\in\O_C\st t\,s_{n,m}\in(x_n)\}$.  As
 $x\,s_{n,m}=v_{n,m}-v_n$, it will suffice to show that $x\in I$. 

 As $n$ and $m$ are coprime, an element $a\in C$ has $a=0$ iff
 $na=ma=0$.  This means that the vanishing locus of $x$ is the
 intersection of the vanishing loci of $x_n$ and $x_m$, or in other
 words that $(x)=(x_n)+(x_m)$.  Visibly $x_n\in I$, so it will suffice
 to show that $x_m\in I$.  This follows from the identity
 $x_m\,s_{n,m}=x_ns_{m,n}$. 
\end{proof}

\begin{corollary}\label{cor-bicyclic}
 If $n$ and $m$ are coprime then $v_{nm}=v_{n,m}v_{m,n}\pmod{x_{nm}}$.
\end{corollary}
\begin{proof}
 Lemma~\ref{lem-nm} tells us that $v_{nm}=v_nv_{m,n}$, so it will
 suffice to show that $(v_n-v_{n,m})v_{m,n}$ is divisible by the
 element $x_{nm}=x_nv_{m,n}$, or to show that $v_n-v_{n,m}$ is
 divisible by $x_n$.  This is just Lemma~\ref{lem-coprime}. 
\end{proof}

\begin{proposition}\label{prop-vn-invariant}
 Let $x'$ be another coordinate on $C$, and let $v'_n$ be the unique
 function such that $x'_n(a)=x'(na)=x'(a)v'_n(a)$ for all $a$.  Then
 $v'_n(a)=v_n(a)$ whenever $na=0$.
\end{proposition}
\begin{proof}
 Both $x$ and $x'$ generate the ideal of functions that vanish on the
 zero section, so we must have $x=px'$ and $x'=qx$ for some functions
 $p$ and $q$.  This implies that $(1-pq)x=0$ but $x$ is regular, so
 $pq=1$.  Now 
 \[ v'_n(a)x'(a) = x'(na)=q(na)x(na)=q(na)v_n(a)x(a) = 
     q(na)p(a)v_n(a)x'(a),
 \]
 and $x'$ is regular so $v'_n(a)=q(na)p(a)v_n(a)$ for all $a$.  If
 $na=0$ this gives $v'_n(a)=q(0)p(a)v_n(a)$.  Moreover, we have
 $p(a)=p(0)\pmod{x(a)}$ and $x(a)v_n(a)=x(na)=0$ so
 $p(a)v_n(a)=p(0)v_n(a)$.  We also have $q(0)p(0)=1$ so
 $v'_n(a)=v_n(a)$ as claimed.
\end{proof}

Given any coordinate $x$ we could apply the above to the coordinate
$x'(a)=x(-a)$.  For that case, however, some additional things can be
said, as we now explain.
\begin{proposition}
 If $na=0$ then $v_{-n}(a)=-v_n(-a)$.
\end{proposition}
\begin{proof}
 Take $m=-1$ 
\end{proof}

\subsection{Transfer elements}

\begin{definition}
 Let $U$ be a finite abelian group, and let $\phi\:U\tm S\to C$ be a
 map of group schemes over $S$.  We write $I(U)$ or $I(\phi)$ or
 $I(U,\phi)$ for the ideal in $\OS$ generated by the elements
 $x(\phi(u))$ for $u\in U$, so $\spec(\OS/I(U))$ is the largest closed
 subscheme of $S$ over which $\phi$ vanishes.
\end{definition}

\begin{theorem}\label{thm-transfer}
 There is a unique way to define elements $t(U)=t(U,\phi)\in\OS$ with
 properties as follows. 
 \begin{enumerate}
  \item[(a)] If $\phi=0$ then $t(U,0)=|U|$; in particular, $t(0,0)=1$.
  \item[(b)] If $U=U_0\oplus U_1$ then $t(U)=t(U_0)t(U_1)$.
  \item[(c)] If $U$ is cyclic of order $n$, generated by $\al$, then
   $t(U)=v_n(\phi(\al))$.
  \item[(d)] If $\lm\:U'\to U$ is an isomorphism then
   $t(U',\phi\lm)=t(U,\phi)$. 
  \item[(e)] $t(U)I(U)=0$.
  \item[(f)] Suppose that $V\leq U$, and that $\ov{\phi}$ is the
   induced map $U/V\to C$ defined over $\spec(\OS/I(V))$.  Then
   $t(U,\phi)=t(V,\phi)t(U/V,\ov{\phi})$.  (Note that~(e) makes the
   right hand side well-defined.)
  \item[(g)] If $V,W\leq U$ then $t(V)t(W)=|V\cap W|t(V+W)$.
 \end{enumerate}
\end{theorem}
The proof will be given at the end of this section, after some
preliminary results.  Note that the elements $v_n(\phi(\al))$ in
clause~(c) are independent of the coordinate, by
Proposition~\ref{prop-vn-invariant}. 

It is clear that there is at most one way to satisfy~(a) to~(d); the
problem is that this gives a well-defined answer which also has
properties~(e) to~(g).  The next definition makes this more formal.

\begin{definition}\label{defn-presentation}
 A \emph{presentation} of a subgroup $U$ is a set $P$ of nonzero
 elements of $U$ such that the resulting map
 $\bigoplus_{\al\in P}\Z/\ord(\al)\to U$ is an isomorphism.  Given a
 presentation $P$, we put 
 \[ t(P) = t(P,\phi) = \prod_{\al\in P} v_{\ord(\al)}(\phi(\al)). \]
 We say that $U$ is \emph{canonical} if $t(P)$ is independent of $P$;
 if so, we write $t(U)$ instead of $t(P)$.
\end{definition}

\begin{lemma}\label{lem-tI}
 For any presentation $P$ of $U$ we have $t(P)I(U)=0$.
\end{lemma}
\begin{proof}
 Suppose that $\al\in P$ has order $n$, and put $a=\phi(\al)$.  Then
 $v_n(a)$ is a factor in $t(P)$ and $x(a)v_n(a)=x(na)=x(0)=0$, so
 $x(a)t(P)=0$.  These elements $x(a)$ generate $I(U)$, so
 $t(P)I(U)=0$. 
\end{proof}

\begin{lemma}\label{lem-cyclic-i}
 If $U$ is cyclic of order $n$ and $\al$ is a generator, then
 $t(\{\al\})$ is independent of $\al$.  (We can thus write $t_1(U)$
 for $t(\{\al\})$.)
\end{lemma}
\begin{proof}
 If $\bt$ is another generator then $\bt=m\al$ for some $m$ that is
 coprime to $n$.  If $a=\phi(\al)$ then $t(\{\al\})=v_n(a)$ and
 $t(\{\bt\})=v_{n,m}(a)$, but these are the same by
 Lemma~\ref{lem-coprime}.  
\end{proof}

\begin{lemma}\label{lem-cyclic-ii}
 If $U$ is cyclic and $U=\prod_{i=1}^rU_i$ then
 $t_1(U)=\prod_it_1(U_i)$. 
\end{lemma}
\begin{proof}
 By induction we can reduce to the case $r=2$.  Put $n=|U_1|$ and
 $m=|U_2|$.  As $U$ is cyclic, these must be coprime.  Choose a
 generator $\al$ for $U$ and put $a=\phi(\al)$.  Note that $m\al$
 generates $U_1$ and $n\al$ generates $U_2$, so $t_1(U_1)=v_{n,m}(a)$
 and $t_1(U_2)=v_{m,n}(a)$, whereas $t_1(U)=v_{nm}(a)$.  The claim now
 follows from Corollary~\ref{cor-bicyclic}.
\end{proof}

\begin{corollary}\label{cor-cyclic-canonical}
 Cyclic groups are canonical.
\end{corollary}
\begin{proof}
 Let $U$ be cyclic, and let $P=\{\al_1,\dotsc,\al_r\}$ be a
 presentation, and let $U_i$ be the subgroup generated by $\al_i$.
 Directly from the definitions we see that $t(P)=\prod_it_1(U_i)$, but
 the lemma tells us that this is the same as $t_1(U)$ and thus is
 independent of $P$.
\end{proof}

\begin{lemma}\label{lem-split-canonical}
 Suppose that $U=\prod_pU(p)$, where $U(p)$ is a $p$-group, and that
 each $U(p)$ is canonical.  Then $U$ is canonical, with
 $t(U)=\prod_pt(U(p))$. 
\end{lemma}
\begin{proof}
 Let $P=\{\al_1,\dotsc,\al_d\}$ be a presentation of $U$.  Let $U_i$
 be the cyclic subgroup generated by $\al_i$, and put $n_i=|U_i|$ and
 $t_i=t(U_i)=v_{n_i}(\phi(\al_i))$.  Clearly $t(P)=\prod_it_i$, so it
 will suffice to prove that this is equal to $\prod_pt(U(p))$.

 Now let $U_i(p)$ be the $p$-torsion part of $U_i$.  Put
 $n_{i,p}=|U_i(p)|$ and choose an element $\bt_{i,p}$ generating
 $U_i(p)$.  Put $t_{i,p}=t(U_i(p))=v_{n_{i,p}}(\phi(\bt_{i,p})$.  If
 we fix $i$ then the elements $\bt_{i,p}$ give a presentation of
 $U(i)$, so $\prod_pt_{i,p}=t(U_i)=t_i$, so
 $\prod_{i,p}t_{i,p}=\prod_it_i$.  On the other hand, if we fix $p$
 then the elements $\bt_{i,p}$ present $U(p)$, so
 $\prod_it_{i,p}=t(U(p))$, so $\prod_{i,p}t_{i,p}=\prod_pt(U(p))$.
 Thus $\prod_it_i=\prod_pt(U(p))$ as required.
\end{proof}

\begin{lemma}\label{lem-elementary-canonical}
 Elementary abelian $p$-groups are canonical.
\end{lemma}
\begin{proof}
 Suppose $U=(\Z/p)^d$.  The element $t(\{\al_1,\dotsc,\al_d\})$ is
 clearly unchanged by permuting the elements $\al_i$, and
 Lemma~\ref{lem-cyclic-i} tells us that it is also unchanged if we
 replace each $\al_i$ by $m_i\al_i$ for some $m_i\in(\Z/p)^\tm$.
 Next, we claim that if $a,b\in C$ and $pa=pb=0$ then
 $v_p(a+b)v_p(b)=v_p(a)v_p(b)$.  Indeed, we have
 $x(b)v_p(b)=x_p(b)=0$, so it will suffice to show that
 $v_p(a+b)=v_p(a)\pmod{x(b)}$.  This is clear because the vanishing
 locus of $x$ is the zero section.  It now follows that
 $t(\{\al_1,\dotsc,\al_d\})$ is unchanged by substitutions of the form
 $\al_i\mapsto\al_i+\al_j$ with $j\neq i$.  These substitutions,
 together with permutations and diagonal matrices, generate
 $GL_d(\Z/p)$, and the claim follows.
\end{proof}

\begin{lemma}\label{lem-p-canonical}
 If $U$ is an abelian $p$-group then $U$ is canonical.  Moreover, if
 we put $V=\{\al\in U\st p\al=0\}$, then
 $t(U,\phi)=t(V,\phi)t(U/V,\ov{\phi})$.
\end{lemma}
\begin{proof}
 By induction on the exponent (starting with
 Lemma~\ref{lem-elementary-canonical}), we may assume that $V$ and
 $U/V$ are canonical.  We therefore have well-defined elements
 $t(V,\phi)\in\OS$ and $t(U/V,\ov{\phi})\in\OS/I(V,\phi)$, with
 $t(V,\phi)I(V,\phi)=0$ by Lemma~\ref{lem-tI}.  This means that there
 is a well-defined product $t_1=t(V,\phi)t(U/V,\ov{\phi})\in\OS$.  It
 will suffice to show that $t(P)=t_1$ for any presentation
 $P=\{\al_1,\dotsc,\al_d\}$ of $U$.  To see this, let the order of
 $\al_i$ be $p^{r_i}$, so $r_i>0$.  We will assume for simplicity that
 in fact $r_i>1$ for all $i$; only minor notational adjustments are
 needed for the remaining cases.  Put $\bt_i=p^{v_i-1}\al_i$, so that
 $\{\bt_1,\dotsc,\bt_d\}$ is a presentation of $V$.  Let $\ov{\al}_i$
 be the image of $\al_i$ in $U/V$, so
 $\{\ov{\al}_1,\dotsc,\ov{\al}_d\}$ is a presentation for $U/V$.  Put
 $a_i=\phi(\al_i)\in C$.  Now $t(P)$ is the product of the terms
 $v_{p^{v_i}}(a_i)$, whereas $t(U/V,\ov{\phi})$ is the product of the
 terms $v_{p^{v_i-1}}(a_i)$, and $t(V,\phi)$ is the product of the
 terms $v_p(p^{v_i-1}a_i)$.  Lemma~\ref{lem-nm} tells us that
 $v_{p^{v_i}}(a_i)=v_{p^{v_i-1}}(a_i)v_p(p^{v_i-1}a_i)$, so
 $t(P)=t(U/V,\ov{\phi})t(V,\phi)$ as required.
\end{proof}

\begin{corollary}
 All finite abelian groups are canonical.
\end{corollary}
\begin{proof}
 This is immediate from Lemmas~\ref{lem-p-canonical}
 and~\ref{lem-split-canonical}. 
\end{proof}

\begin{definition}
 Let $U$ be a finite abelian group, and let $V$ be a subgroup.  We say
 that the pair $(U,V)$ is \emph{good} if for all $C$ and $\phi$ as
 above, we have $t(U,\phi)=t(V,\phi)t(U/V,\ov{\phi})$.
\end{definition}
\begin{remark}
 From now on we will just write $t(U)$ for $t(U,\phi)$ and so on, and
 not mention explicitly the quotient rings in which these elements
 lie. 
\end{remark}

\begin{lemma}\label{lem-good-omni}
 \begin{itemize}
  \item[(a)] If $W\leq V\leq U$ andthe pairs $(V,W)$ and $(U,W)$ and
   $(U/W,V/W)$ are good then so is $(U,V)$.
  \item[(b)] If the pairs $(U_i,V_i)$ are good then so is
   $(\bigoplus_iU_i,\bigoplus_iV_i)$.
  \item[(c)] If the pairs $(U(p),V(p))$ are all good (where
   $U(p)=\{\al\in U\st p^N\al=0\text{ for }N\gg 0\}$) then so is
   $(U,V)$. 
 \end{itemize}
\end{lemma}
\begin{proof}
 For~(a), we can use the goodness of $(U,W)$, $(U/W,V/W)$ and $(V,W)$
 in turn to see that
 \[ t(U) = t(W)t(U/W) = t(W)t(V/W)t(U/V) = t(V)t(U/V). \]
 Part~(b) is clear, and part~(c) is a special case.
\end{proof}

\begin{lemma}\label{lem-good-start}
 If $|V|=p$ then $(U,V)$ is good.
\end{lemma}
\begin{proof}
 We can use Lemma~\ref{lem-good-omni}(c) to reduce to the case where
 $U$ is a $p$-group.  If $U$ is cyclic, the claim then follows from
 Lemma~\ref{lem-p-canonical}.

 If $U$ is not cyclic, let the exponent of $U$ be $p^u$, and choose a
 cyclic summand $W\leq U$ of order $p^u$.  If $V\leq W$ then we choose
 a complementary group $X$ with $U=W\op X$, note that the pairs
 $(W,V)$ and $(X,0)$ are good, and apply Lemma~\ref{lem-good-omni}(b) to
 these pairs.  

 Suppose instead that $V\not\leq W$, so $V\op W\leq U$.  Note that all
 groups mentioned are modules over the ring $R=\Z/p^u$.  Now
 $R\simeq\Hom(R,\Q/\Z)$, so $R$ is self-injective, so $W$ is an
 injective $R$-module.  This means that the projection
 $\pi\:V\op W\to W$ can be extended to give a map $\pi\:U\to W$.  If
 we put $X=\ker(\pi)$ we see that $U=W\op X$ with $V\leq X$.  The pair
 $(X,V)$ can be assumed to be good by induction on the order, and the
 pair $(W,0)$ is good, so the pair $(U,V)=(W,0)\op(X,V)$ is good.
\end{proof}

\begin{corollary}\label{cor-good}
 All pairs $(U,V)$ are good.
\end{corollary}
\begin{proof}
 We can use Lemma~\ref{lem-good-omni}(c) to reduce to the case where
 $U$ is a $p$-group.   We then argue by induction on the order of $V$,
 starting with Lemma~\ref{lem-good-start}.  Choose a subgroup
 $W\leq V$ of order $p$.  The pairs $(U,W)$ and $(V,W)$ are good by
 the lemma, and $(U/W,V/W)$ is good by induction.  Part~(a) of
 Lemma~\ref{lem-good-omni} therefore tells us that $(U,V)$ is good.
\end{proof}

\begin{proposition}\label{prop-t-prod}
 If $V,W\leq U$ then $t(V)t(W)=|V\cap W|t(V+W)$.
\end{proposition}
\begin{proof}
 Put $X=V\cap W$.  We have $t(V)=t(V/X)t(X)$ and $t(W)=t(W/X)t(X)$ and
 $t(X)^2=|X|t(X)$ so $t(V)t(W)=|X|t(V/X)t(W/X)t(X)=|X|t(V/X)t(W)$.
 Now $V/X$ can be identified with $(V+W)/W$ and
 $t((V+W)/W)t(W)=t(V+W)$, so $t(V)t(W)=|X|t(V+W)$ as claimed.
\end{proof}

\begin{proof}[Proof of Theorem~\ref{thm-transfer}]
 Part~(a) follows from the definitions if we remember that $v_n(0)=n$.
 Parts (b) to (d) are built in to the definition of $t$.  Part~(e)
 follows from Lemma~\ref{lem-tI}, and parts~(f) and~(g) are
 Corollary~\ref{cor-good} and Proposition~\ref{prop-t-prod}.
\end{proof}

We now change notation, and consider an $A$-equivariant formal group
$(C,\phi)$, so $\phi\:A^*\tm S\to C$.  

\begin{definition}
 We define a map 
 \[ \eta\:\Om = \Om(A) = \Z\{[A/B]\st B\leq A\} \to \OS \]
 by $\eta([A/B])=t(\ann(B),\phi)$.
\end{definition}

\begin{proposition}\label{prop-eta-ring}
 $\eta$ is a ring map.
\end{proposition}
\begin{proof}
 For every point in $(A/B)\tm(A/C)$, the isotropy group is $B\cap C$.
 This means that in $\Om$ we have $[A/B].[A/C]=n[A/(B\cap C)]$, where
 $n$ is characterised by the fact that $|A/B|.|A/C|=n|A/(B\cap C)|$.
 Now put $V=\ann(B)=(A/B)^*$, so $|V|=|A/B|$.  Similarly, put
 $W=\ann(C)=(A/C)^*$, so $|W|=|A/C|$.  We see that
 $\ann(B\cap C)=V+W$ and $|V+W|=|A/(B\cap C)|$, so $n$ is defined by
 the equation $|V||W|=n|V+W|$, so $n=|V\cap W|$.  We have
 $\eta([A/B])=t(V)$ and $\eta([A/C])=t(W)$ and 
 \[ \eta([A/B][A/C]) = n\eta([A/(B\cap C)]) = |V\cap W|t(V+W). \]
 The claim now follows from Proposition~\ref{prop-t-prod}.
\end{proof}

\subsection{Comparison with topology}

Now consider the case where $C$ arises from an evenly orientable
$A$-equivariant ring spectrum $E$.  We then have a map
$\eta\:\Om(A)\to\OS=E_0$ as above, and also a topologically defined
map $\eta'\:\Om(A)=\pi_0^A(S)\to E_0$, given by
$\eta'([A/B])=\trf_B^A(1)$.  
\begin{proposition}\label{prop-eta-general}
 We have $\eta([A/B])=\trf_B^A(1)$, so $\eta=\eta'$.
\end{proposition}

We will do the most basic case as a separate lemma, and then do the
general case.
\begin{lemma}\label{lem-eta-special}
 If $\al\in A^*$ has order $m$ and $B=\ker(\al)$ then
 $\trf_B^A(1)=v_m(\phi(\al))$.  
\end{lemma}
This was originally proved by Quillen~\cite{qu:epc}, but we will give
an argument here for convenience.
\begin{proof}
 First, we need to deal with a sign issue.  Recall that any coordinate
 $x$ gives a natural system of Thom classes $u_V$ for all bundles $V$.
 We also have another coordinate $\xb(a)=x(-a)$, and we write $\ub_V$
 for the Thom class defined using $\xb$ instead of $x$.

 Now let $T$ be the tautological bundle over $P\CU$, and let $P\CU^T$
 be the corresponding Thom space.  Restriction to the zero section
 sends $u_T$ to $\xb$ and $\ub_T$ to $x$.

 Define $\tht_0\:P\CU^T\to P\CU^{T^m}$
 by $\tht_0(L,x)=(L,x^{\ot m})$.  We claim that
 $\tht_0^*(\ub_{T^m})=v_mu_T$.  (where $v_m$ is defined by $x_m=v_mx$ as
 previously).  Indeed, we know that $\tE^0(P\CU^T)$ is freely
 generated over $E^0P\CU$ by $\ub_T$, so certainly
 $\tht_0^*(\ub_{T^m})=f\ub_T$ for some $f$.  Next, note that $\tht$ is the
 identity on the zero section, and that restriction to the zero
 section converts Thom classes to Euler classes, and that the Euler
 class of $T^m$ is just $x_m$.  We deduce that $x_m=fx$, so $f=v_m$ as
 claimed.  We now restrict to the point $L_\al\in P\CU$, noting that
 $L_\al^{\ot m}\simeq L_0$ and that the restriction $\OC=E^0P\CU\to
 E^0PL_\al=\OS$ is just evaluation at $\phi(\al)$.  This gives us a
 map $\tht_0\:S^{L_\al}\to S^{L_0}$ with $\tht_0(z)=z^m$ and
 $\tht_0^*(u_{L_0})=v_m(\phi(\al))u_{L_\al}$.  In this context it is
 meaningful to define $\tht\:S^{L_\al}\to S^{L_0}$ by
 $\tht(z)=(z^m-1)/m$; this is equivariantly homotopic to $\tht_0$ and
 so also satisfies $\tht^*(\ub_{L_0})=v_m(\phi(\al))\ub_{L_\al}$. 

 We next need to relate this to the transfer.  Recall that
 $\al\in A^*=\Hom(A,\Q/\Z)$; we define $\al'(a)=\exp(2\pi i\al(a))$ so
 that $\al'\:A\to S^1$.  Let $q\:A\to A/B$ be
 the projection.  We have an equivariant embedding $j_0\:A/B\to L_\al$
 given by $j_0(q(a))=\al'(a)$.  To construct the transfer, we need an
 equivariant embedding $j\:A/B\tm L_\al\to L_\al$ with
 $j(q(a),z)\simeq j_0(q(a))+z=\al'(a)+z$ to first order in $z$.  There is a
 unique such $j$ up to isotopy, so we can choose any $j$ that is
 convenient for our purposes.  We will define a diffeomorphism
 \[ i \: \C\to \{z\in \C\st |z|<1/m\} \]
 by $i(z)=z/(1+m|z|)$, and then define $j$ by 
 \[ j(q(a),z) = \al'(a) (1 + m\,i(z/\al'(a)))^{1/m}. \]
 (Here we note that $1+m\,i(z/\al'(a))$ lies in the right half plane,
 and we implicitly use the principal branch of the $m$'th root giving
 an answer whose argument lies in $(-\pi/(2m),\pi/(2m))$; the
 injectivity of $j$ follows easily.) 

 Now define $\dl\:A/B\tm L_\al\to L_0$ by $\dl(q(a),z)=z/\al'(a)$.
 We then have a diagram as shown on the left below:
 \[ \xymatrix{
     A/B\tm L_\al \rto^j \dto_\dl &
     L_\al \dto^\tht &
     S^{L_\al} \rto^{j^!} \dto_\tht &
     (A/B)_+\Smash S^{L_\al} \rto^{\ep\Smash 1} \dto^\dl &
     S^{L_\al} \dto^{\zt_\al} \\
     L_0 \rto_i &
     L_0 &
     S^{L_0} \rto_{i^!\simeq 1} &
     S^{L_0} \rto_{\zt_0} &
     P\CU^T.
    }
 \]
 One can check directly that the left hand square is a pullback, the
 horizontal maps are open inclusions, and the vertical maps are
 proper.  We can thus apply the Pontrjagin-Thom construction to the
 horizontal maps, and one-point compactification to the vertical maps,
 to give a commutative square, which is the left-hand half of the
 right-hand diagram.  Now let $\ep\:A/B\to 1$ be the projection, and
 let $\zt_0$ and $\zt_\al$ be the evident inclusions, to give the
 remaining square.  We claim that this commutes up to equivariant
 homotopy.  To see this, note that $\CU^B$ is a complete universe for
 the group $A/B$, so $S(\CU^B)$ is connected, and it contains $S(L_0)$
 and $S(L_\al)$.  Choose a path $\lm\:[0,1]\to S(\CU^B)$ with
 $\lm(0)\in S(L_0)$ and $\lm(1)\in S(L_\al)$.  Define $A$-equivariant
 maps $\pi_t\:(A/B)_+\Smash S^{L_\al}\to P\CU^T$ by
 \[ \pi_t(q(a)\Smash z)=(\C.a.\lm(t),z\al'(a)^{-1}a.\lm(t))
      = a.(\C.\lm(t),z\al'(a)^{-1}\lm(t)).
 \]
 We find that $\ep_0=\zt_0\dl$ and $\ep_1=\zt_\al(\ep\Smash 1)$, so
 these maps are homotopic as claimed.  

 We now chase the Thom class $\ub_T\in\tE^0(P\CU^T)$ around the
 diagram.  We have $(\zt_0i^!)^*\ub_T=\zt_0^*\ub_T=\ub_{L_0}$ so 
 \[ (\zt_0i^!\tht)^*(\ub_T)=\tht^*(\ub_{L_0})=
      v_m(\phi(\al))\ub_{L_\al}.
 \]
 On the other hand, we have $\zt_\al^*(\ub_T)=\ub_{L_\al}$.  Moreover,
 $j^!$ is $1_{S^{L_\al}}$ smashed with the transfer map
 $S^0\to (A/B)_+$, so
 $((\ep\Smash 1)j^!)^*(\ub_{L_{\al}})=\trf_B^A(1)\ub_{L_\al}$.  Putting
 this together we find that 
 $\trf_B^A(1)\ub_{L_\al}=v_m(\phi(\al))\ub_{L_\al}$, but
 $\tE^0(S^{L_\al})$ is freely generated over $E^0$ by $\ub_{L_\al}$ so
 $\trf_B^A(1)=v_m(\phi(\al))$ as claimed.
\end{proof}

\begin{proof}[Proof of Proposition~\ref{prop-eta-general}]
 Write $A/B$ as a product of cyclic groups.  Each cyclic factor can
 itself be regarded as a quotient of $A$, so we get a decomposition
 $A/B=\prod_iA/B_i$ say.  Now $[A/B]=\prod_i[A/B_i]$ in $\Om(A)$, and
 both $\eta$ and $\eta'$ are ring maps, so it will suffice to prove
 that $\eta[A/B_i]=\eta'[A/B_i]$.  This is clear from
 Lemma~\ref{lem-eta-special}. 
\end{proof}

\subsection{Mackey structure}

We now want to take this further, by constructing a Mackey functor.
\begin{construction}
 Put $k=\OS$ for convenience.  Given a subgroup $B\leq A$, put
 $k_B=k/I(\ann(B),\phi)$, so that the subscheme $S[B]=\spec(k_B)$ is
 the largest closed subscheme of $S$ over which $\phi\:A^*\tm S\to C$
 factors through $B^*$.  We will also write $C[B]=C\tm_SS[B]$, so we
 have an induced map $\phi_B\:B^*\tm S[B]\to C[B]$, and thus a ring
 map $\eta_B\:\Om(B)\to k_B$.  

 If $C\leq B$ then $\ann(C)\geq\ann(B)$ so
 $I(\ann(C),\phi)\geq I(\ann(B),\phi)$, so we have a quotient map
 $k_B\to k_C$.  We write $\res^B_C$ for this map, and note that the
 kernel is $I(\ann(C)/\ann(B),\phi_B)$.

 Next, note that the element
 $\tau_C^B=\eta_B([B/C])=t(\ann(C)/\ann(B),\phi_B)$ annihilates the
 ideal $I(\ann(C)/\ann(B),\phi_B)$, so multiplication by $\tau_C^B$
 induces a well-defined map $\trf_C^B\:k_C\to k_B$.

 As the third ingredient for a Mackey functor we need to specify
 conjugation maps $\gm_a\:k_B\to k_B$ for all $a\in A$.  (Here the
 target is morally $k_{B'}$, where $B'$ is conjugate to $B$ by $a$,
 but of course $B'=B$ because $A$ is abelian.)  We take all the maps
 $\gm_a$ to be the identity.
\end{construction}

\begin{proposition}\label{prop-mackey}
 The above definitions make the system $\{k_B\}_{B\leq A}$ into Mackey
 functor.  
\end{proposition}
\begin{proof}
 The axioms to be checked are as follows.
 \begin{itemize}
  \item[(a)] $\res^B_B=1$ and $\res^B_D=\res^C_D\res^B_C$ whenever
   $D\leq C\leq B\leq A$.
  \item[(b)] $\trf^B_B=1$ and $\trf^B_D=\trf^B_C\trf^C_D$ whenever
   $D\leq C\leq B\leq A$.
  \item[(c)] If $C,D\leq B\leq A$ then 
   \[ \res^B_C\trf_D^B(r)
        = |B/(C+D)| \trf^C_{C\cap D}\res^D_{C\cap D}(r).
   \]
 \end{itemize}
 (In~(c) we have silently made the obvious simplifications to the double
 coset formula arising from the commutativity of $A$.  We have also
 not bothered to list the compatibility conditions between the
 conjugation maps $\gm_a$ and the other structure; these hold for
 trivial reasons, given that $\gm_a=1$.)   

 Here axiom~(a) follows directly from the definitions, and axiom~(b)
 is a translation of Theorem~\ref{thm-transfer}(f).  Similarly,
 axiom~(c) can be derived from part~(g) of
 Theorem~\ref{thm-transfer}. 
\end{proof}

\begin{remark}\label{lem-green}
 Recall that a \emph{Green ring} is a Mackey functor $R$ with a ring
 structure on each group $R(B)$ such that the restriction maps
 $\res^B_C\:R(B)\to R(C)$ and the conjugation maps $\gm_a$ are ring
 maps, and the Frobenius reciprocity formula holds: for $C\leq B$ and
 $r\in R(B)$ and $s\in R(C)$ we have
 $\trf_C^B(\res^B_C(r)s)=r\trf_C^B(s)$.  The Mackey functor $\{k_B\}$
 clearly has these properties.
\end{remark}


\begin{bibdiv}
\begin{biblist}
%\bibselect{../../BiBTeX/refs,../../BiBTeX/myrefs}

\bib{at:kt}{book}{
    author={Atiyah, Michael~F.},
     title={$K$-Theory},
    series={Advanced Book Classics},
 publisher={Addison Wesley},
      date={1989},
}

\bib{atse:ekc}{article}{
  author={Atiyah, Michael~F.},
  author={Segal, G.~B.},
  title={Equivariant K--theory and completion},
  date={1969},
  journal={J. Differential Geometry},
  volume={3},
  pages={1\ndash 18},
}

\bib{co:cor}{thesis}{
  author={Cole, Michael},
  title={Complex oriented $RO(G)$-graded equivariant cohomology theories and their formal group laws},
  type={Ph.D. Thesis},
  date={1996},
}

\bib{cogrkr:efg}{article}{
  author={Cole, Michael},
  author={Greenlees, J. P.~C.},
  author={Kriz, I.},
  title={Equivariant formal group laws},
  date={2000},
  issn={0024-6115},
  journal={Proc. London Math. Soc. (3)},
  volume={81},
  number={2},
  pages={355\ndash 386},
  review={\MR {2001i:55006}},
}

\bib{cogrkr:uec}{article}{
  author={Cole, Michael},
  author={Greenlees, J. P.~C.},
  author={Kriz, I.},
  title={The universality of equivariant complex bordism},
  date={2002},
  issn={0025-5874},
  journal={Math. Z.},
  volume={239},
  number={3},
  pages={455\ndash 475},
  review={\MR {1 893 848}},
}

\bib{da:ghf}{article}{
  author={Damon, James},
  title={The Gysin homomorphism for flag bundles},
  date={1973},
  journal={Amer. J. Math.},
  volume={95},
  pages={643\ndash 659},
  review={\MR {50 \#1255}},
}

\bib{el:sfp}{article}{
  author={Elmendorf, A. D.},
  title={Systems of fixed point sets},
  journal={Trans. Amer. Math. Soc.},
  volume={277},
  date={1983},
  number={1},
  pages={275--284},
  issn={0002-9947},
  review={\MR {690052 (84f:57029)}},
}

\bib{gr:aie}{article}{
  author={Greenlees, J. P.~C.},
  title={Augmentation ideals of equivariant cohomology rings},
  date={1998},
  issn={0040-9383},
  journal={Topology},
  volume={37},
  number={6},
  pages={1313\ndash 1323},
  review={\MR {99h:55005}},
}

\bib{gr:evr}{article}{
  author={Greenlees, J. P. C.},
  title={Equivariant versions of real and complex connective $K$-theory},
  journal={Homology, Homotopy Appl.},
  volume={7},
  date={2005},
  number={3},
  pages={63--82 (electronic)},
  issn={1532-0073},
  review={\MR {2205170 (2006k:19014)}},
}

\bib{grma:lct}{article}{
  author={Greenlees, John P.~C.},
  author={May, J.~Peter},
  title={Localization and completion theorems for $MU$-module spectra},
  date={1997},
  journal={Annals of Mathematics},
  volume={146},
  pages={509\ndash 544},
}

\bib{hokura:ggc}{article}{
  author={Hopkins, Michael~J.},
  author={Kuhn, Nicholas~J.},
  author={Ravenel, Douglas~C.},
  title={Generalized group characters and complex oriented cohomology theories},
  date={2000},
  issn={0894-0347},
  journal={J. Amer. Math. Soc.},
  volume={13},
  number={3},
  pages={553\ndash 594 (electronic)},
  review={\MR {2001k:55015}},
}

\bib{kama:ame}{book}{
  author={Katz, N.~M.},
  author={Mazur, B.},
  title={Arithmetic moduli of elliptic curves},
  series={Annals of Mathematics Studies},
  publisher={Princeton University Press},
  date={1985},
  volume={108},
}

\bib{kr:zec}{incollection}{
  author={Kriz, Igor},
  title={The ${\bf {z}}/p$-equivariant complex cobordism ring},
  date={1999},
  booktitle={Homotopy invariant algebraic structures (baltimore, md, 1998)},
  publisher={Amer. Math. Soc.},
  address={Providence, RI},
  pages={217\ndash 223},
  review={\MR {2001b:55010}},
}

\bib{lemast:esh}{book}{
  author={Lewis, L.~Gaunce},
  author={May, J.~Peter},
  author={(with contributions~by Jim E.~McClure), M.~Steinberger},
  title={Equivariant stable homotopy theory},
  series={Lecture Notes in Mathematics},
  publisher={Springer--Verlag},
  address={New York},
  date={1986},
  volume={1213},
}

\bib{ma:crt}{book}{
  author={Matsumura, Hideyuki},
  title={Commutative ring theory},
  series={Cambridge Studies in Advanced Mathematics},
  publisher={Cambridge University Press},
  date={1986},
  volume={8},
}

\bib{mama:eos}{article}{
  author={Mandell, Michael~A.},
  author={May, J.~Peter},
  title={Equivariant orthogonal spectra and $S$-modules},
  date={2002},
  journal={Mem. Amer. Math. Soc.},
  volume={159},
  pages={108},
}

\bib{qu:epc}{article}{
  author={Quillen, Daniel~G.},
  title={Elementary proofs of some results of cobordism theory using Steenrod operations},
  date={1971},
  journal={Advances in Mathematics},
  volume={7},
  pages={29\ndash 56},
}

\bib{qu:fgl}{article}{
  author={Quillen, Daniel~G.},
  title={On the formal group laws of unoriented and complex cobordism},
  date={1969},
  journal={Bulletin of the American Mathematical Society},
  volume={75},
  pages={1293\ndash 1298},
}

\bib{ta:rdc}{article}{
  author={Tate, John},
  title={Residues of differentials on curves},
  date={1968},
  journal={Ann. Sci. \'Ecole Norm. Sup. (4)},
  volume={1},
  pages={149\ndash 159},
  review={\MR {37 \#2756}},
}

\bib{hoanst:esw}{article}{
  author={Ando, M},
  author={Hopkins, M J},
  author={Strickland, Neil P},
  title={Elliptic spectra, the Witten genus and the theorem of the cube},
  date={2001},
  issn={0020-9910},
  journal={Invent. Math.},
  volume={146},
  number={3},
  pages={595\ndash 687},
  review={\MR {1 869 850}},
}

\bib{anst:wpm}{article}{
  author={Ando, Matthew},
  author={Strickland, Neil P},
  title={Weil pairings and Morava $K$-theory},
  date={2000},
  issn={0040-9383},
  journal={Topology},
  volume={40},
  number={1},
  pages={127\ndash 156},
}

\bib{st:kld}{article}{
  author={Strickland, Neil P},
  title={$K(n)$-local duality for finite groups and groupoids},
  date={2000},
  issn={0040-9383},
  journal={Topology},
  volume={39},
  number={4},
  pages={733\ndash 772},
  review={\MR {1 760 427}},
}

\bib{host:mkl}{article}{
  author={Hovey, Mark},
  author={Strickland, Neil P},
  title={Morava $K$-theories and localisation},
  date={1999},
  issn={0065-9266},
  journal={Mem. Amer. Math. Soc.},
  volume={139},
  number={666},
  pages={104},
}

\bib{st:fsfg}{incollection}{
  author={Strickland, Neil P},
  title={Formal schemes and formal groups},
  date={1999},
  booktitle={Homotopy invariant algebraic structures (baltimore, md, 1998)},
  publisher={Amer. Math. Soc.},
  address={Providence, RI},
  pages={263\ndash 352},
  review={\MR {2000j:55011}},
}

\bib{st:ebc}{unpublished}{
  author={Strickland, Neil P},
  title={Equivariant Bousfield classes},
  date={2008},
  note={In preparation},
}


\end{biblist}
\end{bibdiv}
\end{document}